\def\DateTime{22/September/2002, 3:30(JP)}
\def\Version{Version $1.0$}
\def\yes{\if00}
\def\no{\if01}
\def\iftwelvept{\yes}
\def\ifusepdf{\no}
\def\ifpsfont{\no}

\iftwelvept
\documentclass[leqno,12pt]{amsart}
\else
\documentclass[leqno]{amsart}
\fi
\usepackage{amssymb}
\usepackage{amscd}
\usepackage{latexsym}
\usepackage{verbatim}
\usepackage[all]{xy}

\ifusepdf
\usepackage{hyperref}
\else\fi
\ifpsfont
\usepackage[T1]{fontenc}
\usepackage{times}
\else\fi

\iftwelvept
\else\fi

\theoremstyle{plain}
\newtheorem{Theorem}{Theorem}[section]
\newtheorem{Proposition}[Theorem]{Proposition}
\newtheorem{Lemma}[Theorem]{Lemma}
\newtheorem{Corollary}[Theorem]{Corollary}
\newtheorem{Sublemma}[Theorem]{Sublemma}
\newtheorem{Claim}{Claim}[Theorem]

\theoremstyle{definition}

\newtheorem{Remark}[Theorem]{Remark}

\renewcommand{\theTheorem}{\arabic{section}.\arabic{subsection}.\arabic{Theorem}}
\renewcommand{\theClaim}{\arabic{section}.\arabic{subsection}.\arabic{Theorem}.\arabic{Claim}}
\renewcommand{\theequation}{\arabic{section}.\arabic{subsection}.\arabic{Theorem}.\arabic{Claim}}

\def\rom{\textup}
\newcommand{\ZZ}{{\mathbb{Z}}}
\newcommand{\QQ}{{\mathbb{Q}}}
\newcommand{\RR}{{\mathbb{R}}}
\newcommand{\CC}{{\mathbb{C}}}
\newcommand{\PP}{{\mathbb{P}}}

\newcommand{\DD}{{\mathbb{D}}}
\newcommand{\FF}{{\mathbb{F}}}
\newcommand{\OO}{{\mathcal{O}}}
\newcommand{\XX}{{\mathcal{X}}}

\newcommand{\LL}{{\mathcal{L}}}

\newcommand{\Proj}{\operatorname{Proj}}
\newcommand{\aPic}{\widehat{\operatorname{Pic}}}
\newcommand{\aCH}{\widehat{\operatorname{CH}}}
\newcommand{\Pic}{\operatorname{Pic}}
\newcommand{\Spec}{\operatorname{Spec}}
\newcommand{\Div}{\operatorname{Div}}
\newcommand{\Supp}{\operatorname{Supp}}
\newcommand{\acherncl}{\widehat{{c}}}

\newcommand{\adeg}{\widehat{\operatorname{deg}}}

\newcommand{\zeros}{\operatorname{div}}
\newcommand{\ord}{\operatorname{ord}}

\newcommand{\trdeg}{\operatorname{tr.deg}}
\newcommand{\rank}{\operatorname{rk}}
\newcommand{\lc}{\operatorname{lc}}
\newcommand{\crest}{{\,\vert\,}}

\newcommand{\Proof}{{\sl Proof.}\quad}
\newcommand{\QED}{{\unskip\nobreak\hfil\penalty50\quad\null\nobreak\hfil
{$\Box$}\parfillskip0pt\finalhyphendemerits0\par\medskip}}
\newcommand{\rest}[2]{\left.{#1}\right\vert_{{#2}}}

\begin{document}

\title[The number of algebraic cycles with bounded degree]%
{The number of algebraic cycles with bounded degree}
\author{Atsushi Moriwaki}
\address{Department of Mathematics, Faculty of Science,
Kyoto University, Kyoto, 606-8502, Japan}
\email{moriwaki@kusm.kyoto-u.ac.jp}
\date{\DateTime, (\Version)}
\begin{abstract} 
Let $X$ be a projective scheme over a finite field.
In this paper, we consider the asymptotic
behavior of the number of effective cycles on $X$ with bounded degree as it goes
to the infinity. By this estimate, we can define a certain kind of zeta functions
associated with groups of cycles.
We also consider an analogue in Arakelov geometry.
\end{abstract}


\maketitle


\renewcommand{\theTheorem}{\Alph{Theorem}}
\section*{Introduction}
Let $X$ be a projective scheme over a finite field $\FF_q$ and
$H$ an ample line bundle on $X$.
For a non-negative integer $k$, 
we denote by $n_k(X,H,l)$
the number of all effective $l$-dimensional cycles $V$ on $X$ 
with $\deg_{H}(V) = k$,
where $\deg_{H}(V)$ is the degree of $V$ with respect to $H$
given by
\[
 \deg_{H}(V) = \deg\left(H^{\cdot l} \cdot V\right).
\]
One of the main results of this paper is to give an estimate of
$n_k(X,H,l)$ as $k$ goes to the infinity,
namely,

\begin{Theorem}[Geometric version]
\begin{enumerate}
\renewcommand{\labelenumi}{(\arabic{enumi})}
\item
If $H$ is very ample, then
there is a constant $C$ depending only on $l$ and
$\dim_{\FF_q} H^0(X, H)$ such that
$\log_q n_k(X,H,l) \leq C k^{l+1}$
for all $k \geq 0$.
\item
If $l \not= \dim X$, then
${\displaystyle \limsup_{k\to\infty}
\frac{\log_q n_k(X,H,l)}{k^{l+1}} > 0}$.
\end{enumerate}
\end{Theorem}
As a consequence of the above theorem, we can define
a certain kind of zeta functions of algebraic cycles as follows.
We note that Weil's
zeta function $Z_{X/\FF_q}$ is given by
\[
\sum_{k=0}^{\infty} n_k(X,H,0) T^{k};
\]
accordingly we define a zeta function $Z(X,H,l)$ 
of $l$-dimensional cycles on
a polarized scheme $(X,H)$ over $\FF_q$ to be
\[
Z(X,H,l)(T) = \sum_{k=0}^{\infty} n_k(X,H,l) T^{k^{l+1}}.
\]
Then, by the above theorem, we can see that
$Z(X,H,l)(T)$ is a convergent power series at the origin.

Further, using the same techniques, we can estimate the number of rational points
defined over a function field.
Let $C$ be a projective smooth curve over $\FF_q$ and $F$ the function field of $C$.
Let $f : X \to C$ be a morphism of projective varieties over $\FF_q$ and
$L$ an $f$-ample line bundle on $X$. Let $X_{\eta}$ be the generic fiber of $f$.
For $x \in X_{\eta}(\overline{F})$, we define the height of $x$ with respect to $L$
to be
\[
h_L(x) = \frac{(L \cdot \Delta_x)}{\deg(\Delta_x \to C)},
\]
where $\Delta_x$ is the Zariski closure of the image
of $\Spec(\overline{F}) \to X_{\eta} \hookrightarrow X$.
Then, we can see that, for a fixed $k$, there is a constant $C$ such that
\[
\# \{ x \in X(\overline{F}) \mid \text{$[F(x) : F] \leq k$ and $h_L(x) \leq h$} \}
\leq q^{C \cdot h}
\]
for all $h \geq 1$. Thus, a series
\[
\sum_{\substack{x \in X_{\eta}(\overline{F}), \\ [F(x) : F] \leq k}} q^{-sh_L(x)}
\]
converges for all $s \in \CC$ with $\Re(s) \gg 0$.
This is a local analogue of Batyrev-Manin-Tschinkel's height zeta functions.

Moreover, let $\mathcal{X} \to \Spec(O_K)$ be a flat and projective scheme over
the ring $O_K$ of integers of a number field $K$ and let $\mathcal{H}$ be an ample
line bundle on $\mathcal{X}$. Then, as a corollary of our estimates, we can see
an infinite product
\[
L(\mathcal{X}, \mathcal{H}, l)(s) =
\prod_{P \in \Spec(O_K) \setminus \{ 0 \}}
Z(\mathcal{X}_P,\mathcal{H}_P,l)(\#\kappa(P)^{-s})
\]
converges for all $s \in \CC$ with $\Re(s) \gg 0$, which looks like
a generalization of the usual $L$-functions.

\bigskip
The next purpose of this paper is to give an analogue in
Arakelov geometry.
Let $X$ be a projective arithmetic variety, i.e,
a flat and projective integral scheme over $\ZZ$.
Let $\overline{H}$
be an ample $C^{\infty}$-hermitian $\QQ$-line bundle on $X$.
For a cycle $V$ of dimension $l$ on $X$,
the arithmetic degree of $V$ is defined by
\[
\adeg_{\overline{H}}(V) = \adeg\left(
\acherncl_1(\overline{H})^{\cdot l}\crest V \right).
\]
For a real number $h$, we denote by $\hat{n}_{\leq h}(X, \overline{H}, l)$
(resp. $\hat{n}_{\leq h}^{\rm hol}(X, \overline{H}, l)$)
the number of effective cycles (resp. horizontal effective cycles) 
$V$ of dimension $l$ on $X$ 
with $\adeg_{\overline{H}}(V) \leq h$.
Then, we have the following analogue.

\begin{Theorem}[Arithmetic version]
\label{thm:intro:Arith:version}
\begin{enumerate}
\renewcommand{\labelenumi}{(\arabic{enumi})}
\item
There is a constant $C$ such that
\[
\log \hat{n}_{\leq h}(X, \overline{H}, l) \leq Ch^{l+1}
\]
for all $ h \geq 0$.

\item
If $l \not= \dim X$, then
${\displaystyle \limsup_{h\to\infty}
\frac{\log \hat{n}^{\rm hol}_{\leq h}(X, \overline{H}, l)}{h^{l+1}}
> 0}$.
\end{enumerate}
\end{Theorem}
Techniques involving the proof of Theorem~\ref{thm:intro:Arith:version} are
much harder than the geometric case, but the outline for the proof is similar to geometric one.
We have also the estimate of rational points defined over a finitely generated field
over $\QQ$ (cf. Theorem~\ref{thm:relative:case:b:geom:degree}).

\bigskip
Finally, we would like to give hearty thanks to Prof. Mori,
Prof. Soul\'{e} and Prof. Wan for their useful comments and suggestions for this paper.

\renewcommand{\theTheorem}{\arabic{section}.\arabic{subsection}.\arabic{Theorem}}

\section{Notations and Conventions}
\renewcommand{\thesubsection}{(\arabic{section}.\arabic{subsection})}

Here, 
we introduce notations and
conventions in this paper.

\subsection{}
For a point $x$ of a scheme $X$,
the residue field at $x$ is denoted by $\kappa(x)$.

\subsection{}
\label{subsec:notation:cycle}
Let $X$ be a Noetherian scheme. For a non-negative integer $l$,
we denote by $\mathcal{C}_l(X)$ the set of all
$l$-dimensional integral closed subschemes on $X$.
We set
\[
Z_l(X) = \bigoplus_{V \in \mathcal{C}_l(X)} \ZZ V,
\quad\text{and}\quad
Z_l^{\rm eff}(X) = \bigoplus_{V \in \mathcal{C}_l(X)} \ZZ_{\geq 0} V,
\]
where $\ZZ_{\geq 0} = \{ z \in \ZZ \mid z \geq 0 \}$.
An element of $Z_l(X)$ (resp. $Z_l^{\rm eff}(X)$)
is called an {\em $l$-dimensional cycle} (resp.
{\em $l$-dimensional effective cycle}) on $X$.

For a subset $\mathcal{C}$ of $\mathcal{C}_l(X)$,
we denote $\bigoplus_{V \in \mathcal{C}} \ZZ V$ and
$\bigoplus_{V \in \mathcal{C}} \ZZ_{\geq 0} V$ by
$Z_l(X;\mathcal{C})$ and $Z_l^{\rm eff}(X;\mathcal{C})$ respectively.
In this paper, we consider the following
$\mathcal{C}(U)$ and $\mathcal{C}(X/Y)$ as a subset of
$\mathcal{C}_l(X)$;
For a Zariski open set $U$ of $X$, we set
\[
\mathcal{C}(U) = 
\{ V \in \mathcal{C}_l(X) \mid V \cap U \not= \emptyset \}.
\]
For a morphism $f : X \to Y$ of Noetherian schemes
with $Y$ irreducible, we set
\[
\mathcal{C}(X/Y) = \{ V \in \mathcal{C}_l(X) \mid \overline{f(V)} = Y \}.
\]
For simplicity,
we denote $Z_l(X;\mathcal{C}(U))$,
$Z_l^{\rm eff}(X;\mathcal{C}(U))$,
$Z_l(X;\mathcal{C}(X/Y))$ and
$Z_l^{\rm eff}(X;\mathcal{C}(X/Y))$ by
$Z_l(X;U)$,
$Z_l^{\rm eff}(X;U)$,
$Z_l(X/Y)$ and
$Z_l^{\rm eff}(X/Y)$ respectively.
In order to show the fixed morphism $f : X \to Y$,
$Z_l(X/Y)$ and $Z^{\rm eff}_l(X/Y)$ are sometimes denoted by
$Z_l(X \overset{f}{\to} Y)$ and $Z_l^{\rm eff}(X \overset{f}{\to} Y)$
respectively.

\subsection{}
Let $R$ be a commutative ring with the unity. 
Let $X$ be the product of projective spaces $\PP^{n_1}_{R},
\ldots, \PP^{n_r}_{R}$ over $R$, that is,
\[
X = \PP^{n_1}_{R} \times_R \cdots \times_R \PP^{n_r}_{R}.
\]
If $R$ is UFD, then, for a divisor $D$ on $X$,
there is the unique sequence $(k_1, \ldots, k_n)$ of
non-negative integers and the unique section
$s \in H^0\left(X, \bigotimes_{i=1}^n p_i^*(\OO(k_i))\right)$
module $R^{\times}$ such that
$\zeros(s) = D$, where
$p_i : X \to \PP^{n_i}_{R}$ is
the projection to the $i$-th factor. 
We denote $k_i$ by $\deg_i(D)$ and
call it the {\em $i$-th degree of $D$}.
Moreover, for simplicity,
we denote
\[
\underbrace{\PP^n_{R} \times_R \cdots \times_R \PP^n_{R}}_{\text{$r$-times}}
\]
by $(\PP^n_R)^r$.
Note that $(\PP^n_{R})^0 = \Spec(R)$.

\subsection{}
\label{subsec:product:P:1:projection}
For a non-negative integer $n$, we set
\[
[n] = \begin{cases}
\{ 1, 2, \ldots, n \} & \text{if $n \geq 1$} \\
\emptyset & \text{if $n=0$}.
\end{cases}
\]
We assume $n \geq 1$. Let us consider the scheme $(\PP^1_{R})^n$ over $R$, where $R$ is a commutative ring.
Let $p_i : (\PP^1_{R})^n \to \PP^1_R$ be the projection to the $i$-th factor.
For a subset $I$ of $[n]$, we define $p_I : (\PP^1_{R})^n \to (\PP^1_{R})^{\#( I )}$ as follows:
If $I = \emptyset$, then $p_I$ is the canonical morphism $(\PP^1_{R})^n \to \Spec(R)$.
Otherwise, we set $I = \{ i_1, \ldots, i_{\#( I )} \}$ with $1 \leq i_1 < \cdots < i_{\#( I )} \leq n$.
Then, $p_I = p_{i_1} \times \cdots \times p_{i_{\#( I )}}$, i.e.,
$p_I(x_1, \ldots, x_n) = (x_{i_1}, \ldots, x_{i_{\#( I )}})$.
Note that $p_{\{ i \}} = p_i$.

\subsection{}
Let us fix a basis $\{X_0, \ldots, X_n\}$ of $H^0(\PP^n_{\CC}, \OO(1))$.
The Fubini-Study metric $\Vert\cdot\Vert_{\rm FS}$ 
of $\OO(1)$ with respect to $\{ X_0, \ldots, X_n \}$ is 
given by
\[
\Vert X_i \Vert_{\rm FS} = 
\frac{\vert X_i \vert}{\sqrt{\vert X_0 \vert^2 + \cdots +  \vert X_n \vert^2}}.
\]
For a real number $\lambda$,
the metric $\exp(-\lambda)\Vert \cdot \Vert_{\rm FS}$ is denoted by
$\Vert \cdot \Vert_{{\rm FS}_{\lambda}}$.
Moreover, the hermitian line bundle
$(\OO(1), \Vert\cdot\Vert_{{\rm FS}_{\lambda}})$ is denoted by
$\overline{\OO}^{{\rm FS}_{\lambda}}(1)$.

\subsection{}
\label{subsubsec:asymp:notation}
Let $f$ and $g$ be real valued functions on a set $S$.
We use the notation `$f \asymp g$'
if there are positive real numbers $a, a'$ and real numbers $b,b'$
such that $g(s) \leq af(s) + b$ and $f(s) \leq a'g(s) + b'$
for all $s \in S$.

\renewcommand{\thesubsection}{\arabic{section}.\arabic{subsection}}

\section{The geometric case}
The main purpose of this section is to find
a universal upper bound of the number of
effective cycles with bounded degree on the projective space
over a finite field (cf. Theorem~\ref{thm:proj:sp:geom:case}), 
namely,
\begin{quote}
Fix non-negative integers $n$ and $l$.
Then, there is a constant $C(n,l)$ depending only on
$n$ and $l$ such that the number of effective $l$-dimensional
cycles on $\PP^n_{\FF_q}$ with degree $k$ is less than or equal to
$q^{C(n,k) k^{l+1}}$. 
\end{quote}
First we consider a similar problem on the product
$(\PP^1_{\FF_q})^n$ of the projective line.
The advantage of $(\PP^1_{\FF_q})^n$ is that
it has a lot of morphisms, so that
induction on its dimension works well.

\subsection{Preliminaries}
\setcounter{Theorem}{0}
Here let us prepare basic tools to count cycles.

Let $\{ T_n \}_{n=n_0}^{\infty} = \{T_{n_0}, T_{n_0 + 1}, \ldots, T_n, \ldots \}$ 
be a sequence of sets.
If it satisfies the following properties (1) -- (4),
then it is called a {\em counting system}.
\begin{enumerate}
\renewcommand{\labelenumi}{(\arabic{enumi})}
\item
For each $n \geq n_0$, there is a function $h_n : T_n \to \RR_{\geq 0}$.

\item
For each $n \geq n_0 + 1$, there are maps $\alpha_n : T_n \to T_{n-1}$ and
$\beta_n : T_n \to T_{n_0}$ such that
\[
h_{n-1}(\alpha_n(x)) \leq h_n(x)\quad\text{and}\quad
h_{n_0}(\beta_n(x)) \leq h_n(x)
\]
for all $x \in T_n$.

\item
There is a function $A : \RR_{\geq 0} \times \RR_{\geq 0} \to \RR$
such that
$A(s, t) \leq A(s', t')$ for all $0 \leq s \leq s'$ and $0 \leq t \leq t'$ and that,
for $y \in T_{n-1}$ and $z \in T_{n_0}$,
\[
\#\{ x \in T_n \mid \text{$\alpha_n(x) = y$ and $\beta_n(x) = z$} \}
\leq A(h_{n-1}(y), h_{n_0}(z)).
\]

\item
There is a function $B: \RR_{\geq 0} \to \RR$ and a non-negative constant $t_0$ such that
\[
\# \{ x \in T_{n_0} \mid h_{n_0}(x) \leq h \} \leq B(h)
\]
for all $h \geq t_0$.
\end{enumerate}

\begin{Lemma}
\label{lem:counting:system}
If $\{ T_n \}_{n=n_0}^{\infty}$ is a counting system as above, then
\[
\# \{ x \in T_n \mid h_n(x) \leq h \} \leq B(h)^{n-n_0+1}A(h, h)^{n-n_0}
\]
for all $h \geq t_0$.
\end{Lemma}

\Proof
For $x \in T_n$ with $h_n(x) \leq h$,
by the property (2), we have $h_{n-1}(\alpha_n(x)) \leq h$ and
$h_{n_0}(\beta_n(x)) \leq h$.
Thus, by using (3) and (4),
\begin{align*}
\# \{ x \in T_n \mid h_n(x) \leq h \} & \leq
\# \{ y \in T_{n-1} \mid h_{n-1}(y) \leq h \} \cdot
\# \{ z \in T_{n_0} \mid h_{n_0}(z) \leq h \} \cdot A(h, h) \\
 & \leq
\# \{ y \in T_{n-1} \mid h_{n-1}(y) \leq h \} \cdot
B(h) \cdot A(h, h).
\end{align*}
Therefore, we get our lemma by using induction on $n$.
\QED

The following lemma will be used to see the above property (3).

\begin{Lemma}
\label{lem:number:0:cycle:product}
Let $X$ and $Y$ be projective schemes
over a field $K$. Let $p : X \times_K Y \to X$ and $q: X \times_K Y \to Y$ be
the projection to the first factor and
the projection to the second factor respectively.
Let $x_1, \ldots, x_s$ {\rm (}resp. $y_1, \ldots, y_t${\rm )}
be closed points of $X$ {\rm (}resp. $Y${\rm )}.
Let us fix an effective $0$-cycle $x = \sum_{i=1}^s a_i x_i$ and
an effective $0$-cycle $y = \sum_{j=1}^t b_j y_j$.
Then, the number of effective $0$-cycles $z$ on $X \times_K Y$ with
$p_*(z) = x$ and $q_*(z) = y$ is less than or equal to
$2^{\alpha_X(x)\alpha_Y(y)}$,
where $\alpha_X(x) = 
\sum_{i=1}^s \sqrt{a_i[\kappa(x_i):K]} $
and $\alpha_Y(y) =
\sum_{j=1}^t \sqrt{b_j[\kappa(y_j):K]}$.
\end{Lemma}

\Proof
Let $z_{ijk}$'s ($k=1, \ldots, l_{ij}$) be all closed points
of $\Spec(\kappa(x_i) \otimes_K \kappa(y_j))$.
Then, an effective $0$-cycle $z$ on $X \times_K Y$ with
$p_*(z) = x$ and $q_*(z) = y$ can be written by the form
$\sum_{ijk} c_{ijk}z_{ijk}$.
Hence,
\[
p_*(z) = \sum_{i} \left(\sum_{j,k}
[\kappa(z_{ijk}):\kappa(x_i)]c_{ijk}\right)x_i
\]
and
\[
q_*(z) = \sum_{j} \left(\sum_{i,k} 
[\kappa(z_{ijk}):\kappa(y_j)]c_{ijk}\right)y_j.
\]
Thus,
\[
c_{ijk} \leq \min \left\{a_i,  b_j \right\} \leq
\sqrt{a_ib_j}.
\]
Therefore,
the number $N(x,y)$ of effective $0$-cycles $z$ on $X \times_K Y$ with
$p_*(z) = x$ and $q_*(z) = y$ is less than or equal to
$\prod_{ij}(1 + \sqrt{a_ib_j})^{l_{ij}}$.
Here note that
\[
l_{ij} \leq
\min\left\{ [\kappa(x_i) : K], [\kappa(x_j):K] \right\} \leq
\sqrt{[\kappa(x_i) : K][\kappa(x_j):K]}.
\]
Moreover, $1+x \leq 2^x$ for $x \in \{ 0 \}
\cup [1,
\infty)$. Hence,
\begin{align*}
N(x,y) & \leq \prod_{ij}(1 + \sqrt{a_i\cdot b_j}
)^{\sqrt{[\kappa(x_i) : K][\kappa(y_j): K]}} \\
& \leq \prod_{ij} 2^{\sqrt{a_i[\kappa(x_i) : K]}\sqrt{b_j[\kappa(y_j): K]}} =
2^{\sum_{ij}\sqrt{a_i[\kappa(x_i) : K]}\sqrt{b_j[\kappa(y_j): K]}}.
\end{align*}
Thus, we get our lemma.
\QED

The following lemma will be also used to count cycles.

\begin{Lemma}
\label{lem:estimate:number:push:forward}
Let $\pi : X' \to X$ be a finite morphism of normal integral schemes.
Let $Z = \sum_{i=1}^n a_i Z_i$ be an effective cycle on $X$,
where $Z_i$'s are integral.
Then the number of effective cycles $Z'$ on $X'$ with $\pi_*(Z') = Z$
is less than or equal to
$2^{\deg(\pi) \sum_{i=1}^n a_i}$.
\end{Lemma}

\Proof
We denote by $\alpha(Z)$ 
the number of effective cycles $Z'$ on $X'$ with $\pi_*(Z') = Z$.
Let $Z'_{i1}, \ldots, Z'_{it_i}$ be all integral subschemes lying over $Z_i$.
Then, $t_i \leq \deg(\pi)$.
Let $Z'$ be an effective cycle $Z'$ on $X'$ with $\pi_*(Z') = Z$.
Then, we can set
$Z' = \sum_{i=1}^{n} \sum_{j=1}^{t_i} a_{ij} Z_{ij}$.
Since $\pi_*(Z') = Z$, 
the number of possible $(a_{i1}, \ldots, a_{it_i})$'s is at most
$(1 + a_i)^{\deg(\pi)}$.
Therefore,
\[
\alpha(D) \leq \prod_{i=1}^n (1 + a_i)^{\deg(\pi)}.
\]
Here note that $1+x \leq 2^x$ for $x  \in \{ 0 \} \cup [1, \infty)$. Hence
we get our lemma.
\QED

\subsection{Cycles on $(\PP^1_{\FF_q})^n$}
\setcounter{Theorem}{0}

Let us begin with the case of divisors.

\begin{Proposition}
\label{prop:estimate:X:product:P:1:divisor}
Let $k_1, \ldots, k_n$ be non-negative integers.
Then
\[
\# \{ D \in \Div^{\rm eff}((\PP^1_{\FF_q})^{n}) \mid 
\deg_i(D) \leq k_i \ \forall i = 1, \ldots, n \} \leq 
\frac{q^{\prod_{i=1}^n (k_i + 1)}-1}{q-1} \cdot \prod_{i=1}^n (k_i + 1).
\]
\end{Proposition}

\Proof
In the following, the symbol $D$ is an effective divisor on 
$(\PP^1_{\FF_q})^{n}$.
\begin{align*}
\# \left\{ D \mid 
\deg_i(D) \leq k_i \ (\forall i) \right\}
& = \sum_{0 \leq e_1 \leq k_1, \ldots, 0 \leq e_n \leq k_n}
 \# \left\{ D \mid 
\text{$\deg_i(D) = e_i$ ($\forall i$)} \right\} \\
& = \sum_{0 \leq e_1 \leq k_1, \ldots, 0 \leq e_n \leq k_n}
\frac{q^{(e_1+1)\cdots(e_{n} + 1)} -1}{q-1} \\
& \leq (k_1 + 1) \cdots (k_n+1) \frac{q^{(k_1 + 1) \cdots (k_n+1)}-1}{q-1}
\end{align*}
\QED

In order to proceed with induction,
the following lemmas are very useful.

\begin{Lemma}
\label{lem:number:product:X:Y}
Let $f : X \to S$ and $g : Y \to S$ be morphisms of
projective schemes over $\FF_q$. We assume that $S$ is
integral and of dimension $l$.
Let $p : X \times_S Y \to X$ and $q : X \times_S Y \to Y$ be
the projections to the first factor and the second factor respectively.
Fix $D \in Z_l^{\rm eff}(X/S)$ and $E \in Z^{\rm eff}_l(Y/S)$
\rom{(}for the definition of $Z_l^{\rm eff}(X/S)$ and $Z_l^{\rm eff}(Y/S)$,
see \rom{\ref{subsec:notation:cycle})}.
\begin{enumerate}
\renewcommand{\labelenumi}{(\arabic{enumi})}
\item
Assume $l \geq 1$.
Let $A_1, \ldots, A_l$ be nef line bundles on $X$, 
$B_1, \ldots, B_l$ nef line bundle on $Y$, and $C_1, \ldots, C_l$ 
nef line bundles on $S$
such that $A_i \otimes f^*(C_i)^{\otimes -1}$ and
$B_i \otimes g^*(C_i)^{\otimes -1}$ are nef for all $i$ and that $\deg(C_1 \cdots C_l) > 0$.
Then, 
\begin{multline*}
\log_q \left( \# \left\{ V \in Z^{\rm eff}_l(X \times_S Y/S) \mid
\text{$p_*(V) = D$ and $q_*(V) = E$} \right\} \right) \\
\leq \min \left\{
\frac{\deg(A_1 \cdots A_l \cdot D)\deg(B_1\cdots B_l \cdot E)}{\deg\left(C_1 \cdots C_l \right)^2},
\right. \\
\left. \frac{\sqrt{\theta(D)\theta(E)\deg(A_1 \cdots A_l \cdot D)
\deg(B_1 \cdots B_l \cdot E)}}{\deg\left(C_1 \cdots C_l \right)}
\right\},
\end{multline*}
where $\theta(D)$ \rom{(}resp. $\theta(E)$\rom{)} is 
the number of irreducible components of $\Supp(D)$ \rom{(}resp. $\Supp(E)$\rom{)}.

\item
Assume $l=0$, so that $S = \Spec(\FF_{q^r})$ for some positive integer $r$.
Then,
\begin{multline*}
\log_q \left( \# \left\{ V \in Z^{\rm eff}_0(X \times_S Y/S) \mid
\text{$p_*(V) = D$ and $q_*(V) = E$} \right\} \right) \\
\leq \min \left\{
\frac{\deg(D)\deg(E)}{r^2}, \frac{\sqrt{\theta(D)\theta(E)\deg(D)
\deg(E)}}{r}
\right\}.
\end{multline*}
\end{enumerate}
\end{Lemma}

\Proof
(1)
We set $D = \sum_{i=1}^s a_i D_i$ and $E = \sum_{j=1}^t b_j E_j$.
Then,
\addtocounter{Claim}{1}
\begin{multline}
\label{eqn:lem:number:product:X:Y:1}
\deg(A_1 \cdots A_l \cdot D) = \sum_{i} a_i \deg(A_1 \cdots A_l \cdot D_i) \\
\geq  \sum_{i=1}^s a_i \deg(f^*(C_1) \cdots f^*(C_l) \cdot D_i)
= \sum_{i=1}^s a_i \deg(D_i \to S) \deg(C_1 \cdots C_l).
\end{multline}
In the same way,
\addtocounter{Claim}{1}
\begin{equation}
\label{eqn:lem:number:product:X:Y:2}
\deg(B_1 \cdots B_l \cdot E) \geq \sum_{j=1}^t b_j  \deg(E_j \to S) \deg(C_1 \cdots C_l).
\end{equation}
Thus,
\begin{multline*}
\frac{\deg(A_1 \cdots A_l \cdot D)}{\deg(C_1 \cdots C_l)} \geq
\sum_{i=1}^s \sqrt{a_i \deg(D_i \to S)}
\quad\text{and}\quad
\frac{\deg(B_1 \cdots B_l \cdot E)}{\deg(C_1 \cdots C_l)} \geq
\sum_{j=1}^t \sqrt{b_j \deg(E_j \to S)}.
\end{multline*}
Moreover, note that
\[
\sqrt{n} \sqrt{x_1 + \cdots + x_n} \geq \sqrt{x_1} + \cdots + \sqrt{x_n}.
\]
Thus, the above inequalities \eqref{eqn:lem:number:product:X:Y:1} and
\eqref{eqn:lem:number:product:X:Y:2} imply
\[
\sqrt{\frac{s\deg(A_1 \cdots A_l \cdot D)}{\deg(C_1 \cdots C_l)}} \geq
\sum_{i=1}^s \sqrt{a_i \deg(D_i \to S)}
\]
and
\[
\sqrt{\frac{t\deg(B_1 \cdots B_l \cdot E)}{\deg(C_1 \cdots C_l)}} \geq
\sum_{j=1}^t \sqrt{b_j \deg(E_j \to S)}.
\]
Therefore, considering
$X$, $Y$ and $X \times_S Y$ over the generic point of $S$,
Lemma~\ref{lem:number:0:cycle:product} implies our assertion.

\medskip
(2)
We set $D = \sum_{i=1}^s a_i x_i$ and $E = \sum_{j=1}^t b_j y_j$.
Then,
\[
\deg(D) = \sum_{i=1}^s a_i [\kappa(x_i) : \FF_{q}] = r \sum_{i=1}^s a_i [\kappa(x_i) : \FF_{q^r}]
\]
and
\[
\deg(E) = \sum_{j=1}^t b_j [\kappa(y_j) : \FF_{q}] = r \sum_{j=1}^t b_j [\kappa(y_j) : \FF_{q^r}].
\]
Thus, in the same way as in (1), we get our assertion.
\QED

Using the above lemma, we have the following.

\begin{Proposition}
\label{prop:growth:P:1:product:l}
There is a constant $C(n,l)$ depending only $n$ and $l$
such that
\[
\# \{ V \in Z^{\rm eff}_{l}( (\PP^1_{\FF_q})^n) \mid \deg_{\OO_{(\PP^1)^n}(1, \ldots, 1)}(V) \leq h \} \leq 
q^{C(n,l)h^{l+1}}
\]
for all $h \geq 1$.
\end{Proposition}

\Proof
First we assume $l=0$.
Let us see
\[
\# \{ V \in Z^{\rm eff}_{0}( (\PP^1_{\FF_q})^n) \mid \deg(V) \leq h \} \leq q^{3nh}
\]
for $h \geq 1$. We prove this by induction on $n$.
If $n=1$, then our assertion follows from Proposition~\ref{prop:estimate:X:product:P:1:divisor},
so that we assume $n > 1$.
Let $q : (\PP^1_{\FF_q})^n \to (\PP^1_{\FF_q})^{n-1}$ be the projection given by
$q(x_1, \ldots, x_n) = (x_1, \ldots, x_{n-1})$.
For a fixed $W \in Z^{\rm eff}_{0}( (\PP^1_{\FF_q})^{n-1})$,
let us estimate the number of
$\{ V \in Z^{\rm eff}_{0}( (\PP^1_{\FF_q})^n) \mid q_*(V) = W \}$.
We set $W = \sum_{i=1}^e a_i y_i$.
For $V \in Z^{\rm eff}_{0}( (\PP^1_{\FF_q})^n)$ with $q_*(V) = W$,
let $V = V_1 + \cdots + V_e$ be the decomposition of effective $0$-cycles
with $q_*(V_i) = a_iy_i$ ($i=1, \ldots, e$).
Then, $V_i \in  Z^{\rm eff}_{0}(\PP^1_{\kappa(y_i)})$ and
the degree of $V_i$ in $\PP^1_{\kappa(y_i)}$ is $a_i$.
Thus, the possible number of $V_i$ is less than or equal to
$\#(\kappa(y_i))^{3a_i}$. Thus,
\[
\# \{ V \in Z^{\rm eff}_{0}( (\PP^1_{\FF_q})^n) \mid q_*(V) = W \} \leq
\prod_{i=1}^e \#(\kappa(y_i))^{3a_i} = \prod_{i=1}^e q^{3[\kappa(y_i) : \FF_q]a_i} = q^{3\deg(W)}.
\]
Therefore, since $\deg(V) = \deg(q_*(V))$, using the hypothesis of induction,
\begin{align*}
\# \{ V \in Z^{\rm eff}_{0}( (\PP^1_{\FF_q})^n) \mid \deg(V) \leq h \} & \leq
\# \{ W \in Z^{\rm eff}_{0}( (\PP^1_{\FF_q})^{n-1}) \mid \deg(W) \leq h \} \cdot q^{3h} \\
& \leq q^{3(n-1)h} \cdot q^{3h} = q^{3nh}.
\end{align*}

\medskip
Next we assume $l \geq 1$.
For a subset $I$ of $[n] = \{ 1, \ldots, n \}$ with $\#( I ) = l$, let us consider
the morphism $p_I : (\PP^1_{\FF_q})^n \to (\PP^1_{\FF_q})^l$ (for the definition of $p_I$,
see \ref{subsec:product:P:1:projection}).
We denote by $Z^{{\rm eff}}_{l}
( (\PP^1_{\FF_q})^n \overset{p_{I}}{\to} (\PP^1_{\FF_q})^l)$ the set of
effective cycles on $ (\PP^1_{\FF_q})^n$ generated by
$l$-dimensional subvarieties which dominates $(\PP^1_{\FF_q})^l$ via 
$p_{I}$.
Then, it is easy to see that
\[
Z^{{\rm eff}}_{l}( (\PP^1_{\FF_q})^n) = 
\sum_{I \subseteq [n], \#( I ) = l} Z^{{\rm eff}}_{l}(
(\PP^1_{\FF_q})^n \overset{p_{I}}{\to} (\PP^1_{\FF_q})^l).
\]
Thus, since
\[
\#(\{ I \mid I \subseteq [n], \#( I ) = l \}) = \binom{n}{l} \leq 2^n,
\]
it is sufficient to see that
there is a constant $C'(n,l)$ depending only on
$n$ and $l$ such that
\[
\{ V \in Z^{{\rm eff}}_{l}(
(\PP^1_{\FF_q})^n \overset{p_I}{\to} (\PP^1_{\FF_q})^l) \mid
\deg_H(V) \leq h \} \leq q^{C'(n,l)h^{l+1}}
\]
for all $h \geq 1$.
By re-ordering the coordinate of $(\PP^1_{\FF_q})^n$,
we can find an automorphism
\[
\iota : (\PP^1_{\FF_q})^n \to (\PP^1_{\FF_q})^n
\]
with $\pi_{[l]} \cdot \iota = \pi_{I}$ and
$\iota^*(\OO_{(\PP^1)^n}(1, \ldots, 1)) = \OO_{(\PP^1)^n}(1, \ldots, 1)$.
Thus, we may assume that $I = [l]$.
We denote $p_{[l]}$ by $p$.
Let $p_i : (\PP^1_{\FF_q})^{n} \to \PP^1_{\FF_q}$ be the projection
to the $i$-th factor. For $n \geq l+1$, we set
\[
T_n = Z^{\rm eff}_l((\PP^1_{\FF_q})^{n} \overset{p}{\to} (\PP^1_{\FF_q})^l)
\]
and
$h_n(V) = \deg_{\OO(1, \ldots, 1)}(V)$ for $V \in T_n$.
Let $a_n : (\PP^1_{\FF_q})^n \to (\PP^1_{\FF_q})^{n-1}$ and
$b_n : (\PP^1_{\FF_q})^n \to (\PP^1_{\FF_q})^{l+1}$ be morphisms
given by $a_n = p_{[n-1]}$ and $b_n = p_{[l] \cup \{ n \}}$, namely,
\[
a_n(x_1, \ldots, x_n) = (x_1, \ldots, x_{n-1})
\quad\text{and}\quad
b_n(x_1, \ldots, x_n) = (x_1, \ldots, x_l, x_n).
\]
Here, $\alpha_n : T_n \to T_{n-1}$ and
$\beta_n : T_n \to T_{l+1}$ are given by
\[
\alpha_n(V) = (a_n)_*(V)
\quad\text{and}\quad
\beta_n(V) = (b_n)_*(V).
\]
Then, since
\[
\OO_{(\PP^1)^{n}}(1, \ldots, 1) =
(a_n)^*(\OO_{(\PP^1)^{n-1}}(1,\ldots,1)) 
\otimes p_{n}^*(\OO_{\PP^1}(1))
\]
and
\[
\OO_{(\PP^1)^{n}}(1, \ldots, 1) =
(b_n)^*(\OO_{(\PP^1)^{l+1}}(1,\ldots,1))
\otimes \bigotimes_{j=l+1}^{n-1} p_j^*(\OO_{\PP^1}(1)),
\]
it is easy to see that
$h_{n-1}(\alpha_n(V)) \leq h_n(V)$ and $h_{l+1}(\beta_n(V)) \leq h_n(V)$.
Note that the diagram
\[
\xymatrix{
& (\PP^1_{\FF_q})^{n} \ar[dl]_{a_n} 
\ar[dr]^{b_n} & \\
(\PP^1_{\FF_q})^{n-1} \ar[dr]_{p_{[l]}} & & (\PP^1_{\FF_q})^{l+1} \ar[dl]^{p_{[l]}} \\
& (\PP^1_{\FF_q})^l & 
}
\]
is a fiber product.
Thus, by Lemma~\ref{lem:number:product:X:Y},
if we set $A(s,t) = q^{st}$, then, for $y \in T_{n-1}$ and $z \in T_{l+1}$,
\[
\# \{ x \in T_n \mid \alpha_n(x) = y, \ \beta_n(x) = z \} \leq
A(h_{n-1}(y), h_{l+1}(z)).
\]
Here,
\[
\{ D \in T_{l+1} \mid h_{l+1}(D) \leq h \} \subseteq 
\{ D \in T_{l+1} \mid \text{$\deg_i(D) \leq h$ for all $i=1, \ldots, l+1$} \}.
\]
Thus, by Proposition~\ref{prop:estimate:X:product:P:1:divisor},
if we set $B(h) = (1+h)^{l+1}q^{(1+h)^{l+1}}$, then
\[
\# \{ D \in T_{l+1} \mid h_{l+1}(D) \leq h \} \leq B(h).
\]
Therefore, by Lemma~\ref{lem:counting:system},
\begin{align*}
\# \{ x \in T_n \mid h_n(x) \leq h \} & \leq
B(h)^{n-l}A(h,h)^{n-l-1} \\
& = (1+h)^{(n-l)(l+1)}q^{(n-l)(1+h)^{l+1}}q^{(n-l-1)h^2} \\
& \leq q^{(n-l)(l+1)h} q^{(n-l)(2h)^{l+1}}q^{(n-l-1)h^2} \leq q^{(n-l)(2^{l+1}+l+2) h^{l+1}}
\end{align*}
for all $h \geq 1$.
\QED

Moreover, we have the following variant of Proposition~\ref{prop:growth:P:1:product:l}.

\begin{Proposition}
\label{prop:relative:P:1:product}
Let $d$, $l$ and $n$ be positive integers with $d \leq l \leq n$.
Let $p_{[d]} : (\PP^1_{\FF_q})^{n} \to (\PP^1_{\FF_q})^d$ be the morphism
given in \rom{\ref{subsec:product:P:1:projection}}.
Let $p_i : (\PP^1_{\FF_q})^{n} \to \PP^1_{\FF_q}$ be the projection to the $i$-th factor.
We set $L_n = \bigotimes_{i=1}^{n}p_i^*(\OO(1))$ and
$H_n = \bigotimes_{i=1}^d p_i^*(\OO(1))$. Then, for a fixed $k$,
there is a constant $C$ such that
\[
\# \left\{ V \in Z^{\rm eff}_l((\PP^1_{\FF_q})^{n} \overset{p_{[d]}}{\to}
(\PP^1_{\FF_q})^d) \left| 
\begin{array}{l}
\deg(L_n^{\cdot l - d} \cdot H_n^{\cdot d} \cdot V) \leq k, \\
\deg(L_n^{\cdot l - d +1} \cdot H_n^{\cdot d-1} \cdot V) \leq h
\end{array}
\right\} \right.
\leq q^{C\cdot h^{d}}
\]
for all $h \geq 1$.
\end{Proposition}

\Proof
We set
\[
\Sigma = \{ I \mid [d] \subseteq I \subseteq [n], \#( I ) = l \}.
\]
Then,
\[
Z^{\rm eff}_l((\PP^1_{\FF_q})^{n} \overset{p_{[d]}}{\to}
(\PP^1_{\FF_q})^d) = \sum_{I \in \Sigma}
Z^{\rm eff}_l((\PP^1_{\FF_q})^{n} \overset{p_{I}}{\to}
(\PP^1_{\FF_q})^l).
\]
Thus, it is sufficient to show that there is a constant $C'$
\[
\# \left\{ V \in Z^{\rm eff}_l((\PP^1_{\FF_q})^{n} \overset{p_I}{\to}
(\PP^1_{\FF_q})^l) \left| 
\begin{array}{l}
\deg(L_n^{\cdot l - d} \cdot H_n^{\cdot d} \cdot V) \leq k, \\
\deg(L_n^{\cdot l - d -1} \cdot H_n^{\cdot d-1} \cdot V) \leq h
\end{array}
\right\} \right.
\leq q^{C' \cdot h^{d}}
\]
for all $h \geq 1$.
Re-ordering the coordinate of $(\PP^1_{\FF_q})^n$, we may assume that
$I = [l]$. We denote $p_I$ by $p$.
Let $a_n : (\PP^1_{\FF_q})^{n} \to (\PP^1_{\FF_q})^{n-1}$ and
$b_n : (\PP^1_{\FF_q})^{n} \to (\PP^1_{\FF_q})^{l+1}$ be
morphisms given by $a_n = p_{[n-1]}$ and $b_n = p_{[l] \cup \{ n \}}$, i.e.,
$a_n(x_1, \ldots, x_{n}) = (x_1, \ldots, x_{n-1})$ and
$b_n(x_1, \ldots, x_{n}) = (x_1, \ldots, x_l, x_{n})$.
Then,
\addtocounter{Claim}{1}
\begin{equation}
\label{eqn:lem:relative:P:1:product:1}
\begin{cases}
a_n^*(L_{n-1}) \otimes p_{n}^*(\OO(1)) = L_n \\
b_n^*(L_{l+1}) \otimes \bigotimes_{i=l+1}^{n-1} p_i^*(\OO(1))  = L_n \\
a_n^*(H_{n-1}) = H_n \\
b_n^*(H_{l+1}) = H_n \\
\end{cases}
\end{equation}
Here, for $n \geq l+1$, we set
\[
T_n = \{ V \in Z^{\rm eff}_l((\PP^1_{\FF_q})^{n} \overset{p}{\to}
(\PP^1_{\FF_q})^l) \mid \deg(L_n^{\cdot l - d} \cdot H_n^{\cdot d} \cdot V) \leq k \}.
\]
Let $h_n : T_n \to \RR$ be a map given by $h_n(V) = \deg(L_n^{\cdot l - d +1} \cdot H_n^{\cdot d-1} \cdot V)$.
Then, by \eqref{eqn:lem:relative:P:1:product:1}, we have
maps $\alpha_n : T_n \to T_{n-1}$ and $\beta_n : T_n \to T_{l+1}$
given by $\alpha_n(V) = (a_n)_*(V)$ and $\beta_n(V) = (b_n)_*(V)$.
Moreover, $h_{n-1}(\alpha_n(V)) \leq h_n(V)$ and $h_{l+1}(\beta_n(V)) \leq h_n(V)$ for
all $V \in T_n$.
As in Lemma~\ref{lem:number:product:X:Y},
we denote by $\theta(V)$ the number of irreducible components of a cycle $V$.
Then, for $V \in T_n$, it is easy to see that $\theta(V) \leq k$.
Further,
\[
p^*(H_l) = H_n
\quad\text{and}\quad
p^*(L_l) \otimes \bigotimes_{i=l+1}^{n} p_i^*(\OO(1)) = L_n.
\]
Thus, by Lemma~\ref{lem:number:product:X:Y},
for $D \in T_{n-1}$ and $E \in T_{l+1}$,
\begin{multline*}
\log_q \# \{ V \in T_n \mid \alpha_n(V) = D, \beta_n(V) = E \} \\
\leq \frac{k\sqrt{\deg(L_{n-1}^{\cdot l-d+1} \cdot H_{n-1}^{\cdot d-1} \cdot D)
\deg(L_{l+1}^{\cdot l-d+1} \cdot H_{l+1}^{\cdot d-1} \cdot E)}}{\deg\left(L_l^{\cdot l-d+1} \cdot H_l^{\cdot d-1} \right)}.
\end{multline*}
Thus, if we set $A(x,y) = q^{k\sqrt{xy}}$, then
\[
\# \{ V \in T_n \mid \alpha_n(V) = D, \beta_n(V) = E \} \leq A(h_{n-1}(D), h_{l+1}(E)).
\]
Here let us estimate $\# \{ D \in T_{l+1} \mid h_{l+1}(D) \leq h \}$. In this case,
$D$ is a divisor on $(\PP^1_{\FF_q})^{l+1}$.
Thus,
\[
\deg(L_{l+1}^{\cdot l - d} \cdot H_{l+1}^{\cdot d} \cdot D) = d! (l-d)!
(\deg_{d+1}(D) + \cdots + \deg_{l+1}(D))
\]
and
\[
\deg(L_{l+1}^{\cdot l - d+1} \cdot H_{l+1}^{\cdot d-1} \cdot D) = (d-1)!(l-d+1)!
(\deg_1(D) + \cdots  + \deg_{l+1}(D)).
\]
Therefore, $\deg_i(D) \leq h$ for $1 = 1, \ldots d$ and
$\deg_j(D) \leq k$ for $j=d+1, \ldots, l+1$.
Hence, by Proposition~\ref{prop:estimate:X:product:P:1:divisor},
if we set $B(h) = q^{C_1 \cdot h^{d}}$ for some constant $C_1$, then
\[
\# \{ D \in T_{l+1} \mid h_{l+1}(D) \leq h \} \leq B(h)
\]
for $h \geq 1$.

\medskip
Gathering the above observations and using Lemma~\ref{lem:counting:system},
\[
\# \{ V \in T_n \mid h_n(V) \leq h \} \leq  B(h)^{n-l} A(h, h)^{n-l-1}
\leq q^{(n-l) C_1 \cdot h^{d} + k(n-l-1)h}.
\]
for $h \geq 1$. Hence, we get our proposition.
\QED

\begin{Remark}
\label{remark:prop:relative:P:1:product}
The following are remarks for Proposition~\ref{prop:relative:P:1:product}.
\begin{enumerate}
\renewcommand{\labelenumi}{(\arabic{enumi})}
\item
$\deg(L_n^{\cdot l-d} \cdot H_n^{\cdot d} \cdot V) = d!
\deg((L_n)_{\eta}^{\cdot l-d} \cdot V_{\eta})$,
where the subscript $\eta$ means the restriction of an object
on $(\PP^1_{\FF_q})^n$ to the generic fiber of $p_{[d]} :
(\PP^1_{\FF_q})^n \to (\PP^1_{\FF_q})^d$.

\item
If we set $L'_n = \bigotimes_{i=d+1}^n p_i^*(\OO(1))$, then, for a fixed $k$,
there is a constant $C'$ such that
\[
\# \left\{ V \in Z^{\rm eff}_l((\PP^1_{\FF_q})^{n} \overset{p_{[d]}}{\to}
(\PP^1_{\FF_q})^d) \left| 
\begin{array}{l}
\deg({L'_n}^{\cdot l - d} \cdot H_n^{\cdot d} \cdot V) \leq k, \\
\deg({L'_n}^{\cdot l - d +1} \cdot H_n^{\cdot d-1} \cdot V) \leq h
\end{array}
\right\} \right.
\leq q^{C' \cdot h^{d}}
\]
for all $h \geq 1$.
($\because$ Note that $L_n = L'_n + H_n$. Thus,
$L_n^{\cdot l-d} \cdot H_n^{\cdot d} = {L'_n}^{\cdot l-d} \cdot H_n^{\cdot d}$
and $L_n^{\cdot l-d+1} \cdot H_n^{\cdot d-1} = {L'_n}^{\cdot l-d+1} \cdot H_n^{\cdot d-1}
+ (l-d+1) {L'_n}^{\cdot l-d} \cdot H_n^{\cdot d}$.)
\end{enumerate}
\end{Remark}

\subsection{Cycles on a projective variety over a finite field}
\setcounter{Theorem}{0}

Here we consider the main problem of this section.
Let us begin with two lemmas.

\begin{Lemma}
\label{lem:comp:degree:field}
Let $F$ be a field.
Let $\phi : \PP^{n}_{F} \dasharrow (\PP^1_{F})^{n}$ be the birational map
given by
\[
(X_0 : \ldots : X_{n}) \mapsto (X_0 : X_1) \times \cdots \times (X_0: X_{n}).
\]
Let $\Sigma$ be the boundary of $\PP^{n}_{F}$, that is, 
$\Sigma = \{ X_0 = 0 \}$.
For $V \in Z^{\rm eff}_l(\PP^{n}_{F} ;
\PP^{n}_{F} \setminus \Sigma)$,
we denote by $V'$ the strict transform of $V$ by $\phi$
\rom{(}For the definition of 
$Z^{\rm eff}_l(\PP^{n}_{F} ;
\PP^{n}_{F} \setminus \Sigma)$,
see \rom{\ref{subsec:notation:cycle}}\rom{)}.
Then,
\[
n^{l} \deg\left(
\OO_{\PP^n_F}(1)^{\cdot l} \cdot V
\right)
\geq \deg\left(
\OO_{(\PP^1_{F})^n}(1,\ldots,1)^{\cdot l} \cdot V'
\right)
\]
for all $V \in Z^{\rm eff}_l(\PP^{n}_{F};\PP^{n}_{F} \setminus \Sigma)$.
\end{Lemma}

\Proof
Let $Y$ ($\subseteq \PP^n_{F} \times (\PP^1_{F})^n)$
be the graph of the rational map
$\phi : \PP^{n}_{F} \dasharrow (\PP^1_{F})^{n}$.
Let $\mu : Y \to \PP^{n}_{F}$ and
$\nu : Y \to  (\PP^1_{F})^{n}$ be the morphisms induced by
the projections
$ \PP^n_{F} \times (\PP^1_{F})^n \to \PP^n_F$ and
$ \PP^n_{F} \times (\PP^1_{F})^n \to (\PP^1_F)^n$
respectively.
Here we claim that
there is an effective Cartier divisor $E$ on $Y$
such that
(1) $\mu(E) \subseteq \Sigma$ and 
(2) ${\mu}^*(\OO(n)) = {\nu}^*(\OO(1,\ldots,1)) \otimes \OO_{Y}(E)$.
Let $Y_i$ ($\subseteq \PP^n_{F} \times \PP^1_{F})$ be the graph
of the rational map
$\PP^n_{F} \dasharrow \PP^1_{F}$ given by
\[
(X_0 : \cdots : X_n) \mapsto (X_0 : X_i).
\]
Let $\mu_i : Y_i \to \PP^n_{F}$ and $\nu_i : Y_i \to \PP^1_{F}$
be the morphisms induced by
the projections $\PP^n_{F} \times \PP^1_{F} \to \PP^n_{F}$
and $\PP^n_{F} \times \PP^1_{F} \to \PP^1_{F}$
respectively.
Let $\pi_i : (\PP^1_{F})^n \to \PP^1_{F}$ be the projection
to the $i$-th factor.
Moreover, let $h_i : Y \to Y_i$ be the morphism induced by
$\operatorname{id} \times \pi_i : 
\PP^n_{F} \times (\PP^1_{F})^n \to \PP^n_{F} \times \PP^1_{F}$.
Consequently, we have the following commutative diagram:
\[
\xymatrix{
& Y \ar[ddl]_{\mu} \ar[d]^{h_i} \ar[dr]^{\nu} \\
& Y_i \ar[dl]^{\mu_i} \ar[dr]_{\nu_i} & (\PP^1_{F})^n \ar[d]^{\pi_i} \\
\PP^n_{F} & & \PP^1_{F} 
}
\]
Note that $Y_i$ is the blowing-up by the ideal sheaf $I_i$ generated
by $X_0$ and $X_i$.
Thus there is an effective Cartier divisor $E_i$ on $Y_i$
with $I_i \OO_{Y_i} = \OO_{Y_i}(-E_i)$ and
$\mu_i^*(\OO_{\PP^n}(1)) \otimes \OO_{Y_i}(-E_i) = \nu_i^*(\OO_{\PP^1}(1))$. 
Thus if we set $E = \sum_{i=1}^n h_i^*(E_i)$, then
\[
\mu^*(\OO(n)) = 
\nu^*(\OO(1, \ldots, 1)) \otimes \OO_Y(E).
\]
Hence we get our claim.

\medskip
For $V \in Z^{\rm eff}_l(\PP^{n}_{F}  ;
\PP^{n}_{F} \setminus \Sigma)$, let $V''$ be the strict transform of $V$ by $\mu$.
Then, by using the projection formula,
\[
\deg\left(
\OO(n)^{\cdot l}  \cdot V
\right) 
= \deg\left(
\mu^*(\OO(n))^{\cdot l} \cdot V''\right).
\]
Moreover, by the following Lemma~\ref{lem:comp:nef:eff}
\[
\deg\left(
\mu^*(\OO(n))^{\cdot l} \cdot V''\right) 
\geq \deg\left(
\nu^*(\OO(1,\ldots,1))^{\cdot l}  \cdot V''\right)
\]
Thus, using the projection formula for $\nu$,
we get our lemma because $\nu_*(V'') = V'$.
\QED

\begin{Lemma}
\label{lem:comp:nef:eff}
Let $X$ be a projective variety over a field $F$ and
$L_1, \ldots, L_{\dim X}, M_1, \ldots, M_{\dim X}$ nef line
bundles on $X$. If $L_i \otimes M_i^{\otimes -1}$ is pseudo-effective
for $i=1, \ldots, n$, then
\[
\deg(L_1 \cdots L_{\dim X}) \geq \deg(M_1 \cdots M_{\dim X}).
\]
\end{Lemma}

\Proof
We set $E_i = L_i \otimes M_i^{\otimes -1}$. Then
\[
\deg(L_1 \cdots L_{\dim X}) = \deg(M_1 \cdots M_{\dim X}) +
\sum_{i=1}^{\dim X} \deg(M_1 \cdots M_{i-1} \cdot
E_i \cdot L_{i+1} \cdots L_{\dim X}).
\]
Thus, we get our lemma.
\QED

\begin{Theorem}
\label{thm:proj:sp:geom:case}
There is a constant $C(n,l)$ depending only on $n$ and $l$
such that
\[
\#\left( \{ V \in Z_l^{\rm eff}(\PP^n_{\FF_q}) \mid
\deg_{\OO(1)}(V) \leq h) \} \right)
\leq q^{C(n,l)h^{l+1}}
\]
for all $h \geq 1$.
\end{Theorem}

\Proof
We prove this theorem by induction on $n$.
Let us consider the birational map 
$\PP^n_{\FF_q} \dasharrow (\PP^1_{\FF_q})^n$ given by
\[
\phi : (X_0 : \cdots : X_n) \mapsto (X_0 : X_1) \times \cdots \times (X_0 : X_n).
\]
We set $U = \PP^n_{\FF_q} \setminus \{ X_0 = 0 \}$.
For $V \in Z_l^{\rm eff}(\PP^n_{\FF_q};U)$, we denote by $V'$ 
the strict transform of $V$ by $\phi$. 
Then, by Lemma~\ref{lem:comp:degree:field},
\[
n^l \deg_{\OO(1)}(V) \geq \deg_{\OO(1, \ldots, 1)}(V').
\]
Moreover, note that
if $V'_1 = V'_2$
for $V_1, V_2 \in Z_l^{\rm eff}(\PP^n_{\FF_q};U)$, then $V_1 = V_2$.
Therefore
\[
\# \{ V \in Z_l^{\rm eff}(\PP^n_{\FF_q};U) 
\mid \deg_{\OO(1)}(V) \leq h \} 
\leq
\# \{ V' \in Z_l^{\rm eff}((\PP^1_{\FF_q})^n) \mid 
\deg_{\OO(1,\ldots,1)}(V') 
\leq n^l h \}.
\]
Here, by Proposition~\ref{prop:growth:P:1:product:l},
there is a constant $C'(n,l)$ depending only $n$ and $l$
such that
\[
\# \{ V' \in Z^{\rm eff}_{l}( (\PP^1_{\FF_q})^n) \mid 
\deg_{\OO(1,\ldots,1)}(V') \leq k \} 
\leq  q^{C'(n,l){k}^{l+1}}.
\]
Hence, we have
\addtocounter{Claim}{1}
\begin{equation}
\label{thm:proj:sp:geom:case:eq:2}
\# \{ V \in Z_l^{\rm eff}(\PP^n_{\FF_q};U) 
\mid \deg_{\OO(1)}(V) \leq h \} 
\leq
q^{C'(n,l)n^{l(l+1)} h^{l+1}}.
\end{equation}
On the other hand, since $\PP^n_{\FF_q} \setminus U \simeq
\PP^{n-1}_{\FF_q}$, 
\begin{multline*}
\# \{ V \in Z_l^{\rm eff}(\PP^n_{\FF_q}) \mid \deg_{\OO(1)}(V) \leq h \} \\
\leq \# \{ V \in Z_l^{\rm eff}(\PP^n_{\FF_q};U) \mid 
\deg_{\OO(1)}(V) \leq h \} \cdot
\# \{ V \in Z_l^{\rm eff}(\PP^{n-1}_{\FF_q}) \mid 
\deg_{\OO(1)}(V) \leq h \}
\end{multline*}
Thus, using the hypothesis of induction, 
if we set $C(n,l) = C(n-1, l) + n^{l(l+1)}C'(n,l)$,
then we have our theorem.
\QED

\begin{Corollary}
\label{cor:general:geom:case}
Let $X$ be a projective variety over a finite field $\FF_q$ and
$H$ a very ample line bundle on $X$. 
Then, for every integer $l$
with $0 \leq l \leq \dim X$,
there is a constant $C$ depending only on
$l$ and $\dim_{\FF_q} H^0(X, H)$ such that
\[
\# \{ V \in Z_l^{\rm eff}(X) \mid \deg_{H}(V) \leq h \}
\leq q^{C h^{l+1}}
\]
\end{Corollary}

\Proof
Since $H$ is very ample,
there is an embedding $\iota : X \to \PP^n_{\FF_q}$
with $\iota^*(\OO(1)) = H$, where $n = \dim_{\FF_q} H^0(X, H) - 1$.
Thus, it follows from Theorem~\ref{thm:proj:sp:geom:case}.
\QED

Finally, let us consider a lower bound
of effective cycles with bounded degree.

\begin{Proposition}
\label{prop:general:geom:case:lower:bound}
Let $X$ be a projective variety over a finite field $\FF_q$ and
$H$ an ample line bundle on $X$. Then, for every integer $l$
with $0 \leq l < \dim X$,
\[
 \limsup_{h\to\infty}
\frac{\log \#\left( \{ V \in Z_l^{\rm eff}(X) \mid \deg_{H}(V) =  k \}
\right)}{k^{l+1}}
> 0.
\]
\end{Proposition}

\Proof
Take $(l+1)$-dimensional subvariety $Y$ of $X$.
Then, 
\[
\#\left( \{ V \in Z_l^{\rm eff}(Y) \mid \deg_{\rest{H}{Y}}(V) =  k \} \right) \leq
\#\left( \{ V \in Z_l^{\rm eff}(X) \mid \deg_{H}(V) =  k \} \right).
\]
Thus, we may assume $l = \dim X - 1$.
Here, note that
\[
\vert H^{\otimes m} \vert \subseteq \{ D \in Z_{d-1}^{\rm eff}(X) \mid 
\deg_{H}(D) = m(H^d) \}
\quad\text{and}\quad
\# \vert H^{\otimes m} \vert = \frac{q^{\dim_{\FF_q} H^0(X, H^{\otimes m})}-1}{q-1},
\]
where $d = \dim X$.
Since $H$ is ample, $\dim_{\FF_q} H^0(X, H^{\otimes m}) = O(m^d)$.
Thus we get our proposition.
\QED

\subsection{Northcott's type}
\setcounter{Theorem}{0}
Let $B$ be a $d$-dimensional projective variety over $\FF_q$.
We assume that $d \geq 1$.
Let $H$ be a nef and big line bundle on $B$.
Let $X$ be a projective variety over $\FF_q$ and
$f : X \to B$ a surjective morphism over $\FF_q$.
Let $L$ be a nef line bundle on $X$.
In the following, the subscript $\eta$ of
an object on $X$ means its restriction on the generic
fiber of $f : X \to B$.

\begin{Theorem}
\label{thm:relative:case:finite:field}
If $L_{\eta}$ is ample, then, for a fixed $k$,
there is a constant $C$ such that
\[
\# \{ V \in Z^{\rm eff}_l(X/B) \mid
\deg(L_{\eta}^{\cdot l - d} \cdot V_{\eta}) \leq k,\ 
\deg(L^{\cdot l - d + 1} \cdot f^*(H)^{\cdot d-1} \cdot V) \leq h
\} \leq q^{C \cdot h^{d}}
\]
for all $h \geq 1$.
\end{Theorem}

\Proof
{\bf Step 1:}
We set $B = (\PP^1_{\FF_q})^d$ and
$X = (\PP^1_{\FF_q})^d \times (\PP^1_{\FF_q})^e = (\PP^1_{\FF_q})^{d+e}$.
Let $p_i : B \to \PP^1_{\FF_q}$ be the projection to the $i$-th factor.
Similarly, let $q_j : X \to \PP^1_{\FF_q}$ be the projection to
the $j$-th factor. Here $f : X \to B$ is given the natural projection
$q_1 \times \cdots \times q_d$, namely,
$f(x_1, \ldots, x_{d+e}) = (x_1, \ldots, x_d)$.
Moreover, we set $H = p_1^*(\OO(1)) \otimes \cdots \otimes p_{d+e}^*(\OO(1))$ and
$L = q_{d+1}^*(\OO(1)) \otimes \cdots \otimes q_{d+e}^*(\OO(1))$.
In this situation, we have our theorem by
Proposition~\ref{prop:relative:P:1:product}
and Remark~\ref{remark:prop:relative:P:1:product}.

\bigskip
{\bf Step 2:}
Let us consider a case where $X = B \times (\PP^1_{\FF_q})^e$,
$f : X \to B$ is given by the projection to the first factor and
$L = q^*(\OO(1, \ldots, 1))$.
Here $q : X \to (\PP^1_{\FF_q})^d$ is the natural projection.
By virtue of Noether's normalization theorem, there is a dominant
rational map $B \dasharrow (\PP^1_{\FF_q})^d$.
Let the following diagram
\[
\xymatrix{
& B' \ar[dl]_{\nu} \ar[dr]^{\nu'}& \\
B & & (\PP^1_{\FF_q})^d
}
\]
be the graph of the rational map $B \dasharrow (\PP^1)^d$.
Here we set $X' = B' \times  (\PP^1_{\FF_q})^e$, $B'' = (\PP^1_{\FF_q})^d$ and
$X'' =  (\PP^1_{\FF_q})^d  \times (\PP^1_{\FF_q})^e$.
Let $f': X' \to B'$, $q' : X' \to (\PP^1_{\FF_q})^e$,
$f'' : X'' \to  (\PP^1_{\FF_q})^d$, and
$q'' : X'' \to  (\PP^1_{\FF_q})^e$ be the natural projections.
Moreover, we set $L' = {q'}^*(\OO(1, \ldots, 1))$ and
$L'' = {q''}^*(\OO(1, \ldots, 1))$.

Let $\alpha : Z^{\rm eff}_l(X/B) \to Z^{\rm eff}_l(X'/B')$ be a homomorphism
given by the strict transform in terms of $\nu \times \operatorname{id} : X' \to X$.
Further, let $\beta : Z^{\rm eff}_l(X'/B') \to Z^{\rm eff}_l(X''/B'')$ be the homomorphism
given by the push forward $(\nu' \times \operatorname{id})_*$ of cycles.
Since $H$ is nef and big, there is a positive integer $a$
such that
$H^0(B', \nu^*(H)^{\otimes a} \otimes {\nu'}^*(\OO(-1, \ldots, -1))) \not= 0$.
Then, by Lemma~\ref{lem:comp:nef:eff},
\begin{align*}
a^{d-1} \deg(L^{\cdot l - d + 1} \cdot f^*(H)^{\cdot d-1} \cdot V) & =
a^{d-1} \deg((\nu \times \operatorname{id})^*(L)^{\cdot l - d + 1}
\cdot (\nu \times \operatorname{id})^*({f}^*(H))^{\cdot d-1} \cdot \alpha(V)) \\
& = \deg({L'}^{\cdot l - d + 1} \cdot {f'}^*(\nu^*(H)^{\otimes a})^{\cdot d-1} \cdot \alpha(V)) \\
& \geq \deg({L'}^{\cdot l - d + 1} \cdot
{f'}^*({\nu'}^*(\OO(1, \ldots, 1))^{\cdot d-1} \cdot \alpha(V)) \\
& = \deg({L''}^{\cdot l - d + 1} \cdot
{f''}^*(\OO(1, \ldots, 1))^{\cdot d-1} \cdot \beta(\alpha(V))).
\end{align*}
Moreover,
\begin{align*}
\deg(L_{\eta}^{\cdot l - d + 1} \cdot V_{\eta})
& = \deg({L'}_{\eta'}^{\cdot l - d + 1} \cdot \alpha(V)_{\eta'}) \\
& = \deg({L''}_{\eta''}^{\cdot l - d + 1} \cdot
\beta(\alpha(V))_{\eta''}),
\end{align*}
where the subscripts $\eta'$ and $\eta''$ means the restrictions of objects
to the generic fibers $f'$ and $f''$ respectively.

For a fixed  $V'' \in Z^{\rm eff}_l(X''/B'')$,
we claim that
\[
\log_q \#\{V' \in Z^{\rm eff}_l(X'/B') \mid \beta(V') = V'' \}
\leq \deg(\nu')\deg({L''}_{\eta''}^{l-d+1} \cdot V''_{\eta''}).
\]
Let $B_0''$ be the maximal Zariski open set of $B''$ such that
$\nu'$ is finite over $B_0''$. We set $B_0' = {\nu'}^{-1}(B_0'')$.
Then, the natural homomorphisms
\[
Z^{\rm eff}_l(X'/B') \to Z^{\rm eff}_l(X'_0/B'_0)
\quad\text{and}\quad
Z^{\rm eff}_l(X''/B'') \to Z^{\rm eff}_l(X''_0/B''_0)
\]
are bijective, where
$X_0' = B_0' \times (\PP^1_{\FF_q})^e$ and
$X_0'' = B_0'' \times (\PP^1_{\FF_q})^e$.
Thus, by virtue of Lemma~\ref{lem:estimate:number:push:forward},
if we set $V'' = \sum_i a_i W_i$, then
\[
\log_q \#\{V' \in Z^{\rm eff}_l(X'/B') \mid \beta(V') = V'' \}
\leq \deg(\nu')\sum_i a_i.
\]
On the other hand,
\[
\sum_i a_i \leq \sum_i a_i \deg({L''}_{\eta''}^{l-d+1} \cdot {W_i}_{\eta''}) = 
\deg({L''}_{\eta''}^{l-d+1} \cdot V''_{\eta''}).
\]
Therefore, we get our claim.

\medskip
Hence, by the above observations and Step~1,
we have our case.

\bigskip
{\bf Step 3:}
Next let us consider a case where 
$X = B \times \PP^e_{\FF_q}$, $f : X \to B$ is given the natural projection and
$L = q^*(\OO(1))$. Here $q : X \to \PP^e_{\FF_q}$ is the natural projection.
We prove our theorem of this situation by induction on $e$.
If $e = l -d$, then our assertion is obvious. Thus we assume that $e > l - d$.
Let the following diagram
\[
\xymatrix{
& Y \ar[dl]_{\mu} \ar[dr]^{\nu}& \\
\PP^e_{\FF_q} & & (\PP^1_{\FF_q})^e
}
\]
be the graph of the rational map $\PP^e_{\FF_q} \dasharrow (\PP^1_{\FF_q})^e$
given by
\[
(X_0 : \cdots : X_e) \mapsto (X_0 : X_1) \times \cdots \times (X_0 : X_e).
\]
Then, as in Lemma~\ref{lem:comp:degree:field},
there is an effective Cartier divisor $E$ on $Y$ such that
$\mu(E) \subset \{ X_0 = 0 \}$ and
$\mu^*(\OO(e)) = \nu^*(\OO(1,\ldots,1)) \otimes \OO_Y(E)$.
Here we set $X' = B \times (\PP^1_{\FF_q})^e$ and
$L' = {q'}^*(\OO(1,\ldots,1))$, where $q' : X \to (\PP^1_{\FF_q})^e$
is the natural projection. Moreover,
$f' : X' \to B$ is given by the natural projection.
Then, for $V \in Z_l^{\rm eff}(X; X \setminus B\times\{X_0=0\})$,
by Lemma~\ref{lem:comp:nef:eff},
\begin{align*}
e^{l-d+1} \deg(L^{\cdot l-d+1} \cdot f^*(H)^{\cdot d-1} \cdot V) & =
\deg((\mu \times \operatorname{id})^*(L^{\otimes e})^{\cdot l-d+1} \cdot 
(\mu \times \operatorname{id})^*f^*(H)^{\cdot d-1} \cdot V') \\
& \geq \deg((\nu \times \operatorname{id})^*(L')^{\cdot l-d+1} \cdot 
(\nu \times \operatorname{id})^*{f'}^*(H)^{\cdot d-1} \cdot V') \\
& = \deg({L'}^{\cdot l-d+1} \cdot 
{f'}^*(H)^{\cdot d-1} \cdot (\nu \times \operatorname{id})_*(V')),
\end{align*}
where $V'$ is the strict transform of $V$ by $\mu \times \operatorname{id}$.
Further,
\begin{align*}
e^{l-d} \deg(L_{\eta}^{\cdot l-d} \cdot V_{\eta}) & =
\deg((\mu \times \operatorname{id})^*(L^{\otimes e})_{\eta''}^{\cdot l-d} \cdot V'_{\eta''}) \\
& \geq \deg((\nu \times \operatorname{id})^*(L')_{\eta''}^{\cdot l-d} \cdot V'_{\eta''}) \\
& = \deg({L'}_{\eta'}^{\cdot l-d} \cdot (\nu \times \operatorname{id})_*(V')_{\eta'}),
\end{align*}
where $\eta'$ and $\eta''$ means the restriction of objects on $X'$ and
$B \times Y$ to the generic fibers $X' \to B$ and
$B \times Y \to B$ respectively.
Here $B\times\{X_0=0\} \simeq B \times \PP^{e-1}_{\FF_q}$. Thus, by hypothesis of induction and
Step~2, we have our case.

\bigskip
{\bf Step 4:}
Finally we consider a general case.
Clearly we may assume that $L_{\eta}$ is very ample.
Thus, there are a positive integer $e$ and a subvariety $X'$ of $B \times \PP^e_{\FF_q}$ with the following properties:
\begin{enumerate}
\renewcommand{\labelenumi}{(\arabic{enumi})}
\item
Let $f' : X' \to B$ (resp. $q : X' \to \PP^e_{\FF_q}$) be the projection to the first factor
(resp. the second factor).
There is a non-empty Zariski open set $B_0$ of $B$ such that
$f^{-1}(B_0)$ is isomorphic to ${f'}^{-1}(B_0)$ over $B_0$.
We denote this isomorphism by $\iota$.

\item
If we set $L' = q^*(\OO(1))$, then
$\rest{L}{f^{-1}(B_0)} = \iota^*\left( \rest{L'}{{f'}^{-1}(B_0)} \right)$.
\end{enumerate}
Let
\[
\xymatrix{
& X'' \ar[dl]_{\mu} \ar[dr]^{\mu'} & \\
X \ar[dr]_{f} & & X' \ar[dl]^{f'} \\
& B & }
\]
be the graph of the rational map induced by $\iota$.
We denote $f \cdot \mu = f' \cdot \mu'$ by $f''$.
By the property (2),
\[
f''_*(\mu^*(L) \otimes {\mu'}^*({L'}^{\otimes -1})) \not= 0.
\]
Thus, we can find an ample line bundle $A$ on $B$ such that
\[
H^0(X'', \mu^*(L \otimes f^*(A)) \otimes {\mu'}^*({L'}^{\otimes -1})) \not= 0.
\]
Let us choose a non-zero element $s$ of
\[
H^0(X'', \mu^*(L \otimes f^*(A)) \otimes {\mu'}^*({L'}^{\otimes -1})).
\]
Since $f''(\Supp(\zeros(s)) \not= B$,
by Lemma~\ref{lem:comp:nef:eff},
\begin{align*}
\deg((L \otimes f^*(A))^{\cdot l-d+1} \cdot f^*(H)^{\cdot d-1} \cdot V) & =
\deg(\mu^*(L \otimes f^*(A))^{\cdot l-d+1} \cdot \mu^* f^*(H)^{\cdot d-1} \cdot V') \\
& \geq
\deg({\mu'}^*(L')^{\cdot l-d+1} \cdot {\mu'}^*{f'}^*(H)^{\cdot d-1} \cdot V') \\
& = \deg({L'}^{\cdot l-d+1} \cdot {f'}^*(H)^{\cdot d-1} \cdot\mu'_*(V')),
\end{align*}
where $V'$ is the strict transform of $V$ by $\mu$.
Moreover,
\begin{multline*}
\deg((L + f^*(A))^{\cdot l-d+1} \cdot f^*(H)^{\cdot d-1} \cdot V) =
\deg(L^{\cdot l-d+1} \cdot f^*(H)^{\cdot d-1} \cdot V) \\
+ (l-d+1) \deg(A \cdot H^{\cdot d-1}) \deg(L_{\eta}^{\cdot l-d} \cdot V_{\eta}).
\end{multline*}
Therefore, by Step~3, we get our theorem.
\QED

\subsection{Geometric height functions defined over a finitely generated field over $\FF_q$}
\label{subsec:geometric:height:function:F:q}
Let $K$ be a finitely generated field over $\FF_q$ with $d = \trdeg_{\FF_q}(K) \geq 1$.
Let $X$ be a projective variety over $K$ and $L$ a line bundle on $X$.
Here we fix a projective variety $B$ and a nef and big line bundle $H$ on $B$ such that
the function field of $B$ is $K$.
We choose a pair $(\mathcal{X}, \mathcal{L})$ with the following properties:
\begin{enumerate}
\renewcommand{\labelenumi}{(\arabic{enumi})}
\item
$\mathcal{X}$ is a projective variety over $\FF_q$ and there is a morphism
$f : \mathcal{X} \to B$ over $\FF_q$ such that $X$ is the generic fiber of $f$.

\item
$\mathcal{L}$ is a $\QQ$-line bundle on $\mathcal{X}$ (i.e., 
$\mathcal{L} \in \Pic(\mathcal{X}) \otimes \QQ$) such that
$\rest{\mathcal{L}}{X}$ coincides with $L$ in $\Pic(X) \otimes \QQ$.
\end{enumerate}
The pair $(\mathcal{X}, \mathcal{L})$ is called a {\em model of $(X, L)$}.

For $x \in X(\overline{K})$, let $\Delta_x$ be the closure of the image
$\Spec(\overline{K}) \overset{x}{\longrightarrow} X \hookrightarrow \mathcal{X}$.
Then, the height function of $(X, L)$ with respect to $(B, H)$ and $(\mathcal{X}, \mathcal{L})$
is defined by
\[
h^{(B,H)}_{(\mathcal{X}, \mathcal{L})}(x) = \frac{\deg(\mathcal{L} \cdot f^*(H)^{d-1} \cdot \Delta_x)}{[K(x):K]}.
\]
It is not difficult to see that
if $(\mathcal{X}', \mathcal{L}')$ is another model of $(X,L)$, then
there is a constant $C$ such that
\[
\vert h^{(B,H)}_{(\mathcal{X}, \mathcal{L})}(x) - h^{(B,H)}_{(\mathcal{X}', \mathcal{L}')}(x)
\vert \leq C
\]
for all $x \in X(\overline{K})$ (cf. \cite[the proof of Proposition~3.3.3]{MoArht}). 
Thus, the height function is uniquely determined modulo
bounded functions. In this sense, we denote the class of $h^{(B,H)}_{(\mathcal{X}, \mathcal{L})}$ 
modulo bounded functions by $h^{(B,H)}_{L}(x)$.
As a corollary of Theorem~\ref{thm:relative:case:finite:field}, we have the following.

\begin{Corollary}
\label{cor:refine:northcott:F:q}
Let $h_L$ be a representative of $h^{(B, H)}_L$.
If $L$ is ample, then, for a fixed $k$, there is a constant $C$ such that
\[
\{ x \in X(\overline{K}) \mid h_L(x) \leq h, [K(x):K] \leq k \} \leq q^{C \cdot h^d}
\]
for all $h \geq 1$.
\end{Corollary}

\Proof
Since $L$ is ample, we can find a model $(\mathcal{X}, \mathcal{L})$ of $(X,L)$ such that
$\mathcal{L}$ is nef (cf. Step~4 of Theorem~\ref{thm:relative:case:finite:field}). 
Thus, our assertion follows from Theorem~\ref{thm:relative:case:finite:field}.
\QED

\section{Preliminaries for the arithmetic case}

\subsection{Arakelov geometry}
In this paper, a flat and quasi-projective integral
scheme over $\ZZ$ is called an {\em arithmetic variety}.
If it is smooth over $\QQ$, then it is said to be {\em generically smooth}.

Let $X$ be a generically smooth arithmetic variety.
A pair $(Z, g)$ is called an {\em arithmetic cycle of codimension $p$}
if $Z$ is a cycle of codimension $p$ and
$g$ is a current of type $(p-1,p-1)$ on $X(\CC)$.
We denote by $\widehat{Z}^p(X)$ the set of all arithmetic cycles
on $X$. We set 
\[
\aCH^p(X) = \widehat{Z}^p(X)/\!\!\sim,
\]
where $\sim$ is the arithmetic linear equivalence.

Let $\overline{L} = (L, \Vert\cdot\Vert)$
be a $C^{\infty}$-hermitian line bundle on $X$.
Then, a homomorphism
\[
\acherncl_1(\overline{L}) \cdot :
\aCH^p(X) \to \aCH^{p+1}(X)
\]
is define by
\[
\acherncl_1(\overline{L}) \cdot (Z, g) =
\left(\text{$\zeros(s)$ on $Z$}, [-\log (\Vert s \Vert^2_Z)]
+ c_1(\overline{L}) \wedge g \right),
\]
where $s$ is a rational section of $\rest{L}{Z}$ and
$[-\log (\Vert s \Vert^2_Z)]$ is a current given by
$\phi \mapsto -\int_{Z(\CC)}\log(\Vert s \Vert^2_Z)\phi$.

Here we assume that $X$ is projective.
Then we can define the arithmetic degree map
\[
\adeg : \aCH^{\dim X}(X) \to \RR
\]
by
\[
\adeg \left( \sum_P n_P P, g \right) = \sum_{P} n_P \log(\#(\kappa(P))) +
\frac{1}{2} \int_{X(\CC)} g.
\]
Thus, if $C^{\infty}$-hermitian line bundles
$\overline{L}_1, \ldots, \overline{L}_{\dim X}$ are given, then
we can get the number
\[
\adeg \left(\acherncl_1(\overline{L}_1) \cdots 
\acherncl_1(\overline{L}_{\dim X}) \right),
\]
which is called the {\em arithmetic intersection number
of $\overline{L}_1, \ldots, \overline{L}_{\dim X}$}.

\medskip
Let $X$ be a projective arithmetic variety.
Note that $X$ is not necessarily generically smooth.
Let $\overline{L}_1, \ldots, \overline{L}_{\dim X}$ be
$C^{\infty}$-hermitian line bundles on $X$.
Choose a birational morphism
$\mu : Y \to X$ such that $Y$ is a generically smooth projective
arithmetic variety.
Then, we can see that the arithmetic intersection number
\[
\adeg \left(
\acherncl_1(\mu^*(\overline{L}_1)) \cdots
\acherncl_1(\mu^*(\overline{L}_{\dim X})) \right)
\]
does not depend on the choice of the generic
resolution of singularities $\mu : Y \to X$.
Thus, we denote this number by
\[
\adeg \left(
\acherncl_1(\overline{L}_1) \cdots
\acherncl_1(\overline{L}_{\dim X}) \right).
\]

\medskip
Let $f : X \to Y$ be a morphism of projective arithmetic varieties.
Let $\overline{L}_1, \ldots, \overline{L}_r$ be
$C^{\infty}$-hermitian line bundles on $X$, and
$\overline{M}_1, \ldots, \overline{M}_s$ 
$C^{\infty}$-hermitian line bundles on $Y$. If $r + s = \dim X$, then
the following formula is called the projection formula:
\renewcommand{\theequation}{\arabic{section}.\arabic{subsection}.\arabic{Theorem}}
\addtocounter{Theorem}{1}
\begin{multline}
\label{eqn:projection:formula:1}
\adeg \left( \acherncl_1(\overline{L}_1) \cdots \acherncl_1(\overline{L}_r)
\cdot \acherncl_1(f^*(\overline{M}_1)) \cdots \acherncl_1(f^*(\overline{M}_s))
\right) \\
= \begin{cases}
0 & \text{if $s > \dim Y$} \\
\deg((L_1)_{\eta} \cdots (L_r)_{\eta}) 
\adeg(\acherncl_1(\overline{M}_1) \cdots \acherncl_1(\overline{M}_s)) &
\text{if $s = \dim Y$ and $r > 0$} \\
\deg(f) \adeg(\acherncl_1(\overline{M}_1) \cdots \acherncl_1(\overline{M}_s))
& \text{if $s = \dim Y$ and $r = 0$},
\end{cases}
\end{multline}
\renewcommand{\theequation}{\arabic{section}.\arabic{subsection}.\arabic{Theorem}.\arabic{Claim}}
where the subscript $\eta$ means the restriction of line bundles to the generic
fiber of $f$.

\medskip
Let $\overline{L}_1, \ldots, \overline{L}_l$ be
$C^{\infty}$-hermitian line bundles on a projective arithmetic
variety $X$.
Let $V$ be an $l$-dimensional integral closed subscheme on $X$.
Then, 
$\adeg\left(\acherncl_1(\overline{L}_1) \cdots \acherncl_1(\overline{L}_l)
\crest V\right)$ is defined by
\[
\adeg\left(\acherncl_1(\rest{\overline{L}_1}{V}) \cdots 
\acherncl_1(\rest{\overline{L}_l}{V})\right).
\]
Note that if $V$ is lying over a prime $p$ with respect to $X \to \Spec(\ZZ)$, then
\[
\adeg\left(\acherncl_1(\overline{L}_1) \cdots 
\acherncl_1(\overline{L}_l) \crest V\right)
= \log(p) \deg(\rest{L_1}{V} \cdots \rest{L_l}{V}).
\]
Moreover, for an $l$-dimensional cycle
$Z = \sum_{i} n_i V_i$ on $X$,
$\sum_i 
n_i \adeg\left(\acherncl_1(\overline{L}_1) \cdots \acherncl_1(\overline{L}_l)
\crest V_i \right)$ is given by
\[
\adeg\left(\acherncl_1(\overline{L}_1) \cdots \acherncl_1(\overline{L}_l)
\crest Z \right).
\]
Let $f : X \to Y$ be a morphism of projective arithmetic varieties.
Let $\overline{M}_1, \ldots, \overline{M}_l$ be $C^{\infty}$-hermitian line
bundles on $Y$. Then, as a consequence of \eqref{eqn:projection:formula:1},
we have
\renewcommand{\theequation}{\arabic{section}.\arabic{subsection}.\arabic{Theorem}}
\addtocounter{Theorem}{1}
\begin{equation}
\label{eqn:projection:formula:2}
\adeg \left( \acherncl_1(f^*(\overline{M}_1)) \cdots
\acherncl_1(f^*(\overline{M}_l)) \crest Z \right) =
\adeg \left( \acherncl_1(\overline{M}_1) \cdots
\acherncl_1(\overline{M}_l) \crest f_*(Z) \right)
\end{equation}
\renewcommand{\theequation}{\arabic{section}.\arabic{subsection}.\arabic{Theorem}.\arabic{Claim}}
for all $l$-dimensional cycles $Z$ on $X$.

\subsection{The positivity of $C^{\infty}$-hermitian $\QQ$-line bundles}
Let $X$ be a projective arithmetic variety and
$\overline{L}$ a $C^{\infty}$-hermitian $\QQ$-line bundle on $X$.
Let us consider several kinds of
the positivity of $C^{\infty}$-hermitian $\QQ$-line bundles.

$\bullet${\bf ample}:
We say $\overline{L}$ is {\em ample} if
$L$ is ample on $X$, $c_1(\overline{L})$ is positive form on $X(\CC)$, and
there is a positive number $n$ such that
$L^{\otimes n}$ is generated by
the set $\{ s \in H^0(X, L^{\otimes n}) \mid \Vert s \Vert_{\sup} < 1 \}$.

$\bullet${\bf nef}:
We say  $\overline{L}$ is {\em nef} if
$c_1(\overline{L})$ is a semipositive form on $X(\CC)$ and,
for all one-dimensional integral closed subschemes $\Gamma$ of $X$,
$\adeg \left( \acherncl_1(\overline{L}) \crest \Gamma \right) \geq 0$.

$\bullet${\bf big}:
$\overline{L}$ is said to be {\em big} if 
$\rank_{\ZZ} H^0(X, L^{\otimes m}) = O(m^{\dim X_{\QQ}})$ and there is a non-zero
section $s$ of $H^0(X, L^{\otimes n})$ with $\Vert s \Vert_{\sup} < 1$ for
some positive integer $n$.

$\bullet${\bf $\pmb{\QQ}$-effective}:
$\overline{L}$
is said to be {\em $\QQ$-effective} if there is a positive integer $n$ and
a non-zero $s \in H^0(X, L^{\otimes n})$ with $\Vert s \Vert_{\sup} \leq 1$.

$\bullet${\bf pseudo-effective}:
$\overline{L}$
is said to be {\em pseudo-effective} if there are
(1) a sequence $\{ \overline{L}_n\}_{n=1}^{\infty}$ of $\QQ$-effective 
$C^{\infty}$-hermitian $\QQ$-line bundles,
(2) $C^{\infty}$-hermitian $\QQ$-line bundles 
$\overline{E}_1, \ldots, \overline{E}_r$ and
(3) sequences $\{ a_{1, n} \}_{n=1}^{\infty}, \ldots, \{a_{r, n}\}_{n=1}^{\infty}$
of rational numbers such that
\[
\acherncl_1(\overline{L}) = \acherncl_1(\overline{L}_n) + 
\sum_{i=1}^r a_{i,n}\acherncl_1(\overline{E}_i)
\]
in $\aCH(X) \otimes \QQ$ and
$\lim_{n\to\infty} a_{i, n} = 0$ for all $i$.
If $\overline{L}_1 \otimes \overline{L}_2^{\otimes -1}$ is pseudo-effective
for $C^{\infty}$-hermitian $\QQ$-line bundles $\overline{L}_1, \overline{L}_2$
on $X$, then we denote this by $\overline{L}_1 \succsim \overline{L}_2$.

$\bullet${\bf of surface type}:
$\overline{L}$
is said to be {\em of surface type} if there are
a morphism $\phi : X \to X'$ of
projective arithmetic varieties and
a $C^{\infty}$-hermitian $\QQ$-line bundle $\overline{L}'$ on $X'$
such that $\dim X'_{\QQ} = 1$ (i.e. $X'$ is a projective arithmetic surface),
$\overline{L}'$ is nef and big, and that
$\phi^*(\overline{L}') = \overline{L}$ in $\aPic(X) \otimes \QQ$.

\medskip
Here let us consider three lemmas which will be used later.

\begin{Lemma}
\label{lem:comp:nef:eff:div}
Let $X$ be a projective arithmetic variety. Then, we have the following.
\begin{enumerate}
\renewcommand{\labelenumi}{(\arabic{enumi})}
\item
Let $\overline{L}_1, \ldots, \overline{L}_{\dim X},
\overline{M}_1, \ldots, \overline{M}_{\dim X}$ be
nef $C^{\infty}$-hermitian $\QQ$-line bundles on $X$.
If $\overline{L}_i \otimes \overline{M}_i^{\otimes -1}$ is
pseudo-effective for every $i$, then
\[
\adeg\left( \acherncl_1(\overline{L}_1) \cdots 
\acherncl_1(\overline{L}_{\dim X})
\right) \geq
\adeg\left( \acherncl_1(\overline{M}_1) \cdots 
\acherncl_1(\overline{M}_{\dim X})
\right).
\]

\item
Let $V$ be an effective cycle of dimension $l$ and
let $\overline{L}_1, \ldots, \overline{L}_{l},
\overline{M}_1, \ldots, \overline{M}_{l}$ be
nef $C^{\infty}$-hermitian $\QQ$-line bundles on $X$ such that,
for each $i$, there is a non-zero global section
$s_i \in H^0(X, L_i \otimes M_i^{\otimes -1})$ with
$\Vert s_i \Vert_{\rm sup} \leq 1$.
Let $V = \sum_{j} a_j V_j$ be the irreducible decomposition as a cycle.
If $\rest{s_i}{V_j} \not= 0$ for all $i,j$, then
\[
\adeg\left( \acherncl_1(\overline{L}_1) \cdots 
\acherncl_1(\overline{L}_{l})  \crest V
\right) \geq
\adeg\left( \acherncl_1(\overline{M}_1) \cdots 
\acherncl_1(\overline{M}_{l}) \crest V
\right).
\]
\end{enumerate}
\end{Lemma}

\Proof
(1)
This lemma follows from \cite[Proposition~2.3]{MoArht} and the following formula:
\begin{multline*}
\adeg\left( \acherncl_1(\overline{L}_1) \cdots 
\acherncl_1(\overline{L}_{\dim X})
\right) =
\adeg\left( \acherncl_1(\overline{M}_1) \cdots 
\acherncl_1(\overline{M}_{\dim X})
\right) + \\
\sum_{i=1}^{\dim X} \adeg\left( \acherncl_1(\overline{M}_1) \cdots 
\acherncl_1(\overline{M}_{i-1}) \cdot 
\acherncl_1(\overline{L}_i \otimes \overline{M}_i^{\otimes -1}) \cdot
\acherncl_1(\overline{L}_{i+1}) \cdots
\acherncl_1(\overline{L}_{\dim X})
\right).
\end{multline*}

\medskip
(2)
This is a consequence of (1).
\QED

\begin{Lemma}
\label{lem:power:sum:surface:type}
Let $X$ be a projective arithmetic variety and
$d$ an integer with $1 \leq d \leq \dim X$.
Let $X_1, \ldots, X_d$ be projective arithmetic surfaces
\rom{(}i.e. $2$-dimensional projective arithmetic varieties\rom{)} and 
$\phi_i : X \to X_i$ \rom{(}$i=1, \ldots, d$\rom{)} surjective morphisms.
Let $\overline{L}_1, \ldots, \overline{L}_d$ be 
$C^{\infty}$-hermitian $\QQ$-line bundles on $X_1, \ldots, X_d$
respectively with
$\deg((L_i)_{\QQ}) > 0$ \rom{(}$i=1, \ldots, d$\rom{)}, and
let $\overline{H}_{d+1}, \ldots, \overline{H}_{\dim X}$ be
$C^{\infty}$-hermitian $\QQ$-line bundles on $X$.
We set $\overline{H}_i = \phi_i^*(\overline{L}_i)$ and
$\overline{H} =  \bigotimes_{i=1}^d \overline{H}_i$. Then,
\begin{multline*}
\adeg\left( \acherncl_1(\overline{H})^{\cdot d} \cdot \acherncl_1(\overline{H}_{d+1})
\cdots \acherncl_1(\overline{H}_{\dim X}) \right) =
d! \adeg\left( \acherncl_1(\overline{H}_1) \cdots
\acherncl_1(\overline{H}_d) \cdot \acherncl_1(\overline{H}_{d+1})
\cdots \acherncl_1(\overline{H}_{\dim X}) \right) \\
+ \frac{d!}{2} \sum_{i\not= j}
\frac{\adeg\left( \acherncl_1(\overline{L}_i)^{\cdot 2} \right)
\deg\left( \prod_{\substack{ l \not = j \\ 1 \leq l \leq d}}
(H_l)_{\QQ} \cdot (H_{d+1})_{\QQ} \cdots (H_{\dim X})_{\QQ} \right)}{
\deg((L_i)_{\QQ})}.
\end{multline*}
\end{Lemma}

\Proof
First of all,
\begin{multline*}
\adeg\left( \acherncl_1(\overline{H})^{\cdot d} \cdot \acherncl_1(\overline{H}_{d+1})
\cdots \acherncl_1(\overline{H}_{\dim X}) \right) \\
= \sum_{\substack{a_1 + \cdots + a_d = d \\ a_1 \geq 0, \ldots, a_d \geq 0}}
\frac{d!}{a_1! \cdots a_d!}
\adeg\left( \acherncl_1(\overline{H}_1)^{\cdot a_1} \cdots 
\acherncl_1(\overline{H}_{d})^{\cdot a_{d}} 
\cdot \acherncl_1(\overline{H}_{d+1})
\cdots \acherncl_1(\overline{H}_{\dim X}) \right).
\end{multline*}

\begin{Claim}
If $(a_1, \ldots, a_d) \not= (1, \ldots, 1)$ and
\[
\adeg\left(\prod_{l=1}^d 
\acherncl_1(\overline{H}_l)^{\cdot a_l} \cdot
\cdot \acherncl_1(\overline{H}_{d+1})
\cdots \acherncl_1(\overline{H}_{\dim X})
\right) \not= 0,
\]
then there are $i, j \in \{1, \ldots, d\}$ such that
$a_i = 2$, $a_j= 0$ and $a_l = 1$ for all $l \not= i, j$.
\end{Claim}

Clearly, $a_l \leq 2$ for all $l$. Thus, there is $i$ with $a_i = 2$.
Suppose that $a_j = 2$ for some $j \not= i$.
Then,
\begin{multline*}
\adeg\left(\prod_{l=1}^d 
\acherncl_1(\overline{H}_l)^{\cdot a_l} \cdot
\cdot \acherncl_1(\overline{H}_{d+1})
\cdots \acherncl_1(\overline{H}_{\dim X})
\right) \\
= \adeg\left(\acherncl_1(\phi_i^*(\overline{L}_i))^{\cdot 2} \cdot
\acherncl_1(\phi_j^*(\overline{L}_j))^{\cdot 2} \cdot
\prod_{l=1, l\not=i,j}^d \acherncl_1(\overline{H}_l)^{\cdot a_l}
\cdot \acherncl_1(\overline{H}_{d+1})
\cdots \acherncl_1(\overline{H}_{\dim X})
\right).
\end{multline*}
Thus, using the projection formula with respect to $\phi_i$,
\begin{multline*}
\adeg\left(\acherncl_1(\phi_i^*(\overline{L}_i))^{\cdot 2} \cdot
\acherncl_1(\phi_j^*(\overline{L}_j))^{\cdot 2} \cdot
\prod_{l=1, l\not=i,j}^d \acherncl_1(\overline{H}_l)^{\cdot a_l}
\cdot \acherncl_1(\overline{H}_{d+1})
\cdots \acherncl_1(\overline{H}_{\dim X})
\right) \\
= \adeg(\acherncl_1(\overline{L}_i)^{\cdot 2})
\deg \left(
\phi_j^*(L_j)_{\eta_i}^{\cdot 2} \cdot 
\prod_{l=1, l\not=i,j}^d (H_l)_{\eta_i}^{\cdot a_l}
\cdot (H_{d+1})_{\eta_i}
\cdots (H_{\dim X})_{\eta_i}\right),
\end{multline*}
where $\eta_i$ means the restriction of line bundles to the generic fiber of $\phi_i$.
Here $(X_j)_{\QQ}$ is projective curve.
Thus, we can see
\[
\deg \left(
\phi_j^*(L_j)_{\eta_i}^{\cdot 2} \cdot 
\prod_{l=1, l\not=i,j}^d (H_l)_{\eta_i}^{\cdot a_l}
\cdot (H_{d+1})_{\eta_i}
\cdots (H_{\dim X})_{\eta_i}\right) = 0.
\]
This is a contradiction.
Hence, we get our claim.

\medskip
By the above claim, it is sufficient to see that
\begin{multline*}
\adeg\left( \acherncl_1(\phi_i^*(\overline{L}_i))^{\cdot 2} 
\cdot \prod_{l=1,l\not=i,j}^d 
\acherncl_1(\overline{H}_l) \cdot
\acherncl_1(\overline{H}_{d+1})
\cdots \acherncl_1(\overline{H}_{\dim X})
\right) \\
= \frac{\adeg\left( \acherncl_1(\overline{L}_i)^{\cdot 2} \right)
\deg\left( \prod_{\substack{ l \not = j \\ 1 \leq l \leq d}}
(H_l)_{\QQ} \cdot (H_{d+1})_{\QQ} \cdots (H_{\dim X})_{\QQ} \right)}{
\deg((L_i)_{\QQ})}
\end{multline*}

By the projection formula with respect to $\phi_i$,
\begin{multline*}
\adeg\left( \acherncl_1(\phi_i^*(\overline{L}_i))^{\cdot 2} 
\cdot \prod_{l=1,l\not=i,j}^d 
\acherncl_1(\overline{H}_l) \cdot
\acherncl_1(\overline{H}_{d+1})
\cdots \acherncl_1(\overline{H}_{\dim X})
\right) \\
= \adeg\left( \acherncl_1(\overline{L}_i)^{\cdot 2} \right)
\deg\left( \prod_{\substack{ l \not =i, j \\ 1 \leq l \leq d}}
(H_l)_{\eta_i} \cdot (H_{d+1})_{\eta_i} \cdots (H_{\dim X})_{\eta_i} \right).
\end{multline*}
Moreover, using the projection formula with respect to $\phi_i$ again,
\begin{multline*}
\deg\left( \prod_{\substack{ l \not = j \\ 1 \leq l \leq d}}
(H_l)_{\QQ} \cdot (H_{d+1})_{\QQ} \cdots (H_{\dim X})_{\QQ} \right) \\
= \deg((L_i)_{\QQ})
\deg\left( \prod_{\substack{ l \not =i, j \\ 1 \leq l \leq d}}
(H_l)_{\eta_i} \cdot (H_{d+1})_{\eta_i} \cdots (H_{\dim X})_{\eta_i} \right).
\end{multline*}
Thus, we get our lemma.
\QED

Finally let us consider the following technical lemma.

\begin{Lemma}
\label{lem:adeg:comp:birat:product}
Let $\phi : \PP^{n}_{\ZZ} \dasharrow (\PP^1_{\ZZ})^{n}$ be the birational map
given by
\[
(X_0 : \ldots : X_{n_i}) \mapsto (X_0 : X_1) \times \cdots \times (X_0: X_{n}).
\]
Let $\Sigma$ be the boundary of $\PP^{n}_{\ZZ}$, that is, 
$\Sigma = \{ X_0 = 0 \}$.
Let $B$ be a projective arithmetic variety and
$\overline{H}_1, \ldots, \overline{H}_d$ nef $C^{\infty}$-hermitian line bundles on $B$,
where $d = \dim B_{\QQ}$.
For $V \in Z^{\rm eff}_l(\PP^{n}_{\ZZ} \times B  ;
(\PP^{n}_{\ZZ} \setminus \Sigma) \times B)$,
we denote by $V'$ the strict transform of $V$ by $\phi \times \operatorname{id}
: \PP^n_{\ZZ} \times B \dasharrow (\PP^1_{\ZZ})^n \times B$
\rom{(}For the definition of 
$Z^{\rm eff}_l(\PP^{n}_{\ZZ} \times B  ;
(\PP^{n}_{\ZZ} \setminus \Sigma) \times B)$,
see \rom{\ref{subsec:notation:cycle}}\rom{)}.
Let us fix a non-negative real number $\lambda$.
Then,
\begin{multline*}
n^{l-d} \adeg\left(
\acherncl_1(p^*(\overline{\OO}^{{\rm FS}_{\lambda}}(1)))^{\cdot l-d} \cdot
\acherncl_1(q^*(\overline{H}_1)) \cdots  \acherncl_1(q^*(\overline{H}_d)) \crest V
\right) \\
\geq \adeg\left(
\acherncl_1({p'}^*(\overline{\OO}^{{\rm FS}_{\lambda}}(1,\ldots,1)))^{\cdot l-d} \cdot
\acherncl_1({q'}^*(\overline{H}_1)) \cdots \acherncl_1({q'}^*(\overline{H}_d)) \crest V'
\right)
\end{multline*}
for all $V \in Z^{\rm eff}_l(\PP^{n}_{\ZZ} \times B  ;
(\PP^{n}_{\ZZ} \setminus \Sigma) \times B)$,
where $p : \PP^n_{\ZZ} \times B \to \PP^n_{\ZZ}$ and
$p' : (\PP^1_{\ZZ})^n \times B \to (\PP^1_{\ZZ})^n$
\rom{(}resp. $q : \PP^n_{\ZZ} \times B \to B$ and
$q' : (\PP^1_{\ZZ})^n \times B \to B$\rom{)}
are the projections to the first factor \rom{(}the second factor\rom{)}.
Note that in the case $d=0$, we do not use the nef $C^{\infty}$-hermitian line bundles
$\overline{H}_1, \ldots, \overline{H}_d$.
\end{Lemma}

\Proof
Let $Y$ ($\subseteq \PP^n_{\ZZ} \times (\PP^1_{\ZZ})^n)$ be the graph of the rational map
$\phi : \PP^{n}_{\ZZ} \dasharrow (\PP^1_{\ZZ})^{n}$.
Let $\mu : Y \to \PP^{n}_{\ZZ}$ and
$\nu : Y \to  (\PP^1_{\ZZ})^{n}$ be the morphisms induced by
the projections.
Here we claim the following:

\begin{Claim}
There are an effective Cartier divisor $E$ on $Y$,
a non-zero section $s \in H^0(Y, \OO_{Y}(E))$ and
a $C^{\infty}$-metric $\Vert\cdot\Vert_{E}$ of $\OO_{Y}(E)$ such that
\begin{enumerate}
\renewcommand{\labelenumi}{(\arabic{enumi})}
\item
$\zeros(s) = E$, $\mu(E) \subseteq \Sigma$, 

\item
${\mu}^*(\overline{\OO}^{{\rm FS}_{\lambda}}(1))^{\otimes n} =
{\nu}^*(\overline{\OO}^{{\rm FS}_{\lambda}}(1,\ldots,1))
\otimes (\OO_{Y}(E), \Vert\cdot\Vert_{E})$, and that

\item
$\Vert s \Vert_{E}(x) \leq 1$ for all $x \in Y(\CC)$.
\end{enumerate}
\end{Claim}

\medskip
Let $Y_i$ ($\subseteq \PP^n_{\ZZ} \times \PP^1_{\ZZ})$ be the graph
of the rational map
$\PP^n_{\ZZ} \dasharrow \PP^1_{\ZZ}$ given by
\[
(X_0 : \cdots : X_n) \mapsto (X_0 : X_i).
\]
Let $\mu_i : Y_i \to \PP^n_{\ZZ}$ and $\nu_i : Y_i \to \PP^1_{\ZZ}$
be the morphisms induced by
the projections $\PP^n_{\ZZ} \times \PP^1_{\ZZ} \to \PP^n_{\ZZ}$
and $\PP^n_{\ZZ} \times \PP^1_{\ZZ} \to \PP^1_{\ZZ}$
respectively.
Let $\pi_i : (\PP^1_{\ZZ})^n \to \PP^1_{\ZZ}$ be the projection
to the $i$-th factor.
Moreover, let $h_i : Y \to Y_i$ be the morphism induced by
$\operatorname{id} \times \pi_i : 
\PP^n_{\ZZ} \times (\PP^1_{\ZZ})^n \to \PP^n_{\ZZ} \times \PP^1_{\ZZ}$.
Consequently, we have the following commutative diagram:
\[
\xymatrix{
& Y \ar[ddl]_{\mu} \ar[d]^{h_i} \ar[dr]^{\nu} \\
& Y_i \ar[dl]^{\mu_i} \ar[dr]_{\nu_i} & (\PP^1_{\ZZ})^n \ar[d]^{\pi_i} \\
\PP^n_{\ZZ} & & \PP^1_{\ZZ} 
}
\]
Note that $Y_i$ is the blowing-up by the ideal sheaf $I_i$ generated
by $X_0$ and $X_i$.
Thus there is an effective Cartier divisor $E_i$ on $Y_i$
with $I_i \OO_{Y_i} = \OO_{Y_i}(-E_i)$ and
$\mu_i^*(\OO_{\PP^n}(1)) \otimes \OO_{Y_i}(-E_i) = \nu_i^*(\OO_{\PP^1}(1))$. 
Let $s_i$ be the canonical section of $\OO_{Y_i}(E_i)$.
We choose $C^{\infty}$-metric $\Vert\cdot\Vert_i$ of $\OO_{Y_i}(E_i)$
with
\[
\mu_i^*(\OO_{\PP^n}(1), \Vert\cdot\Vert_{{\rm FS}_{\lambda}}) = 
\nu_i^*(\OO_{\PP^1}(1),\Vert\cdot\Vert_{{\rm FS}_{\lambda}}) \otimes
(\OO_{Y_i}(E_i), \Vert\cdot\Vert_i).
\]
Let $(T_0 : T_1)$ be a coordinate of $\PP^1_{\ZZ}$.
Then, $\mu_i^*(X_0) = \nu_i^*(T_0) \otimes s_i$.
Thus,
\[
\frac{\exp(-\lambda)\vert X_0 \vert}{\sqrt{\vert X_0 \vert^2 + \cdots + \vert X_n \vert^2}} =
\frac{\exp(-\lambda)\vert T_0 \vert}{\sqrt{\vert T_0 \vert^2 + \vert T_1 \vert^2}}
\Vert s_i \Vert_i,
\]
which implies
\[
\Vert s_i \Vert_i = \frac{\sqrt{\vert X_0 \vert^2 + \vert X_1 \vert^2}}
{\sqrt{\vert X_0 \vert^2 + \cdots + \vert X_n \vert^2}}
\]
because $X_0 T_1 = X_i T_0$. Therefore, $\Vert s_i \Vert_i(x_i) \leq 1$
for all $x_i \in Y_i(\CC)$.
We set $E = \sum_{i=1}^n h_i^*(E_i)$ and give a $C^{\infty}$-metric 
$\Vert\cdot\Vert_E$ to $\OO_Y(E)$ with
\[
(\OO_Y(E), \Vert\cdot\Vert_E) =
\bigotimes_{i=1}^n h_i^*(\OO_{Y_i}(E_i), \Vert\cdot\Vert_i).
\]
Thus, if we set $s = h_1^*(s_1) \otimes \cdots \otimes h_n^*(s_n)$, then
$s \in H^0(Y, \OO_Y(E))$,
$\zeros(s) = E$ and $\Vert s \Vert_E(x) \leq 1$ for all $x \in Y(\CC)$.
Moreover, we have
\[
\mu^*(\overline{\OO}^{{\rm FS}_{\lambda}}(1))^{\otimes n}  = 
\nu^*(\overline{\OO}^{{\rm FS}_{\lambda}}(1, \ldots, 1)) \otimes (\OO_Y(E), \Vert\cdot\Vert_E).
\]
Hence we get our claim.

\bigskip
For $V \in Z^{\rm eff}_l(\PP^{n}_{\ZZ} \times B  ;
(\PP^{n}_{\ZZ} \setminus \Sigma) \times B)$, 
let $V''$ be the strict transform of $V$ by $\mu \times \operatorname{id}
: Y \times B \to \PP^n_{\ZZ} \times B$.
Let $p'' : Y \times B \to Y$ and $q'' : Y \times B \to B$
be the projections to the first factor and the second factors respectively.
Then, by using the projection formula,
\begin{multline*}
\adeg\left(
\acherncl_1(p^*(\overline{\OO}^{{\rm FS}_{\lambda}}(n)))^{\cdot l-d} \cdot
\acherncl_1(q^*(\overline{H}_1)) \cdots  \acherncl_1(q^*(\overline{H}_d)) \crest V
\right) \\
= \adeg\left(
\acherncl_1({p''}^*(\mu^*(\overline{\OO}^{{\rm FS}_{\lambda}}(n))))^{\cdot l-d}
\cdot \acherncl_1({q''}^*(\overline{H}_1)) \cdots  \acherncl_1({q''}^*(\overline{H}_d))
\crest V''
\right).
\end{multline*}
Moreover, by virtue of (2) of Lemma~\ref{lem:comp:nef:eff:div},
\begin{multline*}
\adeg\left(
\acherncl_1({p''}^*(\mu^*(\overline{\OO}^{{\rm FS}_{\lambda}}(n))))^{\cdot l-d}
\cdot \acherncl_1({q''}^*(\overline{H}_1)) \cdots  \acherncl_1({q''}^*(\overline{H}_d))
\crest V''
\right) \\
\geq 
\adeg\left(
\acherncl_1({p''}^*(\nu^*(\overline{\OO}^{{\rm FS}_{\lambda}}(1,\ldots,1))))^{\cdot l-d}
\cdot \acherncl_1({q''}^*(\overline{H}_1)) \cdots  \acherncl_1({q''}^*(\overline{H}_d))
\crest V''
\right)
\end{multline*}
Thus, using the projection formula for $\nu \times \operatorname{id} :
Y \times B \to (\PP^1_{\ZZ})^n \times B$,
we get our lemma because $(\nu \times \operatorname{id})_*(V'') = V'$.
\QED

\subsection{Polarization of a finitely generated field over $\QQ$}
Let $K$ be a finitely generated field over $\QQ$ with
$d = \trdeg_{\QQ}(K)$, and let
$B$ be a projective arithmetic variety such that
$K$ is the function field of $B$. 
Here we fix several notations.

$\bullet${\bf polarization}:
A collection $\overline{B} = (B; \overline{H}_1, \ldots, \overline{H}_d)$
of $B$ and nef $C^{\infty}$-hermitian $\QQ$-line bundles 
$\overline{H}_1, \ldots,
\overline{H}_d$ on $B$ is called a {\em polarization of $K$}.

$\bullet${\bf big polarization}:
A polarization $\overline{B} = (B; \overline{H}_1, \ldots, \overline{H}_d)$
is said to be {\em big}
if $\overline{H}_1, \ldots, \overline{H}_d$ are nef and big.

$\bullet${\bf fine polarization}:
A polarization $\overline{B} = (B; \overline{H}_1, \ldots, \overline{H}_d)$
is said to be {\em fine}
if there are a generically finite morphism
$\mu : B' \to B$ of
projective arithmetic varieties, and
$C^{\infty}$-hermitian $\QQ$-line bundles
$\overline{L}_1, \ldots, \overline{L}_d$
on $B'$ such that
$\overline{L}_1, \ldots, \overline{L}_d$ are of surface type,
$\mu^*(\overline{H}_i) \succsim \overline{L}_i$ for all $i$, and that
$\overline{L}_1 \otimes \cdots \otimes \overline{L}_d$ is nef and big.

\medskip
Let us consider the following proposition.

\begin{Proposition}
\label{prop:fine:pol:PP:1}
If a polarization $\overline{B} = (B; \overline{H}_1, \ldots, \overline{H}_d)$
is fine, then there are generically finite morphisms
$\mu : B' \to B$ and $\nu : B' \to (\PP^1_{\ZZ})^{d}$
with the following property:
for any real number $\lambda$,
there are positive rational numbers
$a_1, \ldots, a_{d}$ such that
\[
\mu^*(\overline{H}_i) \succsim \nu^*(q_i^*(\overline{\OO}^{{\rm FS}_{\lambda}}
(1)))^{\otimes a_i}
\]
for all $i = 1, \ldots, d$,
where $q_i : (\PP^1_{\ZZ})^{d} \to \PP^1_{\ZZ}$
is the projection to the $i$-th factor.
\end{Proposition}

\Proof
By the definition of fineness,
there are a
generically finite morphism
$\mu : B' \to B$ of
projective arithmetic varieties,
morphisms $\phi_i : B' \to B_i$ ($i=1, \ldots, d$)
of projective arithmetic varieties, and
nef and big $C^{\infty}$-hermitian $\QQ$-line bundles
$\overline{Q}_i$ on $B_i$ ($i=1, \ldots, d$) such that
$B_i$'s are arithmetic surfaces,
$\mu^*(\overline{H}_i) \succsim \phi_i^*(\overline{Q}_i)$
for all $i$, and that
$\phi_1^*(\overline{Q}_1) \otimes \cdots \otimes 
\phi_d^*(\overline{Q}_d)$ is nef and big.
Here, there are dominant rational maps
$\psi_i : B_i \dasharrow \PP^1_{\ZZ}$ for $i=1, \ldots, d$.
Replacing $B'$ and $B_i$'s by their suitable birational models,
we may assume $\psi_i$'s are morphisms.
Let $\nu : B' \to (\PP^1_{\ZZ})^d$ be a morphism
given by $\nu(x) = (\psi_1(\phi_1(x)), \ldots, \psi_d(\phi_d(x)))$.
Let us fix a real number $\lambda$.
Then, since $\overline{Q}_i$ is nef and big,
there is a positive rational number $a_i$ with
$\overline{Q}_i \succsim 
\psi_i^*(\overline{\OO}^{{\rm FS}_{\lambda}}(1))^{\otimes a_i}$.
Thus, 
\[
\mu^*(\overline{H}_i) \succsim \phi_i^*(\overline{Q}_i) \succsim
\phi_i^*(\psi_i^*(\overline{\OO}^{{\rm FS}_{\lambda}}(1)))^{\otimes a_i}
= \nu^*(q_i^*(\overline{\OO}^{{\rm FS}_{\lambda}}(1)))^{\otimes a_i}.
\]
Finally, we need to see that $\nu$ is generically finite.
For this purpose, it is sufficient to see that
$\nu^*(q_1^*(\OO(1)) \otimes \cdots \otimes q_d^*(\OO(1)))$ is nef and big
on $B'_{\QQ}$.
Indeed, we can find a positive rational number $a$
such that $\psi_i^*(\OO(1)) \otimes Q_i^{\otimes -a}$ is ample
over $B'_{\QQ}$ for all $i$. Thus,
\[
\bigotimes_{i=1}^d
\phi_i^*(\psi_i^*(\OO(1)) \otimes Q_i^{\otimes -a}) = 
\nu^*\left(\bigotimes_{i=1}^d q_i^*(\OO(1))\right)\otimes
\left(\bigotimes_{i=1}^d \phi_i^*(Q_i)\right)^{\otimes -a}
\]
is semiample on $B'_{\QQ}$.
Thus, $\nu^*(q_1^*(\OO(1)) \otimes \cdots \otimes q_d^*(\OO(1)))$
is nef and big because
$\phi_1(Q_1) \otimes \cdots \otimes \phi_d^*(Q_d)$ is
nef and big.
\QED

\medskip
Finally we would like to give a simple sufficient condition for the fineness
of a polarization.
Let $k$ be a number field, and $O_k$ the ring of integer in $k$.
Let $B_1, \ldots, B_l$ be projective and flat integral schemes over $O_k$ whose
generic fibers over $O_k$ are geometrically irreducible.
Let $K_i$ be the function field of $B_i$ and $d_i$ 
the transcendence degree of $K_i$
over $k$. We set $B = B_1 \times_{O_k} \cdots \times_{O_k} B_l$ and
$d = d_1 + \cdots + d_l$. Then, the function field of $B$ is
the quotient field of $K_1 \otimes_k K_2 \otimes_k \cdots \otimes_k K_l$,
which is denoted by $K$, and
the transcendence degree of $K$ over $k$ is $d$.
For each $i$ ($i = 1, \ldots, l$),
let $\overline{H}_{i, 1}, \ldots, \overline{H}_{i, d_i}$
be nef and big $C^{\infty}$-hermitian $\QQ$-line bundles on $B_i$.
We denote by $q_i$ the projection
$B \to B_i$ to the $i$-th factor.
Then, we have the following.

\begin{Proposition}
\label{prop:large:polarization}
A polarization $\overline{B}$ of $K$ given by
\[
\overline{B} = \left(B; q_1^*(\overline{H}_{1,1}), \ldots, 
q_1^*(\overline{H}_{1, d_1}),
\ldots, q_l^*(\overline{H}_{l,1}), \ldots, q_l^*(\overline{H}_{l, d_l})
\right)
\]
is fine.
In particular, a big polarization is fine.
\end{Proposition}

\Proof
Since there is a dominant rational map 
$B_i \dashrightarrow \left( \PP^1_{\ZZ} \right)^{d_i}$
by virtue of Noether's normalization theorem,
we can find a birational morphism
$\mu_i : B'_i \to B_i$ of projective integral schemes over $O_k$ and
a generically finite morphism $\nu_i : B'_i \to  
\left( \PP^1_{\ZZ} \right)^{d_i}$.
We set $B' = B'_1 \times_{O_k} \cdots \times_{O_k} B'_l$,
$\mu = \mu_1 \times \cdots \times \mu_l$ and 
$\nu = \nu_1 \times \cdots \times \nu_l$.
Let $\overline{L}$ be a $C^{\infty}$-hermitian line bundle on 
$\PP^1_{\ZZ}$ given by
$(\OO_{\PP^1_{\ZZ}}(1), \Vert\cdot\Vert_{FS})$.
Note that $\overline{L}$ is nef and big.
Then, since $\mu_i^*(\overline{H}_{i,j})$ is big,
there is a positive integer $a_{i,j}$ with
$\mu_i^*(\overline{H}_{i,j})^{\otimes a_{i,j}} \succsim 
\nu_i^*\left(  p_j^*(\overline{L}) \right)$
(cf. \cite[Proposition~2.2]{MoArht}), that is,
$\mu_i^*(\overline{H}_{i,j}) \succsim 
\nu_i^*\left(  p_j^*\left(\overline{L}^{\otimes 1/a_{i,j}}\right) \right)$.
Thus, we get our proposition.
\QED

\subsection{Height functions over a finitely generated field}
\label{subsec:height:functions:over:fgf}
\setcounter{Theorem}{0}
Let $K$ be a finitely generated field over $\QQ$ with
$d = \trdeg_{\QQ}(K)$, and let
$\overline{B} = (B; \overline{H}_1, \ldots, \overline{H}_d)$
be a polarization of $K$.
Let $X$ be a geometrically irreducible
projective variety over $K$ and $L$ an ample line bundle on $X$.
Let us take a projective integral scheme $\XX$ over $B$ and
a $C^{\infty}$-hermitian $\QQ$-line bundle $\overline{\LL}$ on $\XX$ such that
$X$ is the generic fiber of $\XX \to B$ and $L$ is equal to $\LL_K$
in $\Pic(X) \otimes \QQ$. 
The pair $(\XX, \overline{\LL})$ is called a model of $(X, L)$.
Then, for $x \in X(\overline{K})$,
we define $h^{\overline{B}}_{(\XX, \LL)}(x)$ to be
\[
h^{\overline{B}}_{(\XX, \overline{\LL})}(x) =
\frac{\adeg \left( \acherncl_1(\overline{\LL}) \cdot
\prod_{j=1}^d \acherncl_1(f^*(\overline{H}_j)) \crest \Delta_x \right)}
{[K(x) : K]},
\]
where $\Delta_x$ is the Zariski closure in $\XX$ of the image
$\Spec(\overline{K}) \to X \hookrightarrow \XX$, and
$f : \XX \to B$ is the canonical morphism.
By virtue of \cite[Corollary~3.3.5]{MoArht},
if $(\XX', \LL')$ is another model of $(X, L)$ over $B$, then
there is a constant $C$ with
$\vert h^{\overline{B}}_{(\XX, \LL)}(x) - 
h^{\overline{B}}_{(\XX', \LL')}(x) \vert \leq C$
for all $x \in X(\overline{K})$.
Hence, we have the unique height function $h^{\overline{B}}_L$
modulo the set of bounded functions.
In the case where $X = \PP^n_K$, if we set \begin{multline*}
h_{nv}^{\overline{B}}(x) = \sum_{\substack{\text{$\Gamma$ is a prime} \\ \text{divisor on $B$}}} 
\max_i \{ -\ord_{\Gamma}(\phi_i) \}
\adeg\left(  \acherncl_1(\overline{H}_1) \cdots \acherncl_1(\overline{H}_d)
 \crest \Gamma \right) \\
+ \int_{B(\CC)} \log \left( \max_i \{ \vert \phi_i \vert \} \right)
c_1(\overline{H}_1) \wedge \cdots \wedge c_1(\overline{H}_d)
\end{multline*}
for $x = (\phi_0 : \cdots : \phi_n) \in \PP^n(K)$,
then $h_{\OO(1)}^{\overline{B}} = h_{nv}^{\overline{B}} + O(1)$ on $\PP^n(K)$.

\medskip
Let $B$ be a projective arithmetic variety with $d = \dim B_{\QQ}$.
Let $\overline{H}_1, \ldots, \overline{H}_d$ be $C^{\infty}$-hermitian
$\QQ$-line bundles of surface type on $B$.
By its definition, for each $i$, there are  a morphism $\phi_i : B \to B_i$ of
flat and projective integral schemes over $\ZZ$ and
a $C^{\infty}$-hermitian $\QQ$-line bundle $\overline{L}_i$ on $B_i$
such that $\dim (B_i)_{\QQ} = 1$,
$\overline{L}_i$ is nef and big, and that
$\phi_i^*(\overline{L}_i) = \overline{H}_i$ in $\aPic(B) \otimes \QQ$.
We set $\overline{H} = \bigotimes_{i=1}^d \overline{H}_i$ and 
\[
\lambda_i = \exp\left(-\frac{\adeg(\acherncl_1(\overline{L}_i)^2)}{
\deg((L_i)_{\QQ})}\right).
\]
Let $K$ be the function field of $B$.
Here we consider several kinds of polarizations of $K$ as follows:
\[
\begin{cases}
\overline{B}_0 = (B; \overline{H}, \ldots, \overline{H}), \\
\overline{B}_1 = (B; \overline{H}_1, \ldots, \overline{H}_d), \\
\overline{B}_{i,j} = (B; \overline{H}_1, \ldots, \overline{H}_{j-1},
(\OO_B, \lambda_i\vert\cdot\vert_{can}), \overline{H}_{j+1}, \ldots,
\overline{H}_d) & \text{for $i \not= j$}.
\end{cases}
\]
Let $X$ be a geometrically irreducible projective variety over $K$, and
$L$ an ample line bundle on $X$. Let $(\XX, \overline{\LL})$ be
a model of $(X, L)$ over $B$. 
Then, for all $x \in X(\overline{K})$,
\renewcommand{\theequation}{\arabic{section}.\arabic{subsection}.\arabic{Theorem}}
\addtocounter{Theorem}{1}
\begin{equation}
\label{eqn:comp:heights}
h^{\overline{B}_0}_{(\XX, \overline{\LL})}(x)
= d! h^{\overline{B}_1}_{(\XX, \overline{\LL})}(x) + \frac{d!}{2}
\sum_{i\not=j} h^{\overline{B}_{i,j}}_{(\XX, \overline{\LL})}(x).
\end{equation}

Indeed, by Lemma~\ref{lem:power:sum:surface:type},
\[
h^{\overline{B}_0}_{(\XX, \overline{\LL})}(x) =
d! h^{\overline{B}_1}_{(\XX, \overline{\LL})}(x) + \frac{d!}{2}
\sum_{i\not=j} \frac{\adeg(\acherncl_1(\overline{L}_i)^2)\deg\left(\LL_{\QQ} \cdot
\prod^d_{l=1,l\not=j} f^* \phi_l^*(\overline{L}_l)_{\QQ} \cdot (\Delta_x)_{\QQ}\right)}{
\deg((L_i)_{\QQ})[K(x):K]},
\]
where $f : \XX \to B$ is the canonical morphism.
Moreover, 
\[
h^{\overline{B}_{i,j}}_{(\XX, \overline{\LL})}(x) =
\frac{{\displaystyle -\log(\lambda_i) \int_{\Delta_x(\CC)} 
c_1(\overline{\LL}) \wedge
\bigwedge^d_{l=1,l\not=j} c_1(f^*\phi_l^*(\overline{L}_l))}}
{[K(x):K]}.
\]
On the other hand,
\[
\int_{\Delta_x(\CC)} c_1(\overline{\LL})\wedge
\bigwedge_{l=1,l\not=j} c_1(\pi^*\phi_l^*(\overline{L}_l))
= \deg\left(\LL_{\QQ} \cdot
\prod^d_{l=1,l\not=j} f^*\phi_l^*(\overline{L}_l)_{\QQ}
\cdot (\Delta_x)_{\QQ} \right).
\]
Thus, we obtain
\[
h^{\overline{B}_{i,j}}_{(\XX, \overline{\LL})}(x) 
= \frac{\adeg(\acherncl_1(\overline{L}_i)^2)\deg\left(\LL_{\QQ} \cdot
\prod^d_{l=1,l\not=j} f^* \phi_l^*(\overline{L}_l)_{\QQ} \cdot (\Delta_x)_{\QQ}\right)
}{
\deg((L_i)_{\QQ})[K(x):K]}.
\]
Therefore, we get \eqref{eqn:comp:heights}.

\medskip
Using \eqref{eqn:comp:heights}, we can find a constant $C$ such that
\renewcommand{\theequation}{\arabic{section}.\arabic{subsection}.\arabic{Theorem}}
\addtocounter{Theorem}{1}
\begin{equation}
\label{eqn:comp:heights:2}
h^{\overline{B}_0}_{L}(x)
\leq C h^{\overline{B}_1}_{L}(x) + O(1)
\end{equation}
\renewcommand{\theequation}{\arabic{section}.\arabic{subsection}.\arabic{Theorem}.\arabic{Claim}}
for all $x \in X(\overline{K})$
because there is a positive integer $m$ such that
\[
\overline{H}_j^{\otimes m} \succsim (\OO_B, \lambda_i\vert\cdot\vert_{can})
\]
for every $i,j$.

\begin{Proposition}
\label{prop:comp:big:large:pol}
Let $X$ be a geometrically irreducible projective variety over $K$, and
$L$ an ample line bundle on $X$. Let
$\overline{B}$ and $\overline{B}'$ be fine polarizations of $K$.
Then $h_L^{\overline{B}'} \asymp h_L^{\overline{B}'}$
on $X(\overline{K})$
\rom{(}For the notation $\asymp$, see
\rom{\ref{subsubsec:asymp:notation}}\rom{)}.
\end{Proposition}

\Proof
It is sufficient to see that
there are a positive real number $a$ and a real number $b$ such that
$h_L^{\overline{B}} \leq a h_L^{\overline{B}'} + b$.
We set $\overline{B} = (B; \overline{H}_1, \ldots, \overline{H}_d)$ and
$\overline{B}' = (B'; \overline{H}'_1, \ldots, \overline{H}'_d)$.
Since $\overline{B}'$ is fine, by Proposition~\ref{prop:fine:pol:PP:1},
there are generically finite morphisms
$\mu' : B'' \to B'$ and $\nu : B'' \to \left(\PP^1_{\ZZ}\right)^d$ of
flat and projective integral schemes over $\ZZ$,
and nef and big $C^{\infty}$-hermitian $\QQ$-line bundles
$\overline{L}_1, \ldots, \overline{L}_d$
on $\PP^1_{\ZZ}$ such that
${\mu'}^*(\overline{H}'_i) \succsim 
\nu^*(p_i^*(\overline{L}_i))$ for all $i$,
where $p_i : \left(\PP^1_{\ZZ}\right)^d \to \PP^1_{\ZZ}$ is the projection
to the $i$-th factor. Changing $B''$ if necessarily, we may assume that there is
a generically finite morphism $\mu : B'' \to B$.

Let us consider polarizations
\[
\overline{B}_1 = (B''; \mu^*(\overline{H}_1), \ldots, \mu^*(\overline{H}_d))
\quad\text{and}\quad
\overline{B}_1' = (B''; {\mu'}^*(\overline{H}'_1), \ldots, {\mu'}^*(\overline{H}'_d))
\]
and compare $h^{\overline{B}}$ with $h^{\overline{B}_1}$ (resp. $h^{\overline{B}'}$ with
$h^{\overline{B}_1'}$).
By virtue of the projection formula, we may assume that
$B = B' = B''$ and $\mu = \mu' = \operatorname{id}$.

We set $\overline{H} = \nu^*\left( \bigotimes_{l=1}^d p_l^*(\overline{L}_l) 
\right)$. 
Then, $(B; \overline{H},\ldots, \overline{H})$ is a big polarization. 
Thus, by \cite[(5) of Proposition~3.3.7]{MoArht},
there is a positive integer $b_1$ such that
\[
h_L^{\overline{B}} \leq b_1 h_L^{(B; \overline{H}, \ldots, \overline{H})} 
+ O(1).
\]
Moreover, by \eqref{eqn:comp:heights:2},
we can find a positive constant $b_2$
with
\[
h_L^{(B; \overline{H}, \ldots, \overline{H})} \leq b_2 
h_L^{(B; \nu^*p_1^*(\overline{L}_1), \ldots, 
\nu^*p_d^*(\overline{L}_d))} + O(1).
\]
On the other hand, since ${\overline{H}'_i} \succsim 
\nu^*(p_i^*(\overline{L}_i))$ for all $i$,
\[
h_L^{(B; \nu^*p_1^*(\overline{L}_1), \ldots, \nu^*p_d^*(\overline{L}_d))}
\leq h_L^{(B; \overline{H}_1, \ldots, \overline{H}_d)} + O(1).
\]
Hence, we get our proposition.
\QED

\subsection{Comparisons of norms of polynomials}
\setcounter{Theorem}{0}

Let $S_n = \CC[z_1, \ldots, z_n]$ be the ring of $n$-variable
polynomials over $\CC$.
We define norms $\vert f \vert_{\infty}$ and
$\vert f \vert_2$ of
$f = \sum_{i_1, \ldots, i_n}
a_{i_1, \ldots, i_n} z_1^{i_1} \cdots z_n^{i_n} \in S_n$
as follows:
\[
\vert f \vert_{\infty} = \max_{i_1, \ldots, i_n} 
\{ \vert  a_{i_1, \ldots, i_n} \vert \} \quad\text{and}\quad
\vert f \vert_2 = \sqrt{\sum_{i_1, \ldots, i_n}
 \vert  a_{i_1, \ldots, i_n} \vert^2}.
\]
Moreover, the degree of $f$ with respect to
the variable $z_i$ is denoted by $\deg_i(f)$.

First of all, we have obvious inequalities:
\renewcommand{\theequation}{\arabic{section}.\arabic{subsection}.\arabic{Theorem}}
\addtocounter{Theorem}{1}
\begin{equation}
\label{eqn:comp:norm:2:infinite}
\vert f \vert_{\infty} \leq \vert f \vert_2 \leq
\sqrt{(\deg_1(f) + 1) \cdots (\deg_n(f) + 1)} \vert f \vert_{\infty}.
\end{equation}
We set
\[
S_n^{(d_1, \ldots, d_n)} =
\{ f \in S_n \mid
\deg_i(f) \leq d_i\quad(\forall i=1, \ldots, n)\}.
\]
Note that 
\addtocounter{Theorem}{1}
\begin{equation}
\dim_{\CC} S_n^{(d_1, \ldots, d_n)}
= (d_1 + 1) \cdots (d_n + 1).
\end{equation}
For $f_1, \ldots, f_l \in S_n$,
we set
\addtocounter{Theorem}{1}
\begin{equation}
\label{def:v:polynomial:eqn}
v(f_1, \ldots, f_l) = \exp \left( \int_{\CC^n} \log\left(
\max_{i} \{ \vert f_i \vert \} \right)
\omega_1 \wedge \cdots \wedge \omega_n \right),
\end{equation}
\renewcommand{\theequation}{\arabic{section}.\arabic{subsection}.\arabic{Theorem}.\arabic{Claim}}
where $\omega_i$'s are 
the $(1,1)$-forms  on $\CC^n$ given by
\[
\omega_i = \frac{\sqrt{-1}dz_i \wedge d\bar{z}_i}%
{2 \pi (1 + \vert z_i \vert^2)^2}.
\]
Let us begin with the following proposition.

\begin{Proposition}
\label{prop:comp:three:norms}
For $f_1, \ldots, f_l \in S_n^{(d_1, \ldots, d_n)}$,
we have the following.
\begin{enumerate}
\renewcommand{\labelenumi}{(\arabic{enumi})}
\item
$\max_{i} \{ \vert f_i \vert_{\infty} \}
\leq 2^{d_1 + \cdots + d_n} v(f_1, \ldots, f_l)$.

\item
$v(f_1, \ldots, f_l)
\leq \sqrt{2}^{d_1 + \cdots + d_n}
\sqrt{(\vert f_1 \vert_2)^2 + \cdots + (\vert f_l \vert_2)^2}$.
\end{enumerate}
\end{Proposition}

\Proof
(1) Since
\[
\max_{i} \left\{ \int_{\CC^n} \log\left(\vert f_i \vert \right)
\omega_1 \wedge \cdots \wedge \omega_n \right\} \leq
\int_{\CC^n} \log\left(
\max_{i} \{ \vert f_i \vert \} \right)
\omega_1 \wedge \cdots \wedge \omega_n,
\]
(1) is a consequence of \cite[Lemma~4.1]{MoArht}.

\medskip
For the proof of (2), we set
\[
\DD_0 = \{ z \in \CC \mid 0 < \vert z \vert < 1 \}
\quad\text{and}\quad
\DD_1 = \{ z \in \CC \mid 1 < \vert z \vert \}.
\]
Then,
\[
\int_{\CC^n} \log\left(
\max_{i} \{ \vert f_i \vert \} \right)
\omega_1 \wedge \cdots \wedge \omega_n
= \sum_{(\epsilon_1, \ldots, \epsilon_n) \in \{0,1\}^n}
\int_{\DD_{\epsilon_1} \times \cdots \times \DD_{\epsilon_n}}
\log\left(\max_{i} \{ \vert f_i \vert \}\right) \omega_1 \wedge \cdots \wedge \omega_n.
\]
For $\epsilon = (\epsilon_1, \ldots, \epsilon_n) \in \{0,1\}^n$,
let us consider a holomorphic map
\[
\varphi_{\epsilon} : \DD_0 \times \cdots \times \DD_0
\to \DD_{\epsilon_1} \times \cdots \times \DD_{\epsilon_n}
\]
given by $\varphi_{\epsilon}(z_1, \ldots, z_n)
= (z_1^{\iota(\epsilon_1)}, \ldots, z_n^{\iota(\epsilon_n)})$,
where $\iota : \{ 0, 1 \} \to \{-1, 1 \}$
is a map given by $\iota(0) = 1$ and $\iota(1) = -1$.
Then, since $\varphi_{\epsilon}^*(\omega_1 \wedge \cdots \wedge \omega_n)
= \omega_1 \wedge \cdots \wedge \omega_n$,
\begin{multline*}
\int_{\DD_{\epsilon_1} \times \cdots \times \DD_{\epsilon_n}}
\log\left(\max_{i} \{ \vert f_i \vert \} \right)
\omega_1 \wedge \cdots \wedge \omega_n \\
= \int_{\DD_{0}^n}
\log \left(\max_{i} \{ 
\vert f_i(z_1^{\iota(\epsilon_1)}, \ldots, z_n^{\iota(\epsilon_n)})\vert \}\right)
\omega_1 \wedge \cdots \wedge \omega_n.
\end{multline*}
Here we can find $f_{i, \epsilon} \in
S_n^{(d_1, \ldots, d_n)}$
such that 
\[
f_i(z_1^{\iota(\epsilon_1)}, \ldots, z_n^{\iota(\epsilon_n)}) =
\frac{f_{i, \epsilon}(z_1, \ldots, z_n)}%
{z_1^{\epsilon_1 d_1} \cdots z_n^{\epsilon_n d_n}}
\]
and $\vert f_i \vert_2 = \vert f_{i, \epsilon} \vert_2$.
Note that
\[
\int_{\DD_{0}^n}
\log(\vert z_i \vert) \omega_1 \wedge \cdots \wedge \omega_n =
- \frac{\log(2)}{2^n}
\]
for all $i$.
Therefore, 
\begin{multline*}
\int_{\DD_{\epsilon_1} \times \cdots \times \DD_{\epsilon_n}}
\log\left(\vert \max_{i} \{ \vert f_i \vert \} \vert\right)
\omega_1 \wedge \cdots \wedge \omega_n  \\
= \int_{\DD_{0}^n}
\log\left(\max_{i} \{ 
\vert f_{i,\epsilon}\vert \}\right)
\omega_1 \wedge \cdots \wedge \omega_n
- \sum_{i=1}^n \epsilon_i d_i
\int_{\DD_{0}^n}
\log(\vert z_i \vert) \omega_1 \wedge \cdots \wedge \omega_n \\
=
\int_{\DD_{0}^n}
\log\left(\max_{i} \{ 
\vert f_{i,\epsilon}\vert \}\right)
\omega_1 \wedge \cdots \wedge \omega_n 
+ \frac{\log(2)}{2^n}\sum_{i=1}^n \epsilon_i d_i.
\end{multline*}
Thus, we have
\begin{multline*}
\int_{\CC^n} \log\left(
\max_{i} \{ \vert f_i \vert \} \right)
\omega_1 \wedge \cdots \wedge \omega_n \\
= \sum_{\epsilon \in \{0,1\}^n}
\int_{\DD_{0}^n}
\log\left(\max_{i} \{ 
\vert f_{i,\epsilon}\vert \}\right)
\omega_1 \wedge \cdots \wedge \omega_n
+ \log(\sqrt{2})(d_1 + \cdots + d_n).
\end{multline*}
Hence, by the lemma below (Lemma~\ref{lem:comp:two:norms}), we can conclude
\begin{align*}
\int_{\CC^n} \log\left(
\max_{i} \{ \vert f_i \vert \} \right)
\omega_1 \wedge \cdots \wedge \omega_n & \leq  \sum_{\epsilon \in \{0,1\}^n} 
\frac{\log\left(
\sqrt{(\vert f_1 \vert_2)^2 + \cdots + (\vert f_l \vert_2)^2}
\right)}{2^n} \\
& \qquad\qquad + \log(\sqrt{2})(d_1 + \cdots + d_n)
\\ & =  \log\left(
\sqrt{(\vert f_1 \vert_2)^2 + \cdots + (\vert f_l \vert_2)^2}
\right) +  \log(\sqrt{2})(d_1 + \cdots + d_n).
\end{align*}
\QED

\begin{Lemma}
\label{lem:comp:two:norms}
For all $f_1, \ldots, f_l \in S_n$,
\[
\exp\left(  \int_{\DD_{0}^n} 
\log\left(\max_{i} \{ 
\vert f_{i}\vert \}\right)
(2\omega_1) \wedge \cdots \wedge (2\omega_n)
\right) \leq \sqrt{(\vert f_1 \vert_2)^2 + \cdots + (\vert f_l \vert_2)^2}.
\]
\end{Lemma}

\Proof
Let us begin with the following sublemma:

\begin{Sublemma}
Let $M$ be a differential manifold and $\Omega$ a volume form
on $M$ with $\int_M \Omega = 1$. Let $\varphi : \RR \to \RR$
be a $C^{\infty}$-function with $\varphi'' \geq 0$.
Let $u$ be a real valued function on $M$.
If $u$ and $\varphi(u)$ are integrable on $M$, then
$\varphi\left(\int_M u \Omega \right) \leq
\int_M \varphi(u) \Omega$.
\end{Sublemma}

\Proof
We set $c = \int_M u \Omega$. Since the second derivative
of $\varphi$ is non-negative, we can see
\[
(x - c) \varphi'(c) \leq \varphi(x) - \varphi(c)
\]
for all $x \in \RR$.
Therefore, we get
\[
\int_M(u - c) \varphi'(c) \Omega \leq 
\int_M(\varphi(u) - \varphi(c)) \Omega.
\]
On the other hand, the left hand side of the above inequality
is zero, and the right hand side is
$\int_M \varphi(u)\Omega - \varphi(c)$. Thus, we have our desired
inequality.
\QED

\bigskip
Let us go back to the proof of Lemma~\ref{lem:comp:two:norms}.
Applying the above lemma to the case $\varphi = \exp$,
\begin{align*}
\exp\left(  \int_{\DD_{0}^n} 
\log\left(\max_{i} \{ 
\vert f_{i}\vert^2 \}\right) (2\omega_1) \wedge \cdots \wedge (2\omega_n)
\right) & \leq \int_{\DD_{0}^n} \max_{i} \left\{
\vert f_i \vert^2 \right\} (2\omega_1) \wedge \cdots \wedge (2\omega_n) \\
& \leq \int_{\DD_{0}^n} \sum_{i}
\vert f_i \vert^2 (2\omega_1) \wedge \cdots \wedge (2\omega_n) 
\end{align*}
We set $f_i = \sum_{e_{1}, \ldots, e_{n}}
a^{(i)}_{e_1, \ldots, e_n} z_1^{e_1} \cdots z_n^{e_n}$ for all $i$. Then
\begin{multline*}
\sum_i \int_{\DD_{0}^n} 
\vert f_i \vert^2 (2\omega_1) \wedge \cdots \wedge (2\omega_n) = \\
\sum_i \sum_{\substack{e_1, \ldots, e_n, \\ e'_1, \ldots, e'_n}}
a^{(i)}_{e_1, \ldots, e_n} \overline{a^{(i)}_{e'_1, \ldots, e'_n}}
\int _{\DD_{0}^n} 
 z_1^{e_1} \bar{z_1}^{e'_1} \cdots  z_n^{e_n} \bar{z_n}^{e'_n}
(2\omega_1) \wedge \cdots \wedge
(2\omega_n).
\end{multline*}
It is easy to see that
\[
\int _{\DD_{0}^n} 
 z_1^{e_1} \bar{z_1}^{e'_1} \cdots  z_n^{e_n} \bar{z_n}^{e'_n}
(2\omega_1) \wedge \cdots \wedge
(2\omega_n) = 0
\]
if $(e_1, \ldots, e_n) \not= (e'_1, \ldots, e'_n)$. Moreover,
\[
\int _{\DD_{0}^n} 
 \vert z_1 \vert^{2e_1} \cdots  \vert z_n \vert^{2e_n}
(2\omega_1) \wedge \cdots \wedge
(2\omega_n) =
\left( \int_{\DD_{0}} \vert z_1 \vert^{2e_1} 2\omega_1 \right)
\cdots \left( \int_{\DD_{0}} \vert z_n \vert^{2e_n} 2\omega_n \right).
\]
Thus, it is sufficient to see that
\[
\int_{\DD_{0}} \vert z \vert^{2e} 
\frac{\sqrt{-1}dz \wedge
d\bar{z}}{\pi (1 + \vert z \vert^2)^2}
\leq 1
\]
for all $e \geq 0$.
We set $z = r \exp(\sqrt{-1}\theta)$, then
\[
\int_{\DD_{0}} \vert z \vert^{2e} 
\frac{\sqrt{-1}dz \wedge
d\bar{z}}{\pi (1 + \vert z \vert^2)^2} =
\int_{0}^{1}
\frac{4r^{2e+1}}{(1+r^2)^2} dr =
\int_{0}^{1} \frac{2 t^e}{(t+1)^2} dt
\]
If $e=0$, then the above integral is $1$.
Further if $e \geq 1$, then
\[
\int_{0}^{1} \frac{2 t^e}{(t+1)^2} dt\leq
\int_{0}^{1} 2 t^e dt = \frac{2}{e+1} \leq 1.
\]
\QED

Next let us consider the following proposition.

\begin{Proposition}
\label{prop:infinite:norm:product:two:polynomials}
For $f, g \in \CC[z_1, \ldots, z_n]$,
\[
\vert f\cdot g \vert_{\infty} \leq \vert f \vert_{\infty} \cdot \vert g \vert_{\infty}
\cdot \prod_i^n (1 + \min\{ \deg_i(f), \deg_i(g) \}).
\]
\end{Proposition}

\Proof
For $I \in (\ZZ_{\geq 0})^n$, the $i$-th entry of $I$ is denoted by
$I(i)$. A partial order `$\leq$' on $(\ZZ_{\geq 0})^n$ is
defined as follows:
\[
\text{$I \leq J$}
\quad\overset{\rm def}{\Longleftrightarrow}\quad
\text{$I(i) \leq J(i)$ for all $i=1, \ldots, n$}
\]
Moreover, for $I \in \ZZ_{\geq 0}^n$,
the monomial $z_1^{I(1)} \cdots z_n^{I(n)}$ is denoted by $z^I$.

Let us fix two non-zero polynomials
\[
f = \sum_{I \in (\ZZ_{\geq 0})^n} a_I z^I
\quad\text{and}\quad
g = \sum_{I \in (\ZZ_{\geq 0})^n} b_I z^I.
\]
We set $I_1 = (\deg_1(f), \ldots, \deg_n(f))$,
$I_2 = (\deg_1(g), \ldots, \deg_n(g))$ and
\[
d =  \prod_i^n (1 + \min\{ \deg_i(f), \deg_i(g) \}).
\]
First, we note that, for a fixed $I \in \ZZ_{\geq 0}^n$,
\[
\#\{ (J, J') \in \ZZ_{\geq 0}^n \times \ZZ_{\geq 0}^n
\mid \text{$J + J' = I$, $J \leq I_1$ and $J' \leq I_2$} \} \leq d.
\]
On the other hand,
\[
f \cdot g = \sum_{I} \left( 
\sum_{\substack{J+J' = I\\ J \leq I_1, J' \leq I_2}} a_{J} b_{J'} \right) z^I.
\]
Thus,
\[
\vert f \cdot g \vert_{\infty} \leq
\max_I \left\{ 
\sum_{\substack{J+J' = I\\ J \leq I_1, J' \leq I_2}} \vert a_{J} b_{J'} \vert 
\right\} \leq \max_I \left\{ 
\sum_{\substack{J+J' = I\\ J \leq I_1, J' \leq I_2}}
\vert f \vert_{\infty} \vert g \vert_{\infty} \right\}
\leq d \vert f \vert_{\infty} \vert g \vert_{\infty}
\]
\QED

For $f \in \CC[z_1, \ldots, z_n]$, we denote by $\lc_i(f)$ the coefficient
of the highest terms of $f$ as a polynomial of $z_i$, that is,
if we set
\[
f = a_n z_i^n + \cdots + a_0
\quad(a_i \in \CC[z_1, \ldots, z_{i-1}, z_{i+1}, \ldots, z_n],
\ a_n \not= 0),
\]
then $\lc_i(f) = a_n$.
Note that $\lc_i(0) = 0$ and
$\lc_i : \CC[z_1, \ldots, z_n] \to 
\CC[z_1, \ldots, z_{i-1},z_{i+1},\ldots,z_n]$.
For an element $\sigma$ of the $n$-th symmetric group $\frak{S}_n$,
we set
\[
\lc_{\sigma}(f) = \lc_{\sigma(n)} \cdots \lc_{\sigma(1)}(f).
\]
Then we have the following proposition.

\begin{Proposition}
\label{prop:lower:estimate:v}
For a non-zero $f \in \CC[z_1, \ldots, z_n]$,
\[
\int_{\CC^n} \log (\vert f \vert ) \omega_1 \wedge \cdots \wedge \omega_n 
\geq \max_{\sigma \in \frak{S}_n} \{ \log (\vert \lc_{\sigma}(f) \vert) \}.
\]
In particular, if $f \in \ZZ[z_1, \ldots, z_n]$, then
\[
\int_{\CC^n} \log (\vert f \vert ) \omega_1 \wedge \cdots \wedge \omega_n 
\geq 0.
\]
\end{Proposition}

\Proof
Changing the order of variables, it is sufficient to see that
\[
\int_{\CC^n} \log (\vert f \vert ) \omega_1 \wedge \cdots \wedge \omega_n 
\geq \log (\vert \lc_n \cdots \lc_1(f) \vert).
\]
We prove this by induction on $n$. First we assume $n = 1$.
Then, for $f = \alpha(z-c_1)\cdots(z-c_l)$,
\[
\int_{\CC} \log (\vert f \vert ) \omega_1 = \log \vert \alpha \vert +
\frac{1}{2} \sum_{i=1}^l \log(1 + \vert c_i \vert^2) \geq \log \vert \alpha
\vert.
\]

Next we consider a general $n$. By the hypothesis of induction,
we can see that
\[
\int_{\CC^{n-1}} \log (\vert f \vert ) 
\omega_1 \wedge \cdots \wedge \omega_{n-1} 
\geq \log \vert \lc_{n-1} \cdots \lc_1(f) \vert
\]
as a function with respect to $z_n$. Thus,
\[
\int_{\CC^n} \log (\vert f \vert ) \omega_1 \wedge \cdots \wedge \omega_n 
\geq \int_{\CC} \log \vert \lc_{n-1} \cdots \lc_1(f) \vert \omega_n \geq
\log \vert \lc_{n} \cdots \lc_1(f) \vert.
\]
\QED

\renewcommand{\theTheorem}{\arabic{section}.\arabic{Theorem}}
\renewcommand{\theClaim}{\arabic{section}.\arabic{Theorem}.\arabic{Claim}}
\renewcommand{\theequation}{\arabic{section}.\arabic{Theorem}.\arabic{Claim}}

\section{The number of arithmetic divisors}

Here let us consider several problems concerning the number of
arithmetic divisors with bounded arithmetic degree.

\begin{Proposition}
\label{prop:estimate:X:product:P:1:divisor:arith}
Let $p_i : (\PP^1_{\ZZ})^n \to \PP^1_{\ZZ}$ be
the projection to the $i$-th factor.
Let us fix a positive real number $\lambda$,
a subset $I$ of $\{ 1, \ldots, l\}$ and
a function $\alpha : I \to \ZZ_{\geq 0}$.
For a divisor $D$ on $(\PP^1_{\ZZ})^n$, we set
\[
\delta_{\lambda}(D) = 
\adeg(\acherncl_1(p_1^*(\overline{\OO}^{{\rm FS}_{\lambda}}(1))) \cdots
\acherncl_1(p_n^*(\overline{\OO}^{{\rm FS}_{\lambda}}(1))) \crest D).
\]
Then, 
there is a constant $C(\lambda, \alpha)$ depending only on
$\lambda$ and $\alpha : I \to \ZZ$ such that
\[
\# \left\{ D \in Z^{\rm eff}_{n}((\PP^1_{\ZZ})^n) \mid 
\text{{\rm $\delta_{\lambda}(D) \leq h$ and $\deg_i(D) \leq \alpha(i)$
for all $i \in I$}} \right\} \leq \exp
\left( C(\lambda, \alpha) h^{n+1 - \#(I)} \right)
\]
for $h \gg 0$.
\rom{(}Note that in the case where $I = \emptyset$,
no condition on $\deg_1(D), \ldots, \deg_n(D)$ is posed.\rom{)}
\end{Proposition}

\Proof
Fix a basis $\{ X_i, Y_i \}$ of $H^0(\PP^1_{\ZZ}, \OO(1))$ of the $i$-th factor
of $\left( \PP^1_{\ZZ} \right)^{n}$.
We denote by $\ZZ[X_1, Y_1, \ldots, X_{n}, Y_{n}]^{(k_1, \ldots, k_{n})}$
the set of homogeneous polynomials of multi-degree $(k_1, \ldots, k_{n})$.
Then,
\[
H^0\left((\PP^1_{\ZZ})^n, \bigotimes_{i=1}^{n} p_i^*(\OO(k_i)) \right) =
\ZZ[X_1, Y_1, \ldots, X_{n}, Y_{n}]^{(k_1, \ldots, k_{n})}.
\]
Let $D$ be an effective divisor on $(\PP^1_{\ZZ})^n$ with
$\deg_i(D) = k_i$ ($i=1, \ldots, n$).
Then there is 
\[
P \in \ZZ[X_1, Y_1, \ldots, X_{n}, Y_{n}]^{(k_1, \ldots, k_{n})} 
\setminus \{ 0 \}
\]
with $\zeros(P) = D$.
Let us evaluate
\[
\adeg\left(\acherncl_1(p_1^*(\overline{\OO}^{{\rm FS}_{\lambda}}(1))) \cdots
\acherncl_1(p_n^*(\overline{\OO}^{{\rm FS}_{\lambda}}(1))) \cdot 
\acherncl_1\left(\bigotimes_{i=1}^{n} 
p_i^*(\overline{\OO}^{{\rm FS}_0}(k_i))\right)\right)
\]
in terms of $P$, namely,
\[
(\lambda + 1/2)(k_1 + \cdots + k_{n})
= \delta_{\lambda}(D) - \int_{(\PP^1_{\CC})^{n}} 
\log \Vert P \Vert_{{\rm FS}_0}
\omega_1 \wedge \cdots \wedge \omega_{n},
\]
where $\omega_i = p_i^*(c_1(\overline{\OO}^{\rm FS}(1)))$ ($i=1,\ldots,n$).
We set $p(x_1, \ldots, x_{n}) = P(x_1, 1, \ldots, x_{n}, 1)$.
Then,
\[
\Vert P \Vert_{{\rm FS}_0} = \frac{\vert p \vert}{(1+x_1^2)^{k_1/2} \cdots
(1 + x_{n}^2)^{k_{n}/2}}.
\]
Note that
\[
\int_{(\PP^1_{\CC})^{n}} \log \left( (1+x_1^2)^{k_1/2} \cdots
(1 + x_{n}^2)^{k_{n}/2} \right) \omega_1 \wedge \cdots \wedge \omega_{n}
= \sum_{i=1}^{n} \frac{k_i}{2}.
\]
Thus,
\addtocounter{Claim}{1}
\begin{equation}
\label{prop:estimate:X:product:P:1:divisor:arith:eq:1}
\int_{(\PP^1_{\CC})^{n}} \log \vert p \vert
\omega_1 \wedge \cdots \wedge \omega_{n} =\delta_{\lambda}(D)
-\lambda(k_1 + \cdots + k_{n}).
\end{equation}
On the other hand, by Proposition~\ref{prop:comp:three:norms},
\begin{multline*}
\log \vert P \vert_{\infty} =
\log \vert p \vert_{\infty} \leq \log(2)(\deg_1(p) + \cdots + \deg_n(p)) +
\int_{(\PP^1_{\CC})^{n}} \log \vert p \vert
\omega_1 \wedge \cdots \wedge \omega_{n} \\ 
\leq \log(2)(k_1 + \cdots + k_n) +
\int_{(\PP^1_{\CC})^{n}} \log \vert p \vert
\omega_1 \wedge \cdots \wedge \omega_{n}.
\end{multline*}
Thus,
\addtocounter{Claim}{1}
\begin{equation}
\label{prop:estimate:X:product:P:1:divisor:arith:eq:2}
\log \vert P \vert_{\infty} \leq 
\delta_{\lambda}(D) + (k_1 + \cdots + k_n)(\log 2 - \lambda).
\end{equation}
We assume that $\delta_{\lambda}(D) \leq h$.
Then, since
\[
\int_{(\PP^1_{\CC})^{n}} \log \vert p \vert
\omega_1 \wedge \cdots \wedge \omega_{n} \geq 0,
\]
by \eqref{prop:estimate:X:product:P:1:divisor:arith:eq:1},
\addtocounter{Claim}{1}
\begin{equation}
\label{prop:estimate:X:product:P:1:divisor:arith:eq:3}
k_1 + \cdots + k_n \leq h/\lambda.
\end{equation}
Moreover, using \eqref{prop:estimate:X:product:P:1:divisor:arith:eq:2},
if $\lambda \leq \log 2$, then
\[
\log \vert P \vert_{\infty} \leq 
h + (k_1 + \cdots + k_n)(\log 2 - \lambda) \leq
h + \frac{h}{\lambda}(\log 2 - \lambda) = \frac{h \log 2}{\lambda}.
\]
Thus, if we set
\[
g(h, \lambda) =
\begin{cases}
\exp(h \log 2 / \lambda) & \text{if $0 < \lambda \leq \log 2$} \\
\exp(h) & \text{if $\lambda > \log 2$}, 
\end{cases}
\]
then
\addtocounter{Claim}{1}
\begin{equation}
\label{prop:estimate:X:product:P:1:divisor:arith:eq:4}
\vert P \vert_{\infty} \leq g(h, \lambda).
\end{equation}
Therefore,
\[
\# \{ P \in \ZZ[X_1, Y_1, \ldots, 
X_{n}, Y_{n}]^{(k_1, \ldots, k_{n})} \setminus \{ 0 \}
\mid \delta_{\lambda}(\zeros(P)) \leq h \}
\leq \left(2 g(h, \lambda) + 1
\right)^{(k_1+1)\cdots(k_n + 1)}
\]
Hence if we set
\[
N_{\alpha}(h) = \# \{ D \in Z^{\rm eff}_{n}((\PP^1_{\ZZ})^n) \mid 
\text{{\rm $\delta_{\lambda}(D) \leq h$ and $\deg_i(D) \leq \alpha(i)$
for all $i \in I$}} \},
\]
then we can see
\begin{align*}
N_{\alpha}(h) & \leq
\sum_{\substack{k_1 + \cdots + k_n \leq h/\lambda \\ 
k_i \leq \alpha(i) (\forall i \in I)}} 
\left(2 g(h, \lambda) + 1
\right)^{(k_1+1)\cdots(k_n + 1)} \\
& \leq 
(h/\lambda +1)^{l-\#(I)}\prod_{i\in I}(\alpha(i) + 1) 
\left(2g(h,\lambda) + 1
\right)^{(h/\lambda+1)^{n-\#(I)}\prod_{i\in I}(\alpha(i) + 1)}.
\end{align*}
Note that in the case where $I = \emptyset$,
the number $\prod_{i\in I}(\alpha(i) + 1)$ in the above
inequality is treated as $1$.
Thus, we get our lemma.
\QED

\begin{Proposition}
\label{prop:lower:bound:product:P:1}
Let us fix a positive real number $\lambda$.
For a divisor $D$ on $(\PP^1_{\ZZ})^n$, we set
\[
\delta_{\lambda}(D) = \adeg(\acherncl_1(p_1^*(\overline{\OO}^{{\rm FS}_{\lambda}}(1))) \cdots
\acherncl_1(p_n^*(\overline{\OO}^{{\rm FS}_{\lambda}}(1))) \crest D),
\]
where $p_i : (\PP^1_{\ZZ})^n \to \PP^1_{\ZZ}$ is the projection to the
$i$-th factor.
Let $x_1, \ldots, x_s$ be  closed points of
$(\PP^1_{\ZZ})^n$. Then, we have the following:
\begin{enumerate}
\renewcommand{\labelenumi}{(\arabic{enumi})}
\item
${\displaystyle \limsup_{h\to\infty}
\frac{\log \# \{ D \in \Div^{\rm eff}((\PP^1_{\ZZ})^n) 
\mid \text{{\rm $\delta_{\lambda}(D) \leq h$ and
$x_i \not\in \Supp(D)$ for all $i$}}\}}{h^{n+1}} > 0.}$

\item
We assume $n \geq 1$. Then
\[
\limsup_{h\to\infty}
\frac{\log \# \{ D \in \Div^{\rm eff}_{\rm hol}((\PP^1_{\ZZ})^n) 
\mid \text{{\rm $\delta_{\lambda}(D) \leq h$ and
$x_i \not\in \Supp(D)$ for all $i$}}\}}{h^{n+1}} > 0,
\]
where $\Div^{\rm eff}_{\rm hol}((\PP^1_{\ZZ})^n)$ is the set of all effective divisors
on $(\PP^1_{\ZZ})^n$
generated by prime divisors flat over $\ZZ$.
\end{enumerate}
\end{Proposition}

\Proof
Let us fix a coordinate $\{ X_i, Y_i \}$ of the $i$-th factor
of $(\PP^1_{\ZZ})^n$. Then, note that
\[
\bigoplus_{k_1 \geq 0, \ldots, k_n \geq 0}
H^0\left((\PP^1_{\ZZ})^n, \bigotimes_{i=1}^n p_i^*(\OO(k_i))\right) =
\ZZ[X_1, Y_1, \ldots, X_n, Y_n].
\]
We set $l = 4\prod_i \#(\kappa(x_i))$. Then,
$l = 0$ in $\kappa(x_i)$ for all $i$.
Since $H = \bigotimes_{i=1}^n p_i^*(\OO(1))$ is ample, 
there is a positive integer $k_0$
with $H^1((\PP^1_{\ZZ})^n, H^{\otimes k_0} \otimes m_{x_1} \otimes
\cdots \otimes m_{x_s}) = 0$,
where $m_{x_i}$ is the maximal ideal at $x_i$.
Thus, the homomorphism
\[
H^0((\PP^1_{\ZZ})^n, H^{\otimes k_0}) \to
\bigoplus_{i=1}^s H^{\otimes k_0} \otimes \kappa(x_i)
\]
is surjective. Hence, there is $P_0 \in H^0((\PP^1_{\ZZ})^n, H^{\otimes k_0})$
with $P_0(x_i) \not= 0$ for all $i$.
Clearly, we may assume that $P_0$ is primitive as a polynomial in
$\ZZ[X_1, Y_1, \ldots, X_n, Y_n]$.

For $m \geq 1$ and
$Q \in H^0((\PP^1_{\ZZ})^n, H^{\otimes m k_0})$, we set
$\alpha_m(Q) = P_0^m + l Q$.
Note that 
$\alpha_m(Q)(x_i) \not= 0$ for all $i$.
Thus, we get a map
\[
\phi_m : H^0((\PP^1_{\ZZ})^n, H^{\otimes m k_0}) \to
\{ D \in \Div^{\rm eff}((\PP^1_{\ZZ})^n) \mid \text{$x_i \not\in \Supp(D)$
for all $i$} \}
\]
given by $\phi_m(Q) = \zeros(\alpha_m(Q))$.
Here we claim that $\phi_m$ is injective.
Indeed, if $\phi_m(Q) = \phi_m(Q')$, then
$\alpha_m(Q) = \alpha_m(Q')$ or $\alpha_m(Q) = -\alpha_m(Q')$.
Clearly, if $\alpha_m(Q) = \alpha_m(Q')$,  then $Q = Q'$,
so that we assume $\alpha_m(Q) = -\alpha_m(Q')$.
Then $P_0^m = 2 \prod_i \#(\kappa(x_i)) (Q + Q')$.
Since $P_0$ is primitive, so is $P_0^m$.
This is a contradiction.

We set $d = (1 + k_0)^n$.
Let us choose a positive number $c$ with
\[
c \geq \max\left\{ \log(2d\vert P_0 \vert_{\infty}),
(\lambda + 1) k_0n \right\}.
\]

\begin{Claim}
If ${\displaystyle \vert Q \vert_{\infty} \leq \frac{\exp(cm)}{2l}}$, 
then $\delta_{\lambda}(\phi_m(Q)) \leq 2cm$.
\end{Claim}

We set $p_0 = P_0(x_1, 1, \ldots, x_n, 1)$ and
$q = Q(x_1, 1, \ldots, x_n, 1)$. Then
$\alpha_m(Q)(x_1, 1, \ldots, x_n, 1) = p_0^m + lq$.
By \eqref{prop:estimate:X:product:P:1:divisor:arith:eq:1},
\[
\delta_{\lambda}(\phi_m(Q)) = \lambda k_0mn +
\int_{(\PP^1_{\CC})^{n}} \log \vert p_0^m + lq \vert
\omega_1 \wedge \cdots \wedge \omega_{n}.
\]
Thus, using (2) of Proposition~\ref{prop:comp:three:norms},
\begin{align*}
\delta_{\lambda}(\phi_m(Q)) & \leq
\lambda k_0mn + \frac{k_0mn}{2} \log 2 + \frac{n \log(1 + k_0m)}{2}
+ \log \vert p_0^m + lq \vert_{\infty} \\
& \leq
(\lambda + 1)k_0nm + \log \vert p_0^m + lq \vert_{\infty}.
\end{align*}
On the other hand, using Lemma~\ref{prop:infinite:norm:product:two:polynomials} and
$\exp(c) \geq 2d\vert p_0 \vert_{\infty}$,
\begin{align*}
\vert p_0^m + lq \vert_{\infty} & \leq d^{m-1} \vert p_0 \vert_{\infty}^m
+ l \vert q \vert_{\infty} \leq
(d\vert p_0 \vert_{\infty})^m + l \vert Q \vert_{\infty} \leq
(d\vert p_0 \vert_{\infty})^m  + \frac{\exp(cm)}{2} \\
& \leq \frac{\exp(cm)}{2^m} + \frac{\exp(cm)}{2} 
\leq \exp(cm).
\end{align*}
Therefore, since $c \geq (\lambda + 1) k_0n$, we have
\[
\delta_{\lambda}(\phi_m(Q)) \leq  (\lambda + 1)k_0nm + cm \leq cm + cm = 2cm.
\]

\medskip
Let us go back to the proof of our proposition.
Since $H^0\left((\PP^1_{\ZZ})^n, \bigotimes_{i=1}^n p_i^*(\OO(k_0 m)\right)$
is a free abelian group of rank
$(1 + k_0 m)^n$,
\begin{align*}
\# \left\{ Q \in H^0\left((\PP^1_{\ZZ})^n \bigotimes_{i=1}^n p_i^*(\OO(k_0 m)\right)
\left| \vert Q \vert_{\infty} \leq \frac{\exp(cm)}{2l} \right\}\right.
= & \left(1 + 2\left[\frac{\exp(cm)}{2l}\right]\right)^{(1 + k_0 m)^n} \\
\geq &\left(1 + \left[\frac{\exp(cm)}{2l}\right]\right)^{(1 + k_0 m)^n} \\
\geq & \left(\frac{\exp(cm)}{2l}\right)^{(1 + k_0 m)^n}.
\end{align*}
Therefore, by the above claim,
\begin{multline*}
\log \#\{ D \in \Div^{\rm eff}((\PP^1_{\ZZ})^n) \mid \text{
$\delta_{\lambda}(D) \leq 2cm$ and $x_i \not\in \Supp(D)$
for all $i$} \} \\ \geq
(1 + k_0 m)^n (cm - \log(2l)).
\end{multline*}
Thus, we get (1).

\bigskip
From now, we assume $n > 0$.
We denote by $\mathcal{D}(h)$ (resp. $\mathcal{D}_{\rm hol}(h)$)
the set
\[
\{ D \in \Div^{\rm eff}((\PP^1_{\ZZ})^n) \mid \text{
$\delta_{\lambda}(D) \leq h$ and $x_i \not\in \Supp(D)$
for all $i$} \}
\]
\[
\left(  \text{resp.}\quad
\{ D \in \Div^{\rm eff}_{\rm hol}((\PP^1_{\ZZ})^n) \mid \text{
$\delta_{\lambda}(D) \leq h$ and $x_i \not\in \Supp(D)$
for all $i$} \}
\right) .
\]
For $D \in \mathcal{D}(h)$,
let $D = D_{\rm hol} + D_{\rm ver}$
be the unique decomposition such that
$D_{\rm hol}$ is horizontal over $\ZZ$ and
$D_{\rm ver}$ is vertical over $\ZZ$.
Note that $\delta_{\lambda}(D) = \delta_{\lambda}(D_{\rm hol}) + \delta_{\lambda}(D_{\rm ver})$,
$\delta_{\lambda}(D_{\rm hol}) \geq 0$ and $\delta_{\lambda}(D_{\rm ver}) \geq 0$.
Thus, $\delta_{\lambda}(D_{\rm hol}) \leq h$ and $\delta_{\lambda}(D_{\rm ver}) \leq h$.
Therefore, we have a map
\[
\beta_h : \mathcal{D}(h) \to \mathcal{D}_{\rm hol}(h)
\]
given by $\beta_h(D) = D_{\rm hol}$.
Since $\beta_h(D) = D$ for $D \in \mathcal{D}_{\rm hol}(h)$,
$\beta_h$ is surjective. Here let us consider
a fiber $\beta_h^{-1}(D)$ for $D \in \mathcal{D}_{\rm hol}(h)$.
First of all, an element $D' \in \beta_h^{-1}(D)$ has a form
\[
D' = D + \zeros(n)
\quad (n \in \ZZ \setminus \{ 0 \}).
\]
Since $\delta_{\lambda}(\zeros(n)) = \log \vert n \vert \leq h$,
we can see that $\# \beta_h^{-1}(D) \leq \exp(h)$.
Thus,
\[
\# \mathcal{D}(h) = \sum_{D \in \mathcal{D}_{\rm hol}(h)}
\# \beta_h^{-1}(D)  \leq \sum_{D \in \mathcal{D}_{\rm hol}(h)} \exp(h)
= \exp(h) \cdot \# \mathcal{D}_{\rm hol}(h).
\]
Hence, we get
\[
\limsup_{h \to \infty} \frac{\log \# \mathcal{D}_{\rm hol}(h)}{h^{n+1}} > 0.
\]
\QED

\begin{Remark}
\label{rem:for:prop:lower:bound:product:P:1}
In Proposition~\ref{prop:lower:bound:product:P:1},
we set $\overline{H}^{\lambda} = 
\bigotimes_{i=1}^n p_i^*(\overline{\OO}^{{\rm FS}_{\lambda}}(1))$.
Then, using Lemma~\ref{lem:power:sum:surface:type}, we can see
\[
\adeg_{\overline{H}^{\lambda}}(D) = n! \left( \delta_{\lambda}(D) + \frac{1+4\lambda}{4} \sum_{i=1}^n
\deg_i(D) \right).
\]
Moreover, using Lemma~\ref{prop:lower:estimate:v} and
\eqref{prop:estimate:X:product:P:1:divisor:arith:eq:1},
\[
\lambda \sum_{i=1}^n \deg_i(D) \leq \delta_{\lambda}(D).
\]
Thus,
\[
\adeg_{\overline{H}^{\lambda}}(D) \leq \frac{(8\lambda+1)n!}{4\lambda}
\delta_{\lambda}(D).
\]
Hence, Proposition~\ref{prop:lower:bound:product:P:1}
implies that if $n \geq 1$, then
\[
\limsup_{h\to\infty}
\frac{\log \# \{ D \in \Div^{\rm eff}_{\rm hol}((\PP^1_{\ZZ})^n) 
\mid \text{{\rm $\adeg_{\overline{H}^{\lambda}}(D) \leq h$ and
$x_i \not\in \Supp(D)$ for all $i$}}\}}{h^{n+1}} > 0.
\]
\end{Remark}

\renewcommand{\theTheorem}{\arabic{section}.\arabic{subsection}.\arabic{Theorem}}
\renewcommand{\theClaim}{\arabic{section}.\arabic{subsection}.\arabic{Theorem}.\arabic{Claim}}
\renewcommand{\theequation}{\arabic{section}.\arabic{subsection}.\arabic{Theorem}.\arabic{Claim}}

\section{The arithmetic case}
\subsection{Arithmetic cycles on the products of $\PP^1_{\ZZ}$}
\setcounter{Theorem}{0}
Here let us consider the number of cycles on $(\PP^1_{\ZZ})^n$.
Let us begin with the following lemma.

\begin{Lemma}
\label{lem:number:product:X:Y:arith}
Let $f : X \to S$ and $g : Y \to S$ be morphisms of
projective arithmetic varieties. We assume that $S$ is
of dimension $l \geq 1$.
Let $\overline{A}_1, \ldots, \overline{A}_l$ be
nef $C^{\infty}$-hermitian line bundles on $X$,
$\overline{B}_1, \ldots, \overline{B}_l$ 
nef $C^{\infty}$-hermitian line bundles on $Y$, and
$\overline{C}_1, \ldots, \overline{C}_l$ 
nef $C^{\infty}$-hermitian line
bundles on $S$ such that $\overline{A}_i \otimes
f^*(\overline{C}_i)^{\otimes -1}$ and
$\overline{B}_i \otimes g^*(\overline{C}_i)^{\otimes -1}$ are nef for all $i$ and that
$\adeg\left(\acherncl_1(\overline{C}_1) \cdots \acherncl_1(\overline{C}_l) \right) > 0$.
Let $p : X \times_S Y \to X$ and $q : X \times_S Y \to Y$ be
the projections to the first factor and the second factor respectively.
Fix $D \in Z_l^{\rm eff}(X/S)$ and $E \in Z^{\rm eff}_l(Y/S)$
\rom{(}for the definition of $Z_l^{\rm eff}(X/S)$ and $Z_l^{\rm eff}(Y/S)$,
see \rom{\ref{subsec:notation:cycle})}.
Then,
\begin{multline*}
\log \left( \# \left\{ V \in Z^{\rm eff}_l(X \times_S Y/S) \mid
\text{$p_*(V) = D$ and $q_*(V) = E$} \right\} \right) \\ 
\leq \min \left\{
\frac{\adeg(\acherncl_1(\overline{A}_1) \cdots \acherncl_1(\overline{A}_l) \crest D)
\adeg(\acherncl_1(\overline{B}_1) \cdots \acherncl_1(\overline{B}_l)\crest E)}
{\adeg\left(\acherncl_1(\overline{C}_1) \cdots \acherncl_1(\overline{C}_l) \right)^2}, \right. \\
\left.
\frac{\sqrt{\theta(D)\theta(E)\adeg(\acherncl_1(\overline{A}_1) \cdots \acherncl_1(\overline{A}_l) \crest D)
\adeg(\acherncl_1(\overline{B}_1)\cdots \acherncl_1(\overline{B}_l) \crest E)}}
{\adeg\left(\acherncl_1(\overline{C}_1) \cdots \acherncl_1(\overline{C}_l)\right)}
\right\},
\end{multline*}
where $\theta(D)$ \rom{(}resp. $\theta(E)$\rom{)} is 
the number of irreducible components of $\Supp(D)$ \rom{(}resp. $\Supp(E)$\rom{)}.
\end{Lemma}

\Proof
We set $D = \sum_{i=1}^s a_i D_i$ and $E = \sum_{j=1}^t b_j E_j$.
Then,
\addtocounter{Claim}{1}
\begin{multline}
\label{eqn:lem:number:product:X:Y:arith:1}
\adeg(\acherncl_1(\overline{A}_1)\cdots \acherncl_1(\overline{A}_l) \crest D) = 
\sum_{i} a_i \adeg(\acherncl_1(\overline{A}_1)\cdots \acherncl_1(\overline{A}_l) \crest D_i) \\
\geq 
\sum_{i=1}^s a_i \adeg(\acherncl_1(f^*(\overline{C}_1))\cdots \acherncl_1(f^*(\overline{C}_l)) \crest D_i) \\
= \sum_{i=1}^s a_i \deg(D_i \to S) \adeg(\acherncl_1(\overline{C}_1) \cdots \acherncl_1(\overline{C}_l))).
\end{multline}
In the same way,
\addtocounter{Claim}{1}
\begin{equation}
\label{eqn:lem:number:product:X:Y:arith:2}
\adeg(\acherncl_1(\overline{B}_1)\cdots \acherncl_1(\overline{B}_l) \crest E) \geq 
\sum_{j=1}^t b_j  \deg(E_j \to S) \adeg(\acherncl_1(\overline{C}_1)\cdots\acherncl_1(\overline{C}_l)).
\end{equation}
Thus, in the same way as in Lemma~\ref{lem:number:product:X:Y},
we have our assertion.
\QED

\begin{Proposition}
\label{prop:growth:P:1:product:l:horizontal:arith}
Let us fix a positive real number $\lambda$.
Let $p_i : (\PP^1_{\ZZ})^n \to \PP^1_{\ZZ}$ 
be the projection to the $i$-th factor.
We set $\overline{H}^{\lambda} = 
\bigotimes_{i=1}^{n} p_i^*(\overline{\OO}^{{\rm FS}_{\lambda}}(1))$.
For $1 \leq l \leq n$, we denote
by $Z^{\rm eff}_{l, {\rm hol}}((\PP^1_{\ZZ})^n)$ 
the set of all effective cycles on $(\PP^1_{\ZZ})^n$ generated by
$l$-dimensional integral closed subschemes of $(\PP^1_{\ZZ})^n$ 
which dominate $\Spec(\ZZ)$
by the canonical morphism $(\PP^1_{\ZZ})^n \to \Spec(\ZZ)$.
Then, there is a constant $C$ such that
\[
\# \{ V \in Z^{\rm eff}_{l, {\rm hol}}((\PP^1_{\ZZ})^n)  \mid
\adeg_{\overline{H}^{\lambda}}(V) \leq h \} \leq \exp(C\cdot h^{l+1})
\]
for all $h \geq 1$.
\end{Proposition}

\Proof
We set $\Sigma = \{ I \mid I \subseteq [n], \#( I ) = l-1 \}$.
Then, it is easy to see that
\[
Z^{{\rm eff}}_{l,{\rm hol}}((\PP^1_{\ZZ})^n) = 
\sum_{I \in \Sigma} Z^{{\rm eff}}_{l}((\PP^1_{\ZZ})^n
\overset{p_I}{\to} (\PP^1_{\ZZ})^{l-1}),
\]
where $p_I$ is the morphism given in \ref{subsec:product:P:1:projection}.
Thus, it is sufficient to show that
there is a constant $C'$ such that
\[
\# \{ V \in Z^{{\rm eff}}_{l}((\PP^1_{\ZZ})^n
\overset{p_I}{\to} (\PP^1_{\ZZ})^{l-1}
) \mid \adeg_{\overline{H}^{\lambda}}(V) \leq h \} \leq \exp(C'\cdot h^{l+1})
\]
for all $h \geq 1$.
By re-ordering the coordinate of $(\PP^1_{\ZZ})^n$, we may assume that
$I = [l-1]$. We denote $p_{[l-1]}$ by $p$.
We set
\[
T_n =Z^{\rm eff}_l((\PP^1_{\ZZ})^{n} \overset{p}{\to} (\PP^1_{\ZZ})^{l-1}).
\]
Let us see that $\{ T_n \}_{n=l}^{\infty}$ is a counting system.
First we define $h_n : T_n \to \RR_{\geq 0}$ to be
\[
h_n(V) = \adeg_{\overline{H}^{\lambda}}(V).
\]
Let $a_n : (\PP^1_{\ZZ})^n \to (\PP^1_{\ZZ})^{n-1}$ and
$b_n : (\PP^1_{\ZZ})^n \to (\PP^1_{\ZZ})^l$ be the morphisms given by
$a_n = p_{[n-1]}$ and $b_n = p_{[l-1]\cup \{ n \}}$.
Then, we have maps $\alpha_n : T_n \to T_{n-1}$ and
$\beta_n : T_n \to T_l$ defined by
$\alpha_n(V) = (a_n)_*(V)$ and
$\beta_n(V) = (b_n)_*(V)$.
Here, it is easy to see that
\[
h_{n-1}(\alpha_n(V)) \leq h_n(V)
\quad\text{and}\quad
h_{l}(\beta_n(V)) \leq h_n(V)
\]
for all $V \in T_n$.
Moreover, by Lemma~\ref{lem:number:product:X:Y:arith},
if we set
\[
e_l = \begin{cases}
\text{$\adeg\left( \acherncl_1(\otimes_{i=1}^{l-1}p_i^*(\overline{\OO}^{{\rm FS}_{\lambda}}(1)))
\right)$ on $(\PP^1_{\ZZ})^{l-1}$} & \text{if $l \geq 2$}  \\
\text{$\adeg(\acherncl_1(\ZZ, \exp(-\lambda)\vert\cdot\vert))$ on
$\Spec(\ZZ)$} & \text{if $l=1$}
\end{cases}
\]
and
\[
A(s, t) = \exp\left( 
\frac{s \cdot t}{e_l}
\right),
\]
then
\[
\# \{ x \in T_n \mid \alpha_n(x) = y, \beta_n(x) = z \}
\leq A(h_{n-1}(y), h_l(z))
\]
for all $y \in T_{n-1}$ and $z \in T_l$.
Further, by Proposition~\ref{prop:estimate:X:product:P:1:divisor:arith},
if we set
\[
B(h) = \exp(C'' \cdot h^{l+1})
\]
for some constant $C''$, then
\[
\{ x \in T_1 \mid h_1(x) \leq h \} \leq B(h)
\]
for all $h \geq 1$.
Thus, we can see that $\{T_n \}_{n=l}^{\infty}$ is a counting system.
Therefore, by virtue of Lemma~\ref{lem:counting:system}, we get our proposition.
\QED

\begin{Proposition}
\label{prop:product:P:1:vertival:arith}
Let $Z_{l, {\rm ver}}^{\rm eff}((\PP^1_{\ZZ})^n)$ be
the set of effective cycles on $(\PP^1_{\ZZ})^n$ generated by
$l$-dimensional integral subschemes which are not flat over $\ZZ$.
Then, there is a constant $B(n,l)$ depending only on 
$n$ and $l$ such that
\[
\# \{ V \in Z_{l, {\rm ver}}^{\rm eff}((\PP^1_{\ZZ})^n) \mid 
\adeg_{\overline{\OO}(1,\ldots,1)}(V) \leq h \}
\leq \exp\left( B(n,l) h^{l+1} \right)
\]
for $h \geq 1$.
\end{Proposition}

\Proof
For simplicity, we denote $(\PP^1_{\ZZ})^n$
and $\overline{\OO}(1,\ldots,1)$ by $X$ and $\overline{H}$ respectively.
Let $\pi : X \to \Spec(\ZZ)$ be the canonical morphism.
Let $k$ be a positive integer and
$k = \prod_i p_i^{a_i}$ the prime decomposition of $k$.
We set $X_{p_i} = \pi^{-1}([p_i])$.
Then,
\begin{multline*}
\# \{ V \in Z_{l, {\rm ver}}^{\rm eff}(X) \mid \adeg_{\overline{H}}(V) = 
\log(k) \}
\\
= \prod_i \#\{
V_i \in Z_{l, {\rm ver}}^{\rm eff}(X) \mid 
\text{$\pi(V_i) = [p_i]$ and $\deg_{\rest{H}{X_{p_i}}}(V_i) = a_i$} \}.
\end{multline*}
Let $C'(n,l)$ be a constant as in Proposition~\ref{prop:growth:P:1:product:l}.
We set $C''(n,l) = \max \{ C'(n,l), 1/\log(2) \}$. Note that
$C''(n,l)\log(p) \geq	 1$ for all primes $p$. Thus,
\begin{align*}
\log \# \{ V \in Z_{l, {\rm ver}}^{\rm eff}(X) \mid 
\adeg_{\overline{H}}(V) = \log(k) \}
& \leq \sum_{i} C''(n,l)\log(p_i) a_i^{l+1} \\
& \leq \left( \sum_{i} C''(n,l)\log(p_i) a_i \right)^{l+1} \\
& = C''(n,l)^{l+1} \log(k)^{l+1}.
\end{align*}
Therefore,
\begin{align*}
\# \{ V \in Z_{l, {\rm ver}}^{\rm eff}(X) \mid \adeg_{\overline{H}}(V) \leq h \} 
& \leq
\sum_{k=1}^{[\exp(h)]} 
\# \{ V \in Z_{l, {\rm ver}}^{\rm eff}(X) \mid \adeg_{\overline{H}}(V) = 
\log(k) \} \\
& \leq \sum_{k=1}^{[\exp(h)]} \exp\left( C''(n,l)^{l+1} \log(k)^{l+1} \right) \\
& \leq \exp(h) \cdot \exp\left( C''(n,l)^{l+1} h^{l+1} \right) \\
& = \exp\left( C''(n,l)^{l+1} h^{l+1} + h \right)
\end{align*}
Thus, we get the proposition.
\QED

By using  Proposition~\ref{prop:growth:P:1:product:l:horizontal:arith}
and Proposition~\ref{prop:product:P:1:vertival:arith},
we have the following:

\begin{Theorem}
\label{thm:growth:P:1:product:l:arith}
For all non-negative integers $l$ and $n$ with $0 \leq l \leq n$,
there is a constant $C$ such that
\[
\# \{ V \in Z^{\rm eff}_{l}((\PP^1_{\ZZ})^n) \mid 
\adeg_{\overline{H}^{\lambda}}(V) \leq h \} \leq \exp(Ch^{l+1})
\]
for all $h \geq 1$.
\end{Theorem}

As a variant of Proposition~\ref{prop:growth:P:1:product:l:horizontal:arith},
we have the following:

\begin{Proposition}
\label{prop:growth:P:1:product:l:horizontal:arith:fixed:g:degree}
Let us fix a positive real number $\lambda$.
Let $n$ and $d$ be non-negative integers with $n \geq d+1$.
Let $p_{[d]} : (\PP^1_{\ZZ})^n \to (\PP^1_{\ZZ})^d$ be the morphism
as in \rom{\ref{subsec:product:P:1:projection}}.
Let $p_i : (\PP^1_{\ZZ})^n \to \PP^1_{\ZZ}$ 
be the projection to the $i$-th factor.
For an integer $l$ with $d +1 \leq l \leq n$, we denote
by $Z^{\rm eff}_{l}((\PP^1_{\ZZ})^n/(\PP^1_{\ZZ})^{d})$ 
the set of all effective cycles on $(\PP^1_{\ZZ})^n$ generated by
$l$-dimensional integral closed subschemes of $(\PP^1_{\ZZ})^n$ 
which dominate $(\PP^1_{\ZZ})^d$
by $p_{[d]}$.
We set
\[
\adeg_{[d]}(V) = \begin{cases}
\adeg \left( 
\acherncl_1\left(\bigotimes_{i=1}^{n} 
p_i^*(\overline{\OO}^{{\rm FS}_{\lambda}}(1))\right)^{\cdot l-d} \cdot 
\prod_{j=1}^d
\acherncl_1(p_{j}^*(\overline{\OO}^{{\rm FS}_{\lambda}}(1)))
 \crest V
\right) & \text{if $d \geq 1$} \\
\adeg \left( 
\acherncl_1\left(\bigotimes_{i=1}^{n} 
p_i^*(\overline{\OO}^{{\rm FS}_{\lambda}}(1))\right)^{\cdot l} \crest V
\right) & \text{if $d = 0$}
\end{cases}
\]
for $V \in Z^{\rm eff}_{l}((\PP^1_{\ZZ})^n/(\PP^1_{\ZZ})^{d})$.
Let $K$ be the function field of $(\PP^1_{\ZZ})^d$.
For $V \in Z^{\rm eff}_{l}((\PP^1_{\ZZ})^n/(\PP^1_{\ZZ})^d)$,
we denote by $\deg_K(V)$ 
the degree of $V$ in the generic fiber of $\pi
: (\PP^1_{\ZZ})^n \to (\PP^1_{\ZZ})^d$ with respect to
$\OO_{\PP^{n-d}_K}(1, \ldots, 1)$.
Then, for a fixed $k$, there is a constant $C$ such that
\[
\# \{ V \in Z^{\rm eff}_{l}((\PP^1_{\ZZ})^n/(\PP^1_{\ZZ})^d) 
\mid \text{$\adeg_{\pi}(V) \leq h$
and $\deg_K(V) \leq k$} \} \leq \exp(C h^{d+1})
\]
for all $h \geq 1$.
\end{Proposition}

\Proof
We set
\[
\Sigma = \{ I \mid [d] \subseteq I \subseteq [n], \#( I ) = l-1 \}.
\]
Then, it is easy to see that
\[
Z^{{\rm eff}}_{l}((\PP^1_{\ZZ})^n/(\PP^1_{\ZZ})^d) = 
\sum_{I \in \Sigma} Z^{{\rm eff}}_{l}((\PP^1_{\ZZ})^n
\overset{p_I}{\to} (\PP^1_{\ZZ})^{l-1}).
\]
Thus, it is sufficient to see that, for each $I$,
there is a constant $C'$ such that
\[
\# \{ V \in Z^{{\rm eff}}_{l}((\PP^1_{\ZZ})^n
\overset{p_I}{\to} (\PP^1_{\ZZ})^{l-1}
) \mid \text{$\adeg_{\pi}(V) \leq h$ and
$\deg_K(V) \leq k$} \} \leq \exp(C' \cdot h^{d+1})
\]
for all $h \geq 1$. By changing the coordinate, we may assume that
$I = [l-1]$.
Here we denote $p_I$ by $p$.
For $n \geq l$, we set
\[
T_n =  \{ V \in Z^{{\rm eff}}_{l}((\PP^1_{\ZZ})^n
\overset{p}{\to} (\PP^1_{\ZZ})^{l-1}
) \mid \deg_K(V) \leq k \}.
\]
Let $a_n : (\PP^1_{\ZZ})^n \to (\PP^1_{\ZZ})^{n-1}$ and
$b_n : (\PP^1_{\ZZ})^n \to (\PP^1_{\ZZ})^{l}$ be morphisms
given by $a_n = p_{[n-1]}$ and $b_n = p_{[l-1]\cup \{ n \}}$.
Then, since
\[
(a_n)_K^*(\OO_{\PP^{n-d-1}_K}(1, \ldots, 1)) \otimes p_n^*(\OO(1))_K = \OO_{\PP^{n-d}_K}(1, \ldots, 1)
\]
and
\[
(b_n)_K^*(\OO_{\PP^{l-d}_K}(1, \ldots, 1)) \otimes \bigotimes_{i=l}^{n-1} p_i^*(\OO(1))_K 
= \OO_{\PP^{n-d}_K}(1, \ldots, 1),
\]
we have maps
$\alpha_n : T_n \to T_{n-1}$ and $\beta_n : T_n \to T_l$ given by
$\alpha_n(V) = (a_n)_*(V)$ and $\beta_n(V) = (b_n)_*(V)$.
Moreover, we set
\[
h_n(V) = \adeg_{[d]}(V).
\]
for $V \in T_n$. 
Then, it is easy to see that
\[
h_{n-1}(\alpha_n(V)) \leq h_n(V)
\quad\text{and}\quad
h_l(\beta_n(V)) \leq h_n(V)
\]
for all $V \in T_n$.
Note that
\[
k \geq \deg_K(V) \geq \theta(V) = \text{the number of irreducible components of $V$}.
\]
Thus,
by Lemma~\ref{lem:number:product:X:Y:arith},
if we set
\[
e_l = \begin{cases}
\text{$\adeg\left(\acherncl_1\left(\bigotimes_{i=1}^{l-1} 
p_i^*(\overline{\OO}^{{\rm FS}_{\lambda}}(1))\right)^{\cdot l-d} 
\prod_{i=1}^{d}\acherncl_1(p_{i}^*(\overline{\OO}^{{\rm FS}_{\lambda}}(1)))\right)$ on
$(\PP^1_{\ZZ})^{l-1}$} & \text{if $l \geq 2$ and $d \geq 1$} \\
\text{$\adeg\left(\acherncl_1\left(\bigotimes_{i=1}^{l-1} 
p_i^*(\overline{\OO}^{{\rm FS}_{\lambda}}(1))\right)^{\cdot l}\right)$ on
$(\PP^1_{\ZZ})^{l-1}$} & \text{if $l \geq 2$ and $d = 0$} \\
\text{$\adeg(\acherncl_1(\ZZ, \exp(-\lambda)\vert\cdot\vert))$ on $\Spec(\ZZ)$} &
\text{if $l=1$}
\end{cases}
\]
and
\[
A(s, t) = \frac{k\sqrt{s\cdot t}}{e_l},
\]
then
\[
\# \{ x \in T_n \mid \alpha_n(x) = y, \ \beta_n(x) = z \} \leq
A(h_n(x),h_l(y))
\]
for all $x \in T_{n-1}$ and $y \in T_l$.
Further, in the case where $n=l$,
\[
\deg_K(V) = \deg_{d+1}(V) + \cdots + \deg_{n}(V).
\]
Therefore, by Proposition~\ref{prop:estimate:X:product:P:1:divisor:arith},
there is a constant $C''$ such that
\[
\# \{ x \in T_l \mid h_l(x) \leq h \} \leq \exp(C'' \cdot h^{d+1})
\]
for all $h \geq 1$.
Hence, by Lemma~\ref{lem:counting:system},
there is a constant $C'$
such that
\[
\# \{ x \in T_n \mid h_n(V) \leq h \} \leq \exp(C' \cdot h^{d+1})
\]
for all $h \geq 1$.
Thus, we get our assertion.
\QED

\subsection{Upper estimate of cycles with bounded arithmetic degree}
\setcounter{Theorem}{0}
Here let us consider the following theorem, which is one of the main results of this paper.

\begin{Theorem}
\label{thm:growth:general:l:arith}
Let us fix a positive real number $\lambda$.
For all non-negative integers $l$ with $0 \leq l \leq \dim X$,
there is a constant $C$ such that
\[
\# \{ V \in Z^{\rm eff}_{l}(\PP^n_{\ZZ}) \mid 
\adeg_{\overline{\OO}^{{\rm FS}_{\lambda}}(1)}(V) \leq h \} \leq
\exp(C h^{l+1})
\]
for all $h \geq 1$.
\end{Theorem}

\Proof
Let us consider the birational map 
$\phi : \PP^n_{\ZZ} \dasharrow (\PP^1_{\ZZ})^n$ given by
\[
(X_0 : \cdots : X_n) \mapsto (X_0 : X_1) \times \cdots \times (X_0 : X_n).
\]
We set $U = \PP^n_{\ZZ} \setminus \{ X_0 = 0 \}$.
For $V \in Z_l^{\rm eff}(\PP^n_{\ZZ};U)$, we denote by $V'$ 
the strict transform of $V$ by $\phi$
(for the definition of
$Z_l^{\rm eff}(\PP^n_{\ZZ};U)$,
see \ref{subsec:notation:cycle}). 
Then, by applying Lemma~\ref{lem:adeg:comp:birat:product}
in the case $B = \Spec(\ZZ)$,
\[
n^l \adeg_{\overline{\OO}^{{\rm FS}_{\lambda}}(1)}(V) 
\geq \deg_{\overline{\OO}^{{\rm FS}_{\lambda}}(1, \ldots, 1)}(V').
\]
Moreover, 
if $V'_1 = V'_2$
for $V_1, V_2 \in Z_l^{\rm eff}(\PP^n_{\FF_q};U)$, then $V_1 = V_2$.
Therefore,
\addtocounter{Claim}{1}
\begin{equation}
\label{thm:growth:general:l:arith:eq:1}
\# \{ V \in Z_l^{\rm eff}(\PP^n_{\ZZ};U) 
\mid \adeg_{\overline{\OO}^{{\rm FS}_{\lambda}}(1)}(V) \leq h \} 
\leq
\# \{ V' \in Z_l^{\rm eff}((\PP^1_{\ZZ})^n) \mid 
\adeg_{\overline{\OO}^{{\rm FS}_{\lambda}}(1,\ldots,1)}(V') 
\leq n^l h \}.
\end{equation}
On the other hand, since $\PP^n_{\ZZ} \setminus U \simeq
\PP^{n-1}_{\ZZ}$, 
\addtocounter{Claim}{1}
\begin{multline}
\label{thm:growth:general:l:arith:eq:2}
\# \{ V \in Z_l^{\rm eff}(\PP^n_{\ZZ}) \mid 
\adeg_{\overline{\OO}^{{\rm FS}_{\lambda}}(1)}(V) \leq h \} \\
\leq \# \{ V \in Z_l^{\rm eff}(\PP^n_{\ZZ};U) \mid 
\adeg_{\overline{\OO}^{{\rm FS}_{\lambda}}(1)}(V) \leq h \} \cdot
\# \{ V \in Z_l^{\rm eff}(\PP^{n-1}_{\ZZ}) \mid 
\adeg_{\overline{\OO}^{{\rm FS}_{\lambda}}(1)}(V) \leq h \}
\end{multline}
Thus, using \eqref{thm:growth:general:l:arith:eq:1},
\eqref{thm:growth:general:l:arith:eq:2}, 
Theorem~\ref{thm:growth:P:1:product:l:arith}
and the hypothesis of induction, we have our theorem.
\QED

\begin{Corollary}
\label{cor:growth:general:l:arith}
Let $X$ be a projective arithmetic variety and $\overline{H}$ an ample
$C^{\infty}$-hermitian line bundle on $X$.
For all non-negative integers $l$ with $0 \leq l \leq \dim X$,
there is a constant $C$ such that
\[
\# \{ V \in Z^{\rm eff}_{l}(X) \mid \adeg_{\overline{H}}(V) \leq h \}
\leq \exp(C h^{l+1})
\]
for all $h \geq 1$.
\end{Corollary}

\Proof
Since $X$ is projective over $\ZZ$, there
is an embedding $\iota : X \hookrightarrow \PP^n_{\ZZ}$ over $\ZZ$.
We fix a positive real number $\lambda$.
Then, there is a positive integer $a$ such that
$\overline{H}^{\otimes a} \otimes 
\iota^*(\overline{\OO}^{{\rm FS}_{\lambda}}(-1))$ is ample.
Thus,
\[
a^l \adeg_{\overline{H}}(V) \geq 
\adeg_{\iota^*(\overline{\OO}^{{\rm FS}_{\lambda}}(1))}(V)
\]
for all $V \in Z^{\rm eff}_{l}(X)$.
Thus, our assertion follows from Theorem~\ref{thm:growth:general:l:arith}.
\QED

\subsection{Lower bound of the number of arithmetic cycles
with bounded degree}
\setcounter{Theorem}{0}
Here we consider the lower bound of the number of cycles.

\begin{Theorem}
Let $X$ be a projective arithmetic variety and
$\overline{H}$ an ample $C^{\infty}$-hermitian line
bundle on $X$.
Then, for $0 \leq l < \dim X$,
\[
\limsup_{h\to\infty}
\frac{\log \# \{ V \in Z_{l}^{\rm eff}(X) 
\mid \text{{\rm $\adeg_{\overline{H}}(V) \leq h$}}\}}{h^{l+1}} > 0.
\]
Moreover, if $0 < l < \dim X$, then
\[
\limsup_{h\to\infty}
\frac{\log \# \{ V \in Z_{l, {\rm hol}}^{\rm eff}(X) 
\mid \text{{\rm $\adeg_{\overline{H}}(V) \leq h$}}\}}{h^{l+1}} > 0.
\]
\end{Theorem}

\Proof
Choose a closed integral subscheme $Y$ of $X$ such that
$\dim Y = l+1$ and
$Y$ is flat over $\ZZ$.
First, we assume that $l = 0$.
Then, the canonical morphism $\pi: Y \to \Spec(\ZZ)$ is finite.
For $n \in \ZZ \setminus \{ 0 \}$,
\[
\adeg(\pi^*(\zeros(n))) = \deg(\pi) \adeg(\zeros(n)) = 
\deg(\pi) \log \vert n \vert.
\]
Thus,
\[
\#\{ V \in Z_{0}^{\rm eff}(Y) 
\mid \adeg(V) \leq h \} \geq
\# \{ \pi^*(\zeros(n)) \mid
\text{$n \in \ZZ \setminus \{ 0 \}$ and 
$\deg(\pi) \log \vert n \vert \leq h$} \}.
\]
Note that
\[
\pi^*(\zeros(n)) = \pi^*(\zeros(n'))
\ \Longrightarrow\  \zeros(n) = \zeros(n')
\ \Longrightarrow\  n = \pm n'.
\]
Thus,
\[
\# \{ \pi^*(\zeros(n)) \mid
\text{$n \in \ZZ \setminus \{ 0 \}$ and 
$\deg(\pi) \log \vert n \vert \leq h$} \}
= \left[\exp(h/\deg(\pi))\right].
\]
Therefore, 
\[
\limsup_{h\to\infty}
\frac{\log \# \{ V \in Z_{0}^{\rm eff}(X) 
\mid \adeg(V) \leq h \}}{h} > 0.
\]

\medskip
From now on, we assume that $l > 0$.
Since
\[
\# \{ D \in \Div^{\rm eff}_{\rm hol}(Y) 
\mid \text{{\rm $\adeg_{\overline{H}}(D) \leq h$}}\} \subseteq
\# \{ V \in Z_{l, {\rm hol}}^{\rm eff}(X) 
\mid \text{{\rm $\adeg_{\overline{H}}(V) \leq h$}}\},
\]
we may assume that
$\dim X = l+1$.

Let us take a birational morphism
$\mu : X' \to X$ of projective arithmetic varieties such that
there is a generically finite morphism
$\nu : X' \to (\PP^1_{\ZZ})^n$, where $n = \dim X_{\QQ}$.
We set $\overline{A} = \bigotimes_{i=1}^n 
p_i^*(\overline{\OO}^{{\rm FS}_{\lambda}}(1))$
on $(\PP^1_{\ZZ})^n$ for some positive real number $\lambda$.
Let us choose a positive rational number $a$ such that
there is a non-zero section $s \in
H^0(X', \nu^*(A)^{\otimes a} \otimes \mu^*(H)^{\otimes -1})$
with $\Vert s \Vert_{\rm sup} \leq 1$.
Let $X_0$ be a Zariski open set of $X$ such that
$\mu$ is an isomorphism over $X_0$.
Moreover, let $B$ be the non-flat locus of $\nu$.
Let 
\[
(X' \setminus \mu^{-1}(X_0)) \cup \Supp(\zeros(s)) \cup B 
= Z_1 \cup \cdots \cup Z_r
\]
be the irreducible decomposition.
Choose a closed point $z_i$ of $Z_i \setminus \bigcup_{j\not=i} Z_j$
for each $i$.
Then, by Proposition~\ref{prop:lower:bound:product:P:1} and 
Remark~\ref{rem:for:prop:lower:bound:product:P:1},
\[
\limsup_{h\to\infty}
\frac{\log \# \{ D' \in \Div^{\rm eff}_{\rm hol}((\PP^1_{\ZZ})^n) 
\mid \text{{\rm $\adeg_{\overline{A}}(D') \leq h$ and
$\nu(z_i) \not\in \Supp(D')$ for all $i$}}\}}{h^{n+1}} > 0.
\]
Let $D'$ be an element of $\Div^{\rm eff}_{\rm hol}((\PP^1_{\ZZ})^n)$
with $\nu(z_i) \not\in \Supp(D')$ for all $i$.
First, we claim that $\nu^*(D')$ is horizontal over $\ZZ$.
Assume the contrary, that is,
$\nu^*(D')$ contains a vertical irreducible component
$\Gamma$. Then, $\nu_*(\Gamma) = 0$, which implies $\Gamma \subseteq B$.
Thus, there is $z_i$ with $z_i \in \Gamma$.
Hence,
\[
z_i \in \Supp(\nu^*(D')) = \nu^{-1}(\Supp(D')),
\]
which contradicts the assumption $\nu(x_i) \not\in \Supp(D')$.

By the above claim, we can consider a map
\[
\phi : \{ D' \in \Div^{\rm eff}_{\rm hol}((\PP^1_{\ZZ})^n) 
\mid \text{{\rm $\nu(z_i) \not\in \Supp(D)$ for all $i$}} \} \to
\Div^{\rm eff}_{\rm hol}(X)
\]
given by $\phi(D') = \mu_*(\nu^*(D'))$. Here we claim that
$\phi$ is injective.
We assume that $\phi(D'_1) = \phi(D'_2)$. 
Since $z_i \not\in \Supp(\nu^*(D'_{\epsilon}))$ for
$\epsilon = 1, 2$ and all $i$,
no component of $\nu^*(D'_{\epsilon})$ is contained in
$X' \setminus \mu^{-1}(X_0)$. Thus, we have
$\nu^*(D'_{1}) = \nu^*(D'_{2})$. Hence
\[
\deg(\nu) D'_1 = \nu_*(\nu^*(D'_1)) = \nu_*(\nu^*(D'_2)) =
\deg(\nu) D'_2.
\]
Therefore, $D'_1 = D'_2$.

Let $D'$ be an element of $\Div^{\rm eff}_{\rm hol}((\PP^1_{\ZZ})^n)$
with $\nu(z_i) \not\in \Supp(D')$ for all $i$.
Since no component of $\nu^*(D')$ is contained in
$\Supp(\zeros(s))$, we can see 
\begin{align*}
\adeg_{\overline{H}}(\phi(D')) & =
\adeg_{\mu^*(\overline{H})}(\nu^*(D')) \leq 
\adeg_{\nu^*(\overline{A})^{\otimes a}}(\nu^*(D')) \\
& = a^n \adeg_{\nu^*(\overline{A})}(\nu^*(D'))
= a^n \deg(\nu) \adeg_{\overline{A}}(D').
\end{align*}
Thus
\begin{multline*}
\# \{ D' \in \Div^{\rm eff}_{\rm hol}((\PP^1_{\ZZ})^n) 
\mid \text{{\rm $\adeg_{\overline{A}}(D') \leq h$ and
$\nu(z_i) \not\in \Supp(D)$ for all $i$}} \} \\
\leq
\# \{ D \in \Div^{\rm eff}_{\rm hol}(X) 
\mid \text{{\rm $\adeg_{\overline{H}}(D) \leq \deg(\nu)a^n h$}} \}.
\end{multline*}
Therefore, we get our theorem.
\QED

\section{The arithmetic case with bounded geometric degree}

\subsection{Northcott's type in the arithmetic case}
The purpose of this subsection is to prove the following theorem,
which is a kind of refined Northcott's theorem.

\begin{Theorem}
\label{thm:relative:case:b:geom:degree}
Let $f : X \to B$ be a morphism of projective arithmetic varieties.
Let $K$ be the function field of $B$. 
Let $\overline{H}_1, \ldots, \overline{H}_d$ be a fine polarization of
$B$, where $d = \dim B_{\QQ}$.
Let $\overline{L}$ be a nef $C^{\infty}$-hermitian line bundle on $X$
such that $L_K$ is ample.
For an integer $l$ with $d+1 \leq l \leq \dim X$,
as in \rom{\ref{subsec:notation:cycle}},
let $Z_{l}^{\rm eff}(X/B)$ be
the set of effective cycles on $X$ generated by
integral closed $l$-dimensional
subschemes $\Gamma$ on $X$ with $f(\Gamma) = B$.
We denote by $Z^{\rm eff}_l(X/B, k, h)$
the set of effective cycle $V \in Z^{\rm eff}_l(X/B)$ with
\[
\adeg\left(  \acherncl_1(\overline{L})^{\cdot l-d} \cdot
\acherncl_1(f^*(\overline{H}_1)) \cdots 
\acherncl_1(f^*(\overline{H}_d)) \cdot V \right) \leq h
\quad\text{and}\quad
\deg(L_K^{\cdot l-d-1} \cdot V_K) \leq k.
\]
Then, for a fixed k, there is a constant $C$ such that
\[
\#Z_l^{\rm eff}(X/B, k, h) \leq \exp(C h^{d+1})
\]
for all $h \geq 1$.
\end{Theorem}

\Proof
First, let us consider a case where
$X = \PP^n_{\ZZ} \times_{\ZZ} B$,
$f$ is the natural projection $X \to B$ and
$\overline{L} = p^*(\overline{\OO}^{{\rm FS}_{1}}_{\PP^n_{\ZZ}}(1))$,
where $p : X \to \PP^n_{\ZZ}$ is the natural projection.
Since the polarization $\overline{H}_1, \ldots, \overline{H}_d$ is fine,
by Proposition~\ref{prop:fine:pol:PP:1},
there are generically finite morphisms
$\mu : B' \to B$ and $\nu : B' \to \left(\PP^1_{\ZZ}\right)^d$ of
projective arithmetic varieties, and positive rational numbers
$a_1, \ldots, a_d$  such that
$\mu^*(\overline{H}_i) \succsim 
\nu^*(r_i^*(\overline{\OO}^{{\rm FS}_{1}}_{\PP^1_{\ZZ}}(a_i)))$ for all $i$,
where $r_i : \left(\PP^1_{\ZZ}\right)^d \to \PP^1_{\ZZ}$ is the projection
to the $i$-th factor.
\[
\begin{CD}
B @<{\mu}<< B' @>{\nu}>> (\PP^1_{\ZZ})^d @>{r_i}>> \PP^1_{\ZZ}
\end{CD}
\]
We set $X' = \PP^n_{\ZZ} \times_{\ZZ} B'$, $B'' = (\PP^1_{\ZZ})^d$ and
$X'' = \PP^n_{\ZZ} \times_{\ZZ} B''$.
Let $p' : X' \to \PP^n_{\ZZ}$ and
$p'' : X'' \to \PP^n_{\ZZ}$ be the projections to the first factor and
$f' : X' \to B'$ and $f'' : X'' \to B''$ the projections to the last factor.
Here we claim the following.

\begin{Claim}
Let $K''$ be the function field of $B''$.
We denote by $Z^{\rm eff}_l(X''/B'', k, h)$
the set of effective cycles $V \in Z^{\rm eff}_l(X''/B'')$ with
\[
\adeg\left(  \acherncl_1(\overline{L}'')^{\cdot l-d} \cdot
\acherncl_1({f''}^*(r_1^*(\overline{\OO}^{{\rm FS}_{1}}_{\PP^1_{\ZZ}}(1)))) \cdots 
\acherncl_1({f''}^*(r_d^*(\overline{\OO}^{{\rm FS}_{1}}_{\PP^1_{\ZZ}}(1)))) \crest V \right) \leq h
\]
and
\[
\deg({L''}_{K''}^{\cdot l-d-1} \cdot V_{K''}) \leq k,
\]
where $\overline{L}'' = {p''}^*(\overline{\OO}_{\PP^n_{\ZZ}}^{{\rm FS}_{1}}(1))$.
Then, for a fixed k, there is a constant $C''$ such that
\[
\#Z_l^{\rm eff}(X''/B'', k, h) \leq \exp(C'' h^{d+1})
\]
for all $h \geq 1$.
\end{Claim}

Fixing $l$,
we prove this lemma by induction on $n$.
If $n = l-d-1$, then the assertion is trivial, so that
we assume $n > l -d -1$.
Let $\psi : \PP^n_{\ZZ} \dasharrow (\PP^1_{\ZZ})^n$ be the rational map 
given by
\[
(X_0 : \cdots : X_n) \mapsto (X_0 : X_1) \times \cdots \times (X_0 : X_n).
\]
We set $\phi = \psi \times \operatorname{id} : 
\PP^n_{\ZZ} \times B'' \dasharrow (\PP^1_{\ZZ})^n \times B''$ and
$U = (\PP^n_{\ZZ} \setminus \{ X_0 = 0 \})\times B''$.
Moreover, let $g : Y = (\PP^1_{\ZZ})^n \times B'' \to B''$ be the projection to the last factor,
and $s_i : Y \to \PP^1_{\ZZ}$ the projection to the $i$-th factor. We set
\[
\overline{M} = \bigotimes_{i=1}^n s_i^*(\overline{\OO}^{{\rm FS}_{1}}(1)).
\]
For $V \in Z_l^{\rm eff}(X''/B'' ; U)$
(i.e. $V \in Z_l^{\rm eff}(X/B)$ and any component of $\Supp(V)$
is not contained in $X \setminus U$), let $V'$ be the strict transform of $V$ via
$\phi$. Then, Lemma~\ref{lem:adeg:comp:birat:product} and Lemma~\ref{lem:comp:degree:field},
\begin{multline*}
n(l-d)\adeg\left(  \acherncl_1(\overline{L}'')^{\cdot l-d} \cdot
\acherncl_1({f''}^*(r_1^*(\overline{\OO}^{{\rm FS}_{1}}(1)))) \cdots 
\acherncl_1({f''}^*(r_d^*(\overline{\OO}^{{\rm FS}_{1}}(1)))) \crest V \right) \\
\geq \adeg\left(  \acherncl_1(\overline{M})^{\cdot l-d} \cdot
\acherncl_1({g}^*(r_1^*(\overline{\OO}^{{\rm FS}_{1}}(1)))) \cdots 
\acherncl_1({g}^*(r_d^*(\overline{\OO}^{{\rm FS}_{1}}(1)))) \crest V' \right)
\end{multline*}
and
\[
n(l-d-1) \deg({L''}_{K''}^{l-d-1} \cdot V_{K''}) \geq	 \deg({M}_{K''}^{l-d-1} \cdot V'_{K''}).
\]
On the other hand, $\{ X_0 = 0 \} \times B = \PP^{n-1}_{\ZZ} \times B$.
Thus, by the hypothesis of induction and
Proposition~\ref{prop:growth:P:1:product:l:horizontal:arith:fixed:g:degree},
we have our claim.

\medskip
We can define a homomorphism
\[
({\rm id} \times \mu)^{\star} :
Z_l(X/B) \to  Z_l(X'/B')
\]
with $({\rm id} \times \mu)_* ({\rm id} \times \mu)^{\star} (V) = \deg(\mu)V$
as follows:
Let $B_0$ be the locus of points of $B$ over which $B' \to B$ is flat.
We set $X_0 = f^{-1}(B_0)$. Then, for $V \in Z_l^{\rm eff}(X/B)$,
no component of $\Supp(V)$ is contained in $X \setminus X_0$.
Thus, $({\rm id} \times \mu)^{\star}(V)$ is defined by
the Zariski closure of $({\rm id} \times \mu)^{*}(\rest{V}{X_0})$.
By virtue of $({\rm id} \times \mu)^{\star}$,
\begin{multline*}
\adeg\left( \acherncl_1
\left({p'}^*(\overline{\OO}^{{\rm FS}_{1}}(1))\right)^{\cdot l-d} \cdot
\acherncl_1({f'}^*\mu^*(\overline{H}_1)) \cdots 
\acherncl_1({f'}^*\mu^*(\overline{H}_d))\crest ({\rm id} \times \mu)^{\star}(V) 
\right) \\
= \deg(\mu)
\adeg\left( 
\acherncl_1
\left(p^*(\overline{\OO}^{{\rm FS}_{1}}(1))\right)^{\cdot l-d} \cdot
\acherncl_1(f^*(\overline{H}_1)) \cdots \acherncl_1(f^*(\overline{H}_d))
\crest V \right).
\end{multline*}
Thus, in order to prove the theorem in our case, we may assume that $B' = B$. 
Then,
\begin{multline*}
\adeg\left( 
\acherncl_1
\left(p^*(\overline{\OO}^{{\rm FS}_{1}}(1))\right)^{\cdot l-d} \cdot
\acherncl_1(f^*(\overline{H}_1)) \cdots \acherncl_1(f^*(\overline{H}_d))
\crest V \right)
\geq a_1 \cdots a_d  \times \\
\adeg\left( 
\acherncl_1
\left({p''}^*(\overline{\OO}^{{\rm FS}_{1}}(1))\right)^{\cdot l-d} \cdot
\acherncl_1({f''}^*(r_{1}^*(\overline{\OO}^{{\rm FS}_{1}}(1))) \cdots 
\acherncl_1({f''}^*(r_{d}^*(\overline{\OO}^{{\rm FS}_{1}}(1)))\crest
({\rm id} \times \nu)_*(V) \right).
\end{multline*}
Thus, by the above Claim,
it is sufficient to see the following:
For a fixed $V'$ on $Y$,
the number of effective cycles $V$ on $X$
with $({\rm id} \times \nu)_*(V) = V'$ is less than or equal to
$\exp(\deg(\nu)k)$.
For, let $V' = \sum e_i V'_i$ be the irreducible decomposition.
Then by Lemma~\ref{lem:estimate:number:push:forward},
the above number is less than or equal to
$\exp(\deg(\nu) \sum e_i)$.
On the other hand, we can see
\[
\sum e_i \leq \sum_i e_i \deg(L_K^{\cdot l-d-1} \cdot (V_i)_K) =
\deg(L_K^{\cdot l-d-1} \cdot V_K) \leq k.
\]
Thus, we get the theorem in our case.

\bigskip
Let us go back to the proof of the theorem in a general case.
Replacing $\overline{L}$ by a positive multiple of it,
we may assume that $L_K$ is very ample.
Thus we have an embedding $\phi : X_K \hookrightarrow \PP^n_K$ with
$\phi^*(\OO(1)) = L_K$. Let $X'$ be the Zariski closure
of $X_K$ in $\PP^n_{\ZZ} \times_{\ZZ} B$ and
$f' : X' \to B$ the induced morphism.
Let $p : \PP^n_{\ZZ} \times_{\ZZ} B \to \PP^n_{\ZZ}$
be the projection to the first factor and
$\overline{L}' = \rest{p^*(\overline{\OO}^{{\rm FS}_{1}}(1))}{X'}$.
Then there are birational morphisms
$\mu : Z \to X$ and $\nu : Z \to X'$ of projective arithmetic varieties.
We set $g = f \cdot \mu = f' \cdot \nu$.
Let $A$ be an ample line bundle on $B$
such that
$g_*(\mu^*(L) \otimes \nu^*(L')^{\otimes -1}) \otimes A$ is generated by
global sections.  Thus there is a non-zero global section
$s \in H^0(Z, \mu^*(L) \otimes \nu^*(L')^{\otimes -1} \otimes g^*(A))$.
Since $(\mu^*(L) \otimes \nu^*(L')^{\otimes -1})_K = \OO_{X_K}$,
we can see that $f(\zeros(s)) \subsetneq B$.
We choose a metric of $A$ with $\Vert s \Vert \leq 1$.
For $V \in Z^{\rm eff}_l(X/B)$, let $V_1$ be the strict transform of $V$ by $\mu$ and
$V' = \nu_*(V_1)$. Then,
\begin{multline*}
\adeg\left( \acherncl_1(\mu^*(\overline{L} \otimes 
f^*(\overline{A}))^{\cdot l+1} \cdot
\acherncl_1(g^*(\overline{H}_1)) \cdots 
\acherncl_1(g^*(\overline{H}_d)) \crest V_1 \right) \\
= \adeg\left( \acherncl_1(\overline{L})^{\cdot l+1} \cdot
\acherncl_1(f^*(\overline{H}_1)) \cdots 
\acherncl_1(f^*(\overline{H}_d)) \crest V \right) \\
+ (l+1)\deg(L_K^{l} \cdot V_K)\adeg\left( \acherncl_1(\overline{A}) \cdot
\acherncl_1(\overline{H}_1) \cdots 
\acherncl_1(\overline{H}_d)\right)
\end{multline*}
Moreover,
\begin{multline*}
\adeg\left( \acherncl_1(\mu^*(\overline{L} \otimes 
f^*(\overline{A}))^{\cdot l+1} \cdot
\acherncl_1(g^*(\overline{H}_1)) \cdots 
\acherncl_1(g^*(\overline{H}_d)) \crest V_1 \right) \\
\geq \adeg\left( \acherncl_1(\nu^*\overline{L}')^{\cdot l+1} \cdot
\acherncl_1(g^*(\overline{H}_1)) \cdots 
\acherncl_1(g^*(\overline{H}_d)) \crest V_1 \right) \\
=\adeg\left( \acherncl_1(\overline{L}')^{\cdot l+1} \cdot
\acherncl_1({f'}^*(\overline{H}_1)) \cdots 
\acherncl_1({f'}^*(\overline{H}_d)) \crest V' \right)
\end{multline*}
Thus, we may assume that there is an embedding $X \hookrightarrow \PP^n_{\ZZ} \times_{\ZZ} B$ and
$\overline{L} = p^*(\overline{\OO}^{{\rm FS}_{1}}(1))$.
Therefore we get our Theorem.
\QED

\subsection{The number of rational points over
a finitely generated field}
\setcounter{Theorem}{0}
Let us consider the polarization
\[
\overline{B}_1 = (\left(\PP^1_{\ZZ}\right)^d;
p_1^*(\overline{\OO}^{{\rm FS}_1}_{\PP^1_{\ZZ}}(1)), \ldots, 
p_d^*(\overline{\OO}^{{\rm FS}_1}_{\PP^1_{\ZZ}}(1)))
\]
of $\QQ(z_1, \ldots, z_d)$, where
$p_i : \left(\PP^1_{\ZZ}\right)^d \to \PP^1_{\ZZ}$ is the projection to
the $i$-th factor.
$\overline{B}_1$ is called {\em the standard polarization of $\QQ(z_1, \ldots, z_d)$}.
First, let us see the following lemma.

\begin{Lemma}
\label{lem:estimate:num:rat:Q:z}
\[
\limsup_{h\to\infty}
\frac{\log \#\{ x \in \PP^n(\QQ(z_1, \ldots, z_d)) \mid 
h_{nv}^{\overline{B}_1}(x) \leq h \}}{h^{d+1}} > 0
\]
\rom{(}See \rom{\ref{subsec:height:functions:over:fgf}} for the definition
$h_{nv}^{\overline{B}_1}$.\rom{)}
\end{Lemma}

\Proof
We set $\overline{H}_i = p_i^*(\overline{\OO}^{{\rm FS}_1}_{\PP^1_{\ZZ}}(1))$
for $i = 1, \ldots, d$.
Clearly we may assume $n=1$.
Let $\Delta_{\infty}$ be the closure of $\infty \in \PP^1_{\QQ}$
in $\PP^1_{\ZZ}$. We set $\Delta^{(i)}_{\infty} = p_i^*(\Delta_{\infty})$.
Then
\addtocounter{Claim}{1}
\begin{equation}
\label{eqn:lem:estimate:num:rat:Q:z:1}
\adeg\left(
\acherncl_1(\overline{H}_1) \cdots
\acherncl_1(\overline{H}_d) \crest \Delta^{(i)}_{\infty}
\right)  = 1
\end{equation}
Let $P$ be a $\QQ(z_1, \ldots, z_d)$-valued point of $\PP^1$.
Then, there are
$f_0, f_1 \in \ZZ[z_1, \cdots, z_d]$
such that $f_0$ and $f_1$ are relatively prime and
$P = (f_0 : f_1)$.
Thus, by \eqref{eqn:lem:estimate:num:rat:Q:z:1}
\begin{multline*}
 h_{nv}^{\overline{B}}(P) =
\sum_i \max \{ \deg_i(f_0), \deg_i(f_1) \}
+ \int_{(\PP^1)^d} \log \left( \max \{ |f_0|, |f_1| \} \right) 
c_1(\overline{H}_1)
\wedge \cdots \wedge c_1(\overline{H}_d).
\end{multline*}
Let $a$ be a positive number with $1 - 2da > 0$.
We set
\[
\mathcal{S}(h) =
\left\{ f \in \ZZ[z_1, \ldots, z_d] \mid \text{$v(1, f) \leq \exp((1-da)h)$ and
$\deg_i(f) \leq [ah]$ for all $i$} \right\}.
\]
(See \eqref{def:v:polynomial:eqn} for the definition $v$.)

First we claim that
$h_{nv}^{\overline{B}_0}((1:f)) \leq h$ for
all $f \in \mathcal{S}(h)$. If $f=0$, then the assertion is obvious.
We assume that $f \not= 0$. Then,
\[
h_{nv}^{\overline{B}_0}((1:f)) = \sum_{i=1}^d \deg_i(f) + \log(v(1,f))
\leq d[ah] + (1-da)h \leq h.
\]

Next we claim that
\[
\vert f \vert_{\infty} \leq \frac{\exp((1-2ad)h)}{\sqrt{2}}
\quad\Longrightarrow\quad
v(1,f) \leq \exp((1-ad)h).
\]
For this purpose, we may assume that $f \not= 0$.
Moreover, note that $\sqrt{1+x^2} \leq \sqrt{2}x$ for $x \geq 1$.
Thus, using
\eqref{eqn:comp:norm:2:infinite} and Proposition~\ref{prop:comp:three:norms},
\begin{align*}
v(1,f) & \leq \sqrt{2}^{d[ah]}\sqrt{1 + \vert f \vert_2^2} \leq
\sqrt{2} \sqrt{2}^{dah} \vert f \vert_2 \leq
\sqrt{2} \sqrt{2}^{dah} (1 + [ah])^{d/2} \vert f \vert_{\infty} \\
& \leq \sqrt{2} \exp(dah/2) \exp(dah/2) \vert f \vert_{\infty}
\leq \sqrt{2}\exp(dah) \frac{\exp((1-2ad)h)}{\sqrt{2}}
= \exp((1-ad)h).
\end{align*}

By the second claim,
\begin{align*}
\#\mathcal{S}(h) & \geq
\left(1 +  2 \left[ \frac{\exp((1-2ad)h)}{\sqrt{2}} \right]\right)^{([ah]+1)^d}
\geq \left(\frac{\exp((1-2ad)h)}{\sqrt{2}}\right)^{(ah)^d} \\
& \geq \exp((1-2ad)h - 1)^{(ah)^d} = \exp(a^d(1-2ad)h^{d+1} - a^d h^d).
\end{align*}
Thus, we get our lemma by the first claim.
\QED

\begin{Theorem}
\label{thm:upper:asymp:estimate:height}
Let $K$ be a finitely generated field over $\QQ$, and
$\overline{B}$ a fine polarization of $K$.
Let $X$ be a projective variety over $K$ and $L$ an ample line
bundle on $X$. Then
\[
\#\{ x \in \PP^n(K) \mid 
h_{\OO(1)}^{\overline{B}}(x) \leq h \} \leq \exp(C h^{d+1})
\]
for all $h \geq 0$, where $d = \trdeg_{\QQ}(K)$.
\end{Theorem}

\Proof
This is a consequence of Theorem~\ref{thm:relative:case:b:geom:degree}.
\QED

\begin{Theorem}
\label{thm:lower:asymp:estimate:height}
Let $K$ be a finitely generated field over $\QQ$ and
$\overline{B}$ a fine polarization of $K$.
Then, 
\[
\limsup_{h\to\infty} 
\frac{\log \#\{ x \in \PP^n(K)\mid h_{\OO(1)}^{\overline{B}}(x) \leq h \}}{h^{d+1}}
> 0,
\]
where $d = \trdeg_{\QQ}(K)$.
\end{Theorem}

\Proof
If we set $d = \trdeg_{\QQ}(K)$, then
there is a subfield $\QQ(z_1, \ldots, z_d)$
such that $K$ is finite over $\QQ(z_1, \ldots, z_d)$.
Let $\overline{B}_1$ be the standard polarization of
$\QQ(z_1, \ldots, z_d)$ as in Lemma~\ref{lem:estimate:num:rat:Q:z}.
Then,
by Lemma~\ref{lem:estimate:num:rat:Q:z},
\[
\limsup_{h\to\infty}
\frac{\log \#\{ x \in \PP^n(K) \mid 
h_{\OO(1)}^{\overline{B}_1}(x) \leq h \}}{h^{d+1}} > 0.
\]
Let $\overline{B}_1^K$ be the polarization of $K$
induced by $\overline{B}_1$.
Then,
\[
[K:\QQ(z_1, \ldots, z_d)]h_{\OO(1)}^{\overline{B}_1}(x) =
h_{\OO(1)}^{\overline{B}_1^{K}}(x) +O(1)
\]
for all $x \in \PP^n(K)$.
Moreover, by Proposition~\ref{prop:comp:big:large:pol},
$h_{\OO(1)}^{\overline{B}_1^{K}} \asymp 
h_{\OO(1)}^{\overline{B}}$.
Thus, we get our theorem.
\QED

\section{Zeta functions of algebraic cycles}
In this section, we would like to propose a kind of zeta functions arising from
the number of algebraic cycles.
First let us consider a local case, i.e.,
the case over a finite field.

\subsection{Local case}
Let $X$ be a projective variety over a finite field $\FF_q$ and
$H$ an ample line bundle on $X$.
For a non-negative integer $k$, 
we denote by $n_k(X,H,l)$
the number of all effective $l$-dimensional cycles $V$ on $X$ 
with $\deg_{H}(V) = k$.
We define a zeta function $Z(X,H,l)$ 
of $l$-dimensional cycles on
a polarized scheme $(X,H)$ over $\FF_q$ to be
\[
Z(X,H,l)(T) = \sum_{k=0}^{\infty} n_k(X,H,l) T^{k^{l+1}}.
\]
Then, we have the following:

\begin{Theorem}
\label{thm:convergent:zeta:geom}
$Z(X,H,l)(T)$ is a convergent power series at the origin.
\end{Theorem}

\Proof
First note that
$n_{m^l k}(X, H^{\otimes m}, l) = n_k(X,H, l)$.
Moreover, if we choose $m > 0$ with
$H^{\otimes m}$ very ample, then, by Corollary~\ref{cor:general:geom:case},
there is a constant $C$ with $n_k(X, H^{\otimes m}, l) \leq q^{Ck^{l+1}}$.
Thus,
\[
n_k(X,H, l) = n_{m^l k}(X, H^{\otimes m}, l) \leq q^{C'k^{l+1}},
\]
where $C' = Cm^{l(l+1)}$.
Therefore, if $\vert q^{C'} T \vert < 1$, then
\[
\sum_{k=0}^{\infty} n_k(X,H,l) \vert T^{k^{l+1}} \vert
\leq \sum_{k=0}^{\infty} \vert q^{C'} T \vert^{k^{l+1}} \leq
\sum_{k=0}^{\infty} \vert q^{C'} T \vert^{k} =
\frac{1}{1 - \vert q^{C'} T \vert}.
\]
Thus, we get our theorem.
\QED

See Remark~\ref{rem:wan:zeta:function} for Wan's zeta functions. Next, let us consider
height zeta functions in the local case, which
is a local analogue of Batyrev-Manin-Tschinkel's height zeta functions
(cf. \cite{BatMan}).

\begin{Theorem}
Let $K$ be a finitely generated field over a finite field $\FF_q$ with
$d = \trdeg_{\FF_q}(K) \geq 1$.
Let $X$ be a projective variety over $K$ and
$L$ a ample line bundle on $X$.
Let $h_L$ be a representative of the class of height functions
associated with $(X, L)$ as in \rom{\ref{subsec:geometric:height:function:F:q}}.
Then, for a fixed $k$, a series
\[
\sum_{\substack{x \in X(\overline{K}), \\ [K(x):K] \leq k}} q^{-s (h_L(x))^d}
\]
converges absolutely and uniformly on the compact set in
$\{ s \in \CC \mid \Re(s) > C \}$ for some $C$.
\end{Theorem}

\Proof
We set
\[
X_n = \{ x \in X(\overline{K}) \mid \text{$n-1 < h_L(x) \leq n$ and $[K(x):K] \leq k$} \}
\]
for $n > 1$ and
\[
X_1 = \{ x \in X(\overline{K}) \mid \text{$h_L(x) \leq 1$ and $[K(x):K] \leq k$} \}.
\]
Then, by Corollary~\ref{cor:refine:northcott:F:q}, there is a constant $C$ such that
$\#(X_n) \leq q^{Cn^d}$
for all $n \geq 1$.
\[
\sum_{\substack{x \in X(\overline{K}), \\ [K(x):K] \leq k}} \vert q^{-s (h_L(x))^d}\vert = \sum_{n=1}^{\infty} \sum_{x \in X_n} q^{-\Re(s)(h_L(x))^d}
\leq \sum_{n=1}^{\infty} q^{Cn^d} q^{-\Re(s)n^d} =
\sum_{n=1}^{\infty} \left( q^{-(\Re(s)-C)} \right)^{n^d}.
\]
Thus, we have our assertion.
\QED

\subsection{Global case}
Let $K$ be a number field and $O_K$ the ring of integers in $K$.
Let $f : X \to \Spec(O_K)$ be a flat and projective scheme over $O_K$ and
$H$ an $f$-ample line bundle on $X$.
For $P \in \Spec(O_K) \setminus \{ 0 \}$, we denote by $X_P$ the fiber
of $f$ at $P$. Here let us consider an infinite product
\[
L(X,H,l)(s) = \prod_{P \in \Spec(O_K) \setminus \{ 0 \}} Z(X_P, H_P, l)(\#(\kappa(P))^{-s})
\]
for $s \in \CC$.
Then, we have the following:

\begin{Theorem}
\label{thm:convergent:product:zeta:geom}
There is a constant $C$ such that
the infinite product
$L(X,H,l)(s)$
converges absolutely and uniformly on the compact set in
$\{ s \in \CC \mid \Re(s) > C \}$.
\end{Theorem}

\Proof
Since $Z(X_P, H_P^{\otimes n}, l)(q^{-s}) =
Z(X_P,H_P,l)(q^{-l(l+1)s})$, replacing
$H$ by $H^{\otimes n}$ for some positive number $n$,
we may assume that $H$ is very $f$-ample.
For non-negative integer $k$,
we set
\[
n_k(X_P,H_P,l) = \#\{ V \in Z^{\rm eff}_l(X_P) \mid \deg_{H_P}(V) = k \}.
\]
Then, 
\[
Z(X_P, H_P,l)(\#(\kappa(P))^{-s}) = 1 + \sum_{k=1}^{\infty}
n_k(X_P,H_P,l) \#(\kappa(P))^{-sk^{l+1}}.
\]
We denote $\sum_{k=1}^{\infty}
n_k(X_P,H_P,l) \#(\kappa(P))^{-sk^{l+1}}$ by $u_P(s)$.
We set $N = \rank(f_*(H)) - 1$.
Then, for each $P$, we have an embedding
$\iota_P : X_P \hookrightarrow \PP^N_{\kappa(P)}$
with $\iota_P^*(\OO(1)) = H_P$. Thus, by Theorem~\ref{thm:proj:sp:geom:case},
there is a constant $C$ depending only on $l$ and $N$
with
$n_k(X_P,H_P,l) \leq \kappa(P)^{C k^{l+1}}$ for all $k \geq 1$.
Thus, for $s \in \CC$ with $\Re(s) > C + 1$,
\begin{align*}
\vert u_P(s) \vert & \leq \sum_{k=1}^{\infty} \#(\kappa(P))^{C k^{l+1}} 
\#(\kappa(P))^{-\Re(s)k^{l+1}}  = \sum_{k=1}^{\infty} 
\#(\kappa(P))^{-(\Re(s)-C)k^{l+1}} \\
& \leq \sum_{i=1}^{\infty} \#(\kappa(P))^{-(\Re(s)-C)i} =
\frac{\#(\kappa(P))^{-(\Re(s)-C)}}{1 - \#(\kappa(P))^{-(\Re(s)-C)}} \leq
\#(\kappa(P))^{-(\Re(s)-C)}.
\end{align*}
Therefore, we have
\begin{align*}
\sum_{P \in \Spec(O_K) \setminus \{ 0 \}} \vert u_P(s) \vert
& \leq \sum_{P \in \Spec(O_K) \setminus \{ 0 \}} \kappa(P)^{-(\Re(s)-C)} \\
& \leq \sum_{n=1}^{\infty} n^{-(\Re(s)-C)} = \zeta(\Re(s) - C).
\end{align*}
Hence, we get our theorem by the criterion of the convergence of infinite products.
\QED

\subsection{Arithmetic case}
Next let us consider an analogue in Arakelov geometry.
Let $\mathcal{X}$ be a projective arithmetic variety and
$\overline{\mathcal{H}}$ an ample $C^{\infty}$-hermitian
$\QQ$-line bundle on $\mathcal{X}$.
For an effective cycle $V$ of $l$-dimension,
the norm of $V$ is defined by
\[
N_{\overline{\mathcal{H}}}(V) = 
\exp\left( \adeg_{\overline{\mathcal{H}}}(V)^{l+1} \right).
\]
Then, the zeta function of $(\mathcal{X},\overline{\mathcal{H}})$ 
for cycles of dimension $l$ is defined by
\[
\zeta(\mathcal{X}, \overline{\mathcal{H}}, l)(s) 
= \sum_{V \in Z^{\rm eff}_l(\mathcal{X})} 
N_{\overline{\mathcal{H}}}(V)^{-s}
\]

\begin{Theorem}
\label{thm:zeta:convergent:arithmetic}
There is a constant $C$ such that the above
$\zeta(\mathcal{X}, \overline{\mathcal{H}}, l)(s)$
converges absolutely and uniformly on the compact set in
$\{ s \in \CC \mid \Re(s) > C \}$.
\end{Theorem}

\Proof
We denote $Z^{\rm eff}_l(\mathcal{X}, \overline{\mathcal{H}}, h)$ the set of
all $l$-dimensional effective cycles on $\mathcal{X}$ with
$\adeg_{\overline{\mathcal{H}}}(V) \leq h$.
By Corollary~\ref{cor:growth:general:l:arith},
there is a constant $C$ such that
\[
\#\left(Z^{\rm eff}_l(\mathcal{X}, \overline{\mathcal{H}}, h)\right) 
\leq \exp(C h^{l+1})
\]
for all $h \geq 1$. We choose a positive constant $C'$
with
\[
\exp(C(h+1)^{l+1}) \leq \exp(C'h^{l+1})
\]
for all $h \geq 1$.
Moreover, for a real number $x$, we set
$[x] = \max \{ n \in \ZZ \mid n \leq x \}$.
Note that if $k = [\adeg_{\overline{\mathcal{H}}}(V)]$, then
$k \leq \adeg_{\overline{\mathcal{H}}}(V) < k+1$.
Thus, for $s \in \CC$ with $\Re(s) > C'$,
\begin{align*}
\sum_{V \in Z^{\rm eff}_l(\mathcal{X})} 
\left| N_{\overline{\mathcal{H}}}(V)^{-s} \right| & =
\sum_{k=0}^{\infty} \sum_{\substack{V \in Z^{\rm eff}_l(\mathcal{X}) \\
[\adeg_{\overline{\mathcal{H}}}(V)] = k}} 
\left| N_{\overline{\mathcal{H}}}(V) \right|^{-\Re(s)} \\
& \leq \sum_{k=0}^{\infty} \#(Z^{\rm eff}_l(\mathcal{X}, 
\overline{\mathcal{H}}, k+1))
\exp(k^{l+1})^{-\Re(s)} \\
& \leq \sum_{k=0}^{\infty} \exp(C(k+1)^{l+1})
\exp(k^{l+1})^{-\Re(s)} \\
& \leq \exp(C) + \sum_{k=1}^{\infty} \exp(-(\Re(s)-C'))^{k^{l+1}} \\
& \leq \exp(C) + \frac{\exp(-(\Re(s)-C'))}{1 - \exp(-(\Re(s)-C'))}.
\end{align*}
Thus, we get our theorem.
\QED

\subsection{Remarks}
Here let us discuss remarks of the previous zeta functions.
The first one is the abscissa of convergence of zeta functions.

\begin{Remark}
Let $I$ be an index set and
$\{ \lambda_i \}_{i \in I}$ a sequence of real numbers
such that  
the set
$I(t) = \{ i \in I \mid \lambda_i \leq t \}$
is finite for every real number $t$.
Then the abscissa $\sigma_0$ of convergence of the Dirichlet series
\[
\sum_{i \in I} \exp(-\lambda_i s) = \lim_{t\to\infty}
\sum_{i \in I(t)} \exp(-\lambda_i s) 
\]
is given by
\[
\sigma_0 = \limsup_{t \to \infty} 
\frac{\log\left( \#(I(t))\right)}{t}.
\]

Let $X$ be a projective scheme over $\FF_q$ and $H$ an ample line bundle
on $X$. Moreover, let $\mathcal{X}$ be a projective arithmetic variety
and $\overline{\mathcal{H}}$ an $C^{\infty}$-hermitian line bundle
on $\mathcal{X}$. We denote by $\sigma_0(X,H,l)$ (resp.
$\sigma_0(\mathcal{X}, \overline{\mathcal{H}}, l)$)
the abscissa of convergence of $Z(X,H,l)(q^{-s})$ (resp.
$\zeta(\mathcal{X}, \overline{\mathcal{H}}, l)(s)$).
Then, $\sigma_0(X,H,l)$ and
$\sigma_0(\mathcal{X}, \overline{\mathcal{H}}, l)$ are given by
\[
\sigma_0(X,H,l) = \limsup_{h\to\infty}
\frac{\log_q \# \left(
\{ V \in Z^{\rm eff}_l(X) \mid \deg_H(V) \leq h \} \right)}{h^{l+1}}
\]
and
\[
\sigma_0(\mathcal{X}, \overline{\mathcal{H}},l) =
\limsup_{h\to\infty} \frac{\log \# \left(
\{ V \in Z^{\rm eff}_l(\mathcal{X}) \mid 
\adeg_{\overline{\mathcal{H}}}(V) \leq h \} \right)}{h^{l+1}}
\]
respectively.

For example, let $X$ be an $n$-dimensional projective scheme over $\FF_q$ with
$\Pic(X) = \ZZ \cdot H$, where $H$ is ample. Then,
\[
\sigma_0(X, H, n-1) = \frac{1}{\deg(H^{n})^{n-1} n!}.
\]
\end{Remark}

\begin{Remark}
\label{rem:wan:zeta:function}
Let $X$ be a projective variety over a finite field $\FF_q$ and
$H$ an ample line bundle on $X$. As before,
the number of all effective $l$-dimensional cycles $V$ on $X$ 
with $\deg_{H}(V) = k$ is denoted by $n_k(X,H,l)$.
In \cite{Wan}, Wan defined a zeta function $\tilde{Z}(X,H,l)$ by
\[
\tilde{Z}(X,H,l)(T) = \sum_{k=0}^{\infty} n_k(X,H,l) T^{k}.
\]
He proved $\tilde{Z}(X,H,l)(T)$ is $p$-adic analytic and proposed
several kinds of conjectures.
Of course, $\tilde{Z}(X,H,l)(T)$ is never analytic as $\CC$-valued functions
if $0 < l < \dim X$.
In order to get classical analytic functions, we need to replace $T^k$ by $T^{k^{l+1}}$.
\end{Remark}

\renewcommand{\thesection}{Appendix \Alph{section}}
\renewcommand{\theTheorem}{\Alph{section}.\arabic{Theorem}}
\renewcommand{\theClaim}{\Alph{section}.\arabic{Theorem}.\arabic{Claim}}
\renewcommand{\theequation}{\Alph{section}.\arabic{Theorem}.\arabic{Claim}}
\setcounter{section}{0}
\section{Bogomolov plus Lang in terms of a fine polarization}

\begin{Theorem}[{\cite[Theorem~4.3]{MoArht}}]
\label{thm:northcott:intro}
We assume that the polarization
$\overline{B}$ is fine.
Let $X$ be a geometrically irreducible projective variety over $K$, and
$L$ an ample line bundle on $X$.
Then, for any number $M$ and any positive integer $e$, the set
\[
\{ x \in X(\overline{K}) \mid h^{\overline{B}}_L(x) \leq M, \quad
[K(x) : K] \leq e \}
\]
is finite.
\end{Theorem}

\begin{Theorem}[{\cite[Theorem~A]{MoBL}}]
\label{thm:BL:conj:fair:pol:intro}
We assume that the polarization
$\overline{B}$ is fine.
Let $A$ be an abelian variety over $K$, and $L$ a symmetric
ample line bundle on $A$.
Let 
\[
\langle\ , \ \rangle_L^{\overline{B}} :
A(\overline{K}) \times A(\overline{K}) \to \RR
\]
be a paring given by
\[
\langle x , y \rangle_L^{\overline{B}} =
\frac{1}{2} \left( \hat{h}_L^{\overline{B}}(x+y)
-  \hat{h}_L^{\overline{B}}(x) -  \hat{h}_L^{\overline{B}}(x) \right).
\]
For $x_1, \ldots, x_l \in A(\overline{K})$,
we denote $\det \left( \langle x_i, x_j \rangle_L^{\overline{B}} \right)$
by $\delta_L^{\overline{B}}(x_1, \ldots, x_l)$.

Let $\Gamma$ be a subgroup of finite rank in $A(\overline{K})$
\rom{(}i.e., $\Gamma \otimes \QQ$ is finite-dimensional\rom{)}, and
$X$ a subvariety of $A_{\overline{K}}$.
Fix a basis $\{\gamma_1, \ldots, \gamma_n \}$ of $\Gamma \otimes \QQ$.
If the set
\[
\{ x \in X(\overline{K}) \mid
\delta_L^{\overline{B}}(\gamma_1, \ldots, \gamma_n, x) \leq \epsilon \}
\]
is Zariski dense in $X$ for every positive number $\epsilon$,
then $X$ is a translation of an abelian subvariety of $A_{\overline{K}}$
by an element of $\Gamma_{div} = \{ x \in A(\overline{K}) \mid
\text{$nx \in \Gamma$ for some positive integer $n$} \}$.
\end{Theorem}

{\bf The proof of Theorem~\ref{thm:northcott:intro} and
Theorem~\ref{thm:BL:conj:fair:pol:intro}:}
Here, let us give the proof of Theorem~\ref{thm:northcott:intro},
Theorem~\ref{thm:BL:conj:fair:pol:intro}. 
Theorem~\ref{thm:northcott:intro} is obvious by
\cite[Theorem~4.3]{MoArht} and
Proposition~\ref{prop:comp:big:large:pol}, or
Theorem~\ref{thm:relative:case:b:geom:degree}.
Theorem~\ref{thm:BL:conj:fair:pol:intro} is a consequence of
\cite{MoBL}, Proposition~\ref{prop:comp:big:large:pol} and
the following lemma.

\begin{Lemma}
Let $V$ be a vector space over $\RR$, and
$\langle\ , \ \rangle$ and $\langle\ , \ \rangle'$ be two inner products on $V$.
If $\langle x , x \rangle \leq \langle x , x \rangle'$ for all $x \in V$, then
$\det \left( \langle x_i , x_j \rangle \right) \leq \det 
\left( \langle x_i , x_j \rangle' \right)$
for all $x_1, \ldots, x_n \in V$.
\end{Lemma}

\Proof
If $x_1, \ldots, x_n$ are linearly dependent, then our assertion is trivial.
Otherwise, it is nothing more than \cite[Lemma~3.4]{MoBG}.
\QED

\begin{Remark}
\label{rem:north:not:hold}
In order to guarantee Northcott's theorem,
the fineness of a polarization is crucial.
The following example shows us that even if the polarization
is ample in the geometric sense, Northcott's theorem does not hold
in general.

Let $k = \QQ(\sqrt{29})$, $\epsilon = (5 + \sqrt{29})/2$, and
$O_k = \ZZ[\epsilon]$. We set
\[
E = \Proj\left( O_k[X, Y, Z]/(Y^2Z + XYZ + \epsilon^2YZ^2 - X^3) \right).
\] 
Then, $E$ is an abelian scheme over $O_k$. Then, as in the proof of
\cite[Proposition~3.1.1]{MoArht},
we can construct a nef $C^{\infty}$-hermitian line bundle $\overline{H}$ on $E$
such that $[2]^*(\overline{H}) = \overline{H}^{\otimes 4}$ and
$H_k$ is ample on $E_k$, $c_1(\overline{H})$ is positive on $E(\CC)$, and that
$\adeg\left(\acherncl_1(\overline{H})^2\right) = 0$.
Let $K$ be the function field of $E$. Then,
$\overline{B} = (E; \overline{H})$ is a polarization of $K$.
Here we claim that Northcott's theorem dose not hold 
for the polarization $(E, \overline{H})$ of $K$.

Let $p_i : E \times_{O_k} E \to E$ be the projection to the $i$-th factor.
Then, considering $p_2 : E \times_{O_k} E \to E$,
$(E \times_{O_k} E, p_1^*(\overline{H}))$ gives rise to
a model of $(E_K, H_K)$.
Let $\Gamma_n$ be the graph of $[2]^n : E \to E$, i.e.,
$\Gamma_n = \{ ([2]^n(x), x) \mid x \in E \}$.
Moreover, let $x_n$ be a $K$-valued point of $E_K$ arising from $\Gamma_n$.
Then, if we denote the section $E \to \Gamma_n$ by $s_n$, then
\begin{align*}
h_{H_K}^{\overline{B}}(x_n) & = \adeg 
\left( p_1^*(\overline{H}) \cdot p_2^*(\overline{H}) \cdot \Gamma_n \right)
= \adeg \left( s_n^*( p_1^*(\overline{H})) \cdot 
s_n^*(p_2^*(\overline{H})) \right) \\
& = \adeg \left( ([2]^n)^*(\overline{H}) \cdot \overline{H} \right)
= \adeg \left( \overline{H}^{\otimes 4^n} \cdot \overline{H} \right)
= 4^n \adeg \left(\overline{H} \cdot \overline{H} \right) = 0.
\end{align*}
On the other hand, $x_n$'s are distinct points in $E_K(K)$.
\end{Remark}

\section{Geometric Northcott's theorem}

\begin{Proposition}
\label{prop:geom:northcott:thm}
Let $X$ be a smooth projective variety over an algebraically closed field $k$
of characteristic zero,
$C$ a smooth projective curve over $k$, and $f : X \to C$
a surjective morphism whose generic fiber is geometrically irreducible.
Let $L$ be an ample line bundle on $X$.
If $\deg(f_*(\omega_{X/C}^{n})) >0$ for some $n > 0$,
then, for any number $A$, the set
\[
\left\{ \Delta \mid \hbox{$\Delta$ is a section of $f : X \to C$ with
$(L \cdot \Delta) \leq A$} 
\right\}
\]
is not dense in $X$.
\end{Proposition}

\Proof
Let us begin with the following lemma.

\begin{Lemma}
\label{lem:direct:image:dualizing:split}
Let $f : X \to Y$ be a surjective morphism
of smooth projective varieties over an algebraically closed field $k$ 
of characteristic zero.
If there are a projective smooth algebraic variety $T$ over $k$ 
and a dominant rational
map $\phi : T \times_k Y \dasharrow X$ over $Y$, then 
the double dual $f_*(\omega_{X/Y}^n)^{\vee\vee}$ of $f_*(\omega_{X/Y}^n)$
is a free $\OO_Y$-sheaf for all $n \geq 0$.
\end{Lemma}

\Proof
Let $A$ be a very ample line bundle on $T$.
If $\dim T > \dim f$ and $T_1$ is a general
member of $|A|$, then 
$\rest{\phi}{T_1 \times Y} : T_1 \times Y \dasharrow X$ still
dominates $X$. Thus, considering induction on $\dim T$,
we may assume that $\dim T = \dim f$.

Let $\mu : Z \to T \times Y$ be a birational morphism of
smooth projective varieties such that 
$\psi = \phi \cdot \mu : Z \to X$ is a morphism.
Then, $\psi$ is generically finite. Thus, there is
a natural injection 
$\psi^*(\omega_{X/Y}) \hookrightarrow \omega_{Z/Y}$.
Hence, $\psi^*(\omega_{X/Y}^n) \hookrightarrow \omega_{Z/Y}^n$
for all $n > 0$. Therefore,
\[
\omega_{X/Y}^n \hookrightarrow \psi_*(\psi^*(\omega_{X/Y}^n))
\hookrightarrow \psi_*(\omega_{Z/Y}^n).
\]
Applying $f_*$ to the above injection, we have
\[
f_*(\omega_{X/Y}^n) \hookrightarrow f_*(\psi_*(\omega_{Z/Y}^n)).
\]
Further, letting $p$ be the natural projection $p : T \times Y \to Y$,
\[
f_*(\psi_*(\omega_{Z/Y}^n)) = p_*(\mu_*(\omega_{Z/Y}^n))
= p_*(\omega_{T \times Y/Y}^n) = H^0(T, \omega_T^n) \otimes_k \OO_Y.
\]
Thus, $f_*(\omega_{X/Y}^n)^{\vee\vee}$ is a subsheaf of the free sheaf
$H^0(T, \omega_T^n) \otimes_k \OO_Y$.

Here we claim
\addtocounter{Claim}{1}
\begin{equation}
\label{eqn:lem:direct:image:dualizing:split:1}
\left(c_1\left(f_*(\omega_{X/Y}^n)^{\vee\vee}\right) 
\cdot H^{d-1}\right) \geq 0,
\end{equation}
where $H$ is an ample line bundle on $Y$ and $d = \dim Y$.
This is an immediate consequence of weak positivity of 
$f_*(\omega_{X/Y}^n)^{\vee\vee}$ due to Viehweg \cite{Vi}.
We can however conclude our claim by a weaker result of Kawamata \cite{Ka1}, 
namely $\deg(f_*(\omega_{X/Y}^n)) \geq 0$ if $\dim Y = 1$.
For, considering complete intersections by general members of $|H^m|$ ($m \gg 0$),
we may assume $\dim Y = 1$. 

We can find a projection $\alpha : H^0(T, \omega_T^n) \otimes_k \OO_Y
\to \OO_Y^{\oplus r_n}$ such that 
$r_n = \rank f_*(\omega_{X/Y}^n)^{\vee\vee}$
and the composition
\[
f_*(\omega_{X/Y}^n)^{\vee\vee} \hookrightarrow
H^0(T, \omega_T^n) \otimes_k \OO_Y \overset{\alpha}{\longrightarrow} \OO_Y^{\oplus r_n}
\]
is injective. Therefore, since $f_*(\omega_{X/Y}^n)^{\vee\vee}$
is reflexive, the above homomorphism is an isomorphism by 
\eqref{eqn:lem:direct:image:dualizing:split:1}.
\QED

Let us go back to the proof of Proposition~\ref{prop:geom:northcott:thm}.
Let $\operatorname{Hom}_k(C, X)$ 
be a scheme consisting of morphisms from $C$ to $X$.
Then, there is a morphism $\alpha : \operatorname{Hom}_k(C, X) \to 
\operatorname{Hom}_k(C, C)$
given by $\alpha(s) = f \cdot s$.
We set $\operatorname{Sec}(f) = \alpha^{-1}(\operatorname{id}_C)$.
Then, there is a natural morphism
$\beta : \operatorname{Sec}(f) \times C \to X$ given by
$\beta(s,y) = s(y)$.
Since $L$ is ample,
\[
\left\{ \Delta \mid \text{$\Delta$ is a section of $f : X \to C$ with
$(L \cdot \Delta) \leq A$} 
\right\}
\]
is a bounded family, so that
there are finitely many connected components
$\operatorname{Sec}(f)_1, \ldots, \operatorname{Sec}(f)_r$ 
of $\operatorname{Sec}(f)$ such that,
for all sections $\Delta$ with $(L \cdot \Delta) \leq A$,
there is $s \in \operatorname{Sec}(f)_i$ for some $i$ with
$\Delta = s(C)$.
On the other hand, by Lemma~\ref{lem:direct:image:dualizing:split},
$\operatorname{Sec}(f)_i \times C \to X$
is not a dominant morphism for every $i$. Thus, we get our proposition.
\QED

\end{document}